\documentclass[12pt]{amsart}
\usepackage{a4wide,amsfonts,graphicx,
amssymb,stmaryrd,amscd,amsmath,latexsym,amsbsy}
\numberwithin{equation}{section}
\theoremstyle{plain}
\newtheorem{thm}{Theorem}[section]

\newtheorem{prop}[thm]{Proposition}

\newtheorem{defi}[thm]{Definition}
\newtheorem{lem}[thm]{Lemma}
\newtheorem{cor}[thm]{Corollary}
\newtheorem{eg}[thm]{Example}

\theoremstyle{remark}
\newtheorem{rema}[thm]{Remark}

\title{Macdonald-Koornwinder polynomials}
\author{J.V. Stokman}
\address{J.V. Stokman, Korteweg-de Vries 
Institute for Mathematics, Universiteit van Amsterdam,
Science Park 904, 1098 XH Amsterdam, The Netherlands.}
\email{j.v.stokman@uva.nl}
\begin{document}
\begin{abstract}
An overview of the basic results on
Macdonald-Koornwinder polynomials
and double affine Hecke algebras is given. 
We develop the theory in such a way that it naturally encompasses
all known cases as well as a new rank two case. 
Among the basic properties of the Macdonald-Koornwinder polynomials we
treat are the quadratic norm formulas, duality and the evaluation formulas.
This text is a provisional version of a chapter on Macdonald-Koornwinder
polynomials for volume 5 of the Askey-Bateman
project, entitled ``Multivariable special functions''. 
\end{abstract}

\maketitle
\tableofcontents
\section{Introduction}
This text is a provisional version of a chapter on Macdonald-Koornwinder 
polynomials for volume 5 of the Askey-Bateman
project, entitled ``Multivariable special functions''. 
The aim is to introduce  
nonsymmetric and symmetric Macdonald-Koornwinder 
polynomials and to present their 
basic properties: (bi-)or\-tho\-go\-na\-li\-ty, norm formulas,
q-difference(-reflection) equations, duality and evaluation formulas.

Symmetric $\textup{GL}$ type Macdonald 
polynomials were introduced by Macdonald \cite{Mbook}
as a two-parameter family of orthogonal multivariate
polynomials interpolating between the Jack polynomials 
and the Hall-Littlewood polynomials. Root system generalizations
were subsequently defined by Macdonald in \cite{Mpol}. They were labelled by 
so called admissible pairs of root systems. Recasting these
data in terms of affine root systems (cf. \cite{C,M}) it 
is natural to speak of an untwisted theory and a twisted theory of 
Macdonald polynomials associated to root systems.
An important further extension for nonreduced root systems was
constructed by Koornwinder \cite{Ko}. In \cite{M,Ha} important steps
were undertaken to develop a general theory covering all these cases.
In this text we pursue this further. We call the corresponding polynomials
Macdonald-Koornwinder polynomials.

The symmetric Macdonald-Koornwinder polynomials are root system generalizations 
of classical one-variable $q$-orthogonal polynomials 
from the $q$-Askey scheme \cite{KS}.
Before giving the contents of this chapter in more detail,
we first illustrate this point of view for the symmetric $\textup{GL}$
Macdonald polynomials and the symmetric Koornwinder polynomials. 
\subsection{Symmetric $\textup{GL}$ type Macdonald polynomials}
The symmetric $\textup{GL}_{n+1}$ type Macdonald \cite{Mbook}
polynomials form a linear
basis of the space of symmetric Laurent polynomials in $n+1$ variables 
$t_1,\ldots,t_{n+1}$, depending on two parameters $0<q,\kappa<1$. They form
an orthogonal basis of the Hilbert space $L^2(T_u,v_+(t)dt)$,
where $dt$ is the normalized Haar measure on
$T_u=\{t=(t_1,\ldots,t_{n+1})\in\mathbb{C}^{n+1}\,\, | \,\, |t_i|=1\}$
and with weight function $v_+(t)$ explicitly given by
\[
v_+(t)=\prod_{1\leq i\not=j\leq n+1}\frac{\bigl(t_i/t_j;q\bigr)_{\infty}}
{\bigl(\kappa^2t_i/t_j;q\bigr)_{\infty}}
\]
(the $q$-shifted factorial is defined by \eqref{qshifted}). In addition,
the symmetric $\textup{GL}_{n+1}$ type Macdonald polynomials are common
eigenfunctions of the commuting trigonometric Ruijsenaars-Macdonald 
\cite{R,Mbook} $q$-difference operators
\[
(D_jf)(t):=\sum_{\stackrel{I\subseteq \{1,\ldots,n+1\}}
{\#I=j}}\left(\prod_{r\in I,s\not\in I}\frac{\kappa^{-1}t_r-\kappa t_s}
{t_r-t_s}\right)f(q^{-\sum_{r\in I}\epsilon_r}t)
\qquad 1\leq j\leq n+1,
\]
where $q^{-\sum_{r\in I}\epsilon_r}$ is the $(n+1)$-vector with 
$q^{-1}$'s at the entries labelled by
$I$ and ones elsewhere. The eigenvalue equation
for $D_{n+1}$ is equivalent to a homogeneity condition for the symmetric
Macdonald polynomials. Up to a monomial factor the symmetric $\textup{GL}_{n+1}$
Macdonald polynomials thus only depend on the
$n$ variables $t_1/t_2,\ldots,t_n/t_{n+1}$. For $n=1$ the Macdonald
polynomials are essentially the continuous $q$-ultraspherical polynomials
in the variable $t_1/t_2$.
The above results then reduce to 
the orthogonality relations and the
second order $q$-difference equation satisfied by the continuous 
$q$-ultraspherical polynomials, see Subsection \ref{GLn} for further details.

\subsection{Symmetric Koornwinder polynomials}
The symmetric Koornwinder \cite{Ko} polynomials form a linear basis
of the space of Laurent polynomials in $n$ variables $t_1,\ldots,t_n$
which are invariant under the hyperoctahedral group 
(acting by permutations and inversions of the variables), depending
on six parameters $0<q,a,b,c,d,k<1$. They form an orthogonal basis of the
Hilbert space $L^2(T_u,v_+(t)dt)$ with $T_u=\{t\in\mathbb{C}^n \, | \, 
|t_i|=1\}$, Haar measure $dt$ on $T_u$, and weight function
\[
v_+(t)=\prod_{i=1}^n\frac{\bigl(t_i^{\pm 2};q\bigr)_{\infty}}
{\bigl(at_i^{\pm 1};q\bigr)_{\infty}
\bigl(bt_i^{\pm 1};q\bigr)_{\infty}
\bigl(ct_i^{\pm 1};q\bigr)_{\infty}
\bigl(dt_i^{\pm 1};q\bigr)_{\infty}}
\prod_{1\leq r\not=s\leq n}\frac{\bigl(t_tt_s^{\pm 1};q\bigr)_{\infty}}
{\bigl(kt_rt_s^{\pm 1};q\bigr)_{\infty}},
\]
where $(uz^{\pm 1};q\bigr)_{\infty}:=(uz;q\bigr)_{\infty}(uz^{-1};q\bigr)_{\infty}$.
The symmetric Koornwinder polynomials
are eigenfunctions of Koornwinder's \cite{Ko} multivariable
extension of the Askey-Wilson \cite{AW} second order $q$-difference operator
\begin{equation*}
\begin{split}
(Df)(t):=&\sum_{i=1}^n\sum_{\xi\in\{\pm 1\}}A_i^\xi(t)
\bigl(f(q^{\xi\epsilon_i}t)-f(t)\bigr),\\
A_i^\xi(t):=&\frac{(1-at_i^{\xi})(1-bt_i^{\xi})(1-ct_i^{\xi})
(1-dt_i^{\xi})}{(1-t_i^{2\xi})(1-qt_i^{2\xi})}
\prod_{j\not=i}\frac{(1-kt_i^\xi t_j)(1-kt_i^\xi t_j^{-1})}
{(1-t_i^\xi t_j)(1-t_i^\xi t_j^{-1})}.
\end{split}
\end{equation*}
For $n=1$ the $k$-dependence drops out and the symmetric Koornwinder polynomials
are the Askey-Wilson polynomials from \cite{AW}. See Subsection 
\ref{CcheckC} for further details.

\subsection{General description of the contents}

Most literature on Macdonald-Koornwinder
polynomials deals with one out of the above mentioned 
four cases of Macdonald-Koornwinder type polynomials (the $\textup{GL}$-case, 
the untwisted case, the twisted case and the Koornwinder case).
Cherednik \cite{C} treats the first three cases separately.
Macdonald's \cite{M} exposition covers
the last three cases, but various steps still need case by case analysis.
We develop the theory 
in such a way that it naturally unifies the above four
cases, but in addition contains a new class of rank
two Macdonald-Koornwinder type polynomials (see Subsection \ref{Except2}). 
The setup will be close to
Haiman's \cite{Ha} approach, which allowed him to give a uniform
proof of the duality of the Macdonald-Koornwinder polynomials.

It is tempting to believe that
the initial data for the Macdonald-Koornwinder polynomials
should be similarity classes of 
irreducible affine root systems $R$ together with a
choice of a deformation parameter $q$ and a 
multiplicity function (playing the role of
the free parameters in the theory). In such a parametrization
the untwisted and twisted cases
should relate to the similarity classes of the irreducible reduced
affine root systems of untwisted and twisted type respectively, 
the $\textup{GL}$ case to the irreducible reduced affine root system of
type $\textup{A}$ with a ``reductive'' extension of the affine
Weyl group, and the Koornwinder case with the nonreduced irreducible
affine root system of type $C^\vee C$ (we refer here to the
classification of affine root systems from \cite{M0}, see also
the appendix). It turns out though that a more subtle labelling is needed
to capture all fundamental properties of the Macdonald-Koornwinder polynomials.

We take as initial data quintuples
$D=(R_0,\Delta_0,\bullet,\Lambda,\Lambda^d)$ with $R_0$ a finite reduced
irreducible root system, $\Delta_0$ an ordered basis of $R_0$, $\bullet\in
\{u,t\}$ ($u$ stands for untwisted and $t$ stands for twisted), 
and $\Lambda$, $\Lambda^d$ two lattices satisfying appropriate compatibility
conditions with respect to the (co)root lattice of $R_0$
(see \eqref{latticeconditions1}).
We build from $D$ an irreducible affine root system $R$ and an extended
affine Weyl group $W$. 
The extended affine
Weyl group $W$ is simply the semi-direct product group
$W_0\ltimes\Lambda^d$ with $W_0$ the Weyl group of $R_0$.
The affine root system $R$ is constructed as follows. We associate to
$R_0$ and $\bullet$ the reduced irreducible affine root system $R^\bullet$ 
of type $\bullet$ with gradient root system $R_0$.
Then $R$ is an irreducible affine root system obtained from $R^\bullet$
by adding $2a$ if $a\in R^\bullet$ has the property that the pairings of
the associated coroot $a^\vee$ to elements of $\Lambda$ take value
in $2\mathbb{Z}$ (see \eqref{S}). The Macdonald-Koornwinder 
polynomials associated to
$D$ are the ones naturally related to $R$ in the labelling
proposed in the previous paragraph.  

In the literature 
the terminology Macdonald polynomials is used when 
the underlying affine root system $R$ is reduced, while 
the terminology Koornwinder polynomials or Macdonald-Koornwinder 
polynomials is used when $R$ is nonreduced
and of type $C^\vee C$. In this text 
we will use the terminology Macdonald-Koornwinder
polynomials when dealing with arbitrary initial data $D$.
We will speak of Macdonald polynomials if the underlying affine root
system $R$ is reduced and of Koornwinder polynomials if $R$ is of type
$C^\vee C$.

Duality is related to a simple involution on
initial data, $D\mapsto D^d=(R_0^d,\Delta_0^d,\bullet,\Lambda^d,\Lambda)$
with $R_0^d$ the coroot system $R_0^\vee$ if $\bullet=u$ and the root
system $R_0$ if $\bullet=t$. This duality is subtle on the level of
affine root systems (it can for instance happen that the affine root
system $R^d$ associated to $D^d$ is reduced
while $R$ is nonreduced). 
The reason that the present choice of initial data is convenient
is the fact that the duality map $D\mapsto D^d$ on initial data
naturally lifts to the duality
antiisomorphism 
of the associated double affine braid group (see \cite{Ha} and
Section \ref{4}).
This duality antiisomorphism 
is the key tool to prove duality, evaluation formulas
and norm formulas for the Macdonald-Koornwinder polynomials.

\subsection{Topics that are not discussed}
We do not discuss the shift operators for the Macdonald-Koornwinder polynomials
(see, e.g., \cite{CAnn,St}), leading
to the explicit evaluations of the constant terms (generalized
Selberg integrals). Other
important developments
involving Macdonald-Koornwinder 
polynomials that are not discussed in this chapter, are
\begin{enumerate}
\item[{\bf --}] connections to combinatorics. This is discussed in 
a separate chapter of the fifth volume of the Askey-Bateman project,
\item[{\bf --}] connections to algebraic geometry, see, e.g.,
\cite{Hno,Hnot,SV},
\item[{\bf --}] connections to representation theory,
see \cite{C,Cr1,Cr2,Ka1,Ion1,Ion2} and references therein.
\item[{\bf --}] applications to harmonic analysis 
on quantum groups, see, e.g., 
\cite{No1,NoS,NDS,L,EK1,EK2,EK3,OSt},
\item[{\bf --}] 
limit cases of symmetric Macdonald-Koornwinder polynomials, see, e.g., 
\cite{Mbook,Mpol,C,HO,HS,CW,Stlim,BF} and the separate chapter of
the fifth volume of the Askey-Bateman project on Heckman-Opdam polynomials,
\item[{\bf --}] 
the theory of Gaussians, Macdonald-Mehta type integrals, and 
basic hypergeometric functions associated to root systems, see, e.g.,
\cite{CMM,StSM,CW,StPr},
\item[{\bf --}] interpolation Macdonald-Koornwinder polynomials,
see, e.g.,
\cite{Kn,KnS,Sah,O,O0,LRW},
\item[{\bf --}]
the relation to
quantum integrable systems, 
see, e.g., \cite{CKZ1,CKZ2,R,SVe,Ch,ChE,P,KP,KT,dFZJ,StKZ}, 
\item[{\bf --}]
special parameter values (e.g. roots of unity), see, e.g., 
\cite{C,Cns,vDSt,KMSV,Ka1,Ka2,FJMM,Ch,ChE},
\item[{\bf --}] affine and  
elliptic generalizations, see, e.g.,
\cite{EK4,Caff,R} and \cite{KH,Ra1,Ra2,CG} respectively.
\end{enumerate}
See the introductory chapter of Cherednik's \cite{C} book for 
an extensive discussion on the ramifications of Macdonald polynomials.

\subsection{Detailed description of the contents}
Precise references to the literature are given in the main text.

In Section \ref{2} we give the definition of the
affine braid group, affine Weyl group and affine Hecke algebra. 
We determine an explicit realization of the affine Hecke algebra, 
which will serve 
as starting point of the Cherednik-Macdonald theory on Macdonald-Koornwinder
polynomials in the next section. 
In addition we introduce the initial data. We introduce the
space of multiplicity functions associated to the fixed initial data $D$. 
We extend the duality on initial data
to an isomorphism of the associated spaces of multiplicity functions.
We give the basic representation of the extended affine Hecke algebra
associated to $D$, using the explicit realization of the affine Hecke
algebra. 

In Section \ref{3} we define and study the nonsymmetric and symmetric
monic Macdonald-Koornwinder polynomials associated to the initial data $D$. 
We first focus on the nonsymmetric monic Macdonald-Koornwinder polynomials.
We characterize them as common eigenfunctions of a family
of commuting $q$-difference
reflection operators. These operators are obtained as the
images under the basic representation of the elements of the
Bernstein-Zelevinsky abelian subalgebra of the extended affine Hecke algebra.
We determine the biorthogonality relations of the nonsymmetric monic 
Macdonald-Koornwinder
polynomials and use the finite Hecke symmetrizer to obtain the
symmetric monic Macdonald-Koornwinder polynomials.
We give the orthogonality relations of the symmetric monic
Macdonald-Koornwinder polynomials and show that they are common eigenfunctions
of the commuting Macdonald $q$-difference operators (also known
as Ruijsenaars operators in the $\textup{GL}$ case). 
We finish this section
by describing three
cases in detail: the $\textup{GL}$ case, the $C^\vee C$ case, and an 
exceptional
nonreduced rank two case not covered in Cherednik's \cite{C}
and Macdonald's \cite{M} treatments.

In Section \ref{4} we introduce the double affine braid group 
and the double affine Hecke algebra associated
to $D$ and $(D,\kappa)$ respectively, where $\kappa$ is a choice of
a multiplicity function on $R$. We lift the 
duality on initial data to a duality antiisomorphism on the level of the
associated double affine braid groups. We show how it descends to the
level of double affine Hecke algebras 
and how it leads to
an explicit evaluation formula for the monic Macdonald-Koornwinder polynomials.
We proceed by
defining the associated normalized nonsymmetric and symmetric 
Macdonald-Koornwinder 
polynomials and deriving their duality and quadratic norms.

In Section \ref{5} we give the norm and evaluation formulas
in terms of $q$-shifted factorials 
for the $\textup{GL}$ case, the $C^\vee C$ case, and for 
the exceptional nonreduced rank two case. 

In the appendix we give a short introduction
to (the classification of) affine root systems, following 
closely \cite{M0} but with some adjustments. We close the appendix
with a list of all the affine Dynkin diagrams.

\section{The basic representation of the extended affine Hecke 
algebra}\label{2}
We first introduce the affine Hecke algebra and
the appropriate initial data for the Cherednik-Macdonald theory on 
Macdonald-Koornwinder polynomials. Then we introduce 
the basic representation of the extended affine Hecke algebra, 
which is fundamental in the development of the Cherednik-Macdonald
theory.

\subsection{Affine Hecke algebras}
A convenient reference for this subsection is \cite{Hu2}.
For unexplained notations and terminology regarding affine Weyl groups
we refer to the appendix.

For a generalized Cartan matrix $A=(a_{ij})_{1\leq i,j\leq r}$
let $M=(m_{ij})_{1\leq i,j\leq r}$ be the matrix with 
$m_{ii}=1$ and for $i\not=j$, $m_{ij}=2,3,4,6,\infty$ if $a_{ij}a_{ji}=0,1,2,3,
\geq 4$
respectively.
\begin{defi}
Let $A=(a_{ij})_{1\leq i,j\leq r}$ be a generalized Cartan matrix.
\begin{enumerate}
\item[(1)] The braid group $\mathcal{B}(A)$ is the group generated by
$T_i$ ($1\leq i\leq r$) with defining relations $T_iT_jT_i\cdots=
T_jT_iT_j\cdots$ ($m_{ij}$ factors on each side)
if $1\leq i\not=j\leq r$ (which should be interpreted as no relation
if $m_{ij}=\infty$).
\item[(2)] The Coxeter group $\mathcal{W}(A)$ associated to $A$ is 
the quotient of $\mathcal{B}(A)$ by the normal subgroup generated
by $T_i^2$ ($1\leq i\leq r$).
\end{enumerate}
\end{defi}
It is convenient to denote by $s_i$ the element
in $\mathcal{W}(A)$ corresponding to $T_i$ for $1\leq i\leq r$.
They are the Coxeter generators of the Coxeter group $\mathcal{W}(A)$.

Let $A=(a_{ij})_{1\leq i,j\leq r}$ be a generalized Cartan matrix.
Suppose that $k_i$ ($1\leq i\leq r$) are nonzero complex numbers
such that $k_i=k_j$
if $s_i$ is conjugate to $s_j$ in $\mathcal{W}(A)$. We write $k$ for the
collection $\{k_i\}_i$. Let 
$\mathbb{C}[\mathcal{B}(A)]$ be the complex group algebra of the braid group
$\mathcal{B}(A)$.
\begin{defi}
The Hecke algebra $H(A,\kappa)$ is the complex unital
associative algebra given by $\mathbb{C}[\mathcal{B}(A)]/I_k$, where
$I_k$ is the two-sided ideal of $\mathbb{C}[\mathcal{B}(A)]$
generated by $(T_i-k_i)(T_i+k_i^{-1})$ for $1\leq i\leq r$.
\end{defi}
If $k_i=1$ for all $1\leq i\leq r$ then the associated affine
Hecke algebra is the complex group
algebra of $\mathcal{W}(A)$.

If $w=s_{i_1}s_{i_2}\cdots s_{i_r}$ is a reduced expression in $\mathcal{W}(A)$,
i.e. a shortest expression of $w$ as product of Coxeter generators,
then 
\[
T_w:=T_{i_1}T_{i_2}\cdots T_{i_r}\in H(A,k)
\]
is well defined (already in the braid group $\mathcal{B}(A)$), 
and the $T_w$ ($w\in\mathcal{W}(A)$) form a complex linear
basis of $H(A,k)$.

Suppose $R_0\subset V$ is a finite crystallographic 
root system with ordered basis $\Delta_0=(\alpha_1,\ldots,
\alpha_n)$ and write $A_0$ for the associated Cartan matrix.
Then $W_0\simeq\mathcal{W}(A_0)$ by mapping the simple reflections
$s_{\alpha_i}\in W_0$ to the Coxeter generators $s_i$ of
$\mathcal{W}(A_0)$ for $1\leq i\leq n$. 
The associated finite Hecke algebra $H(A_0,k)$ depends only on 
$k$ and $W_0$ (as Coxeter group). 
We will sometimes denote it by $H(W_0,k)$.

Similarly, if $R$ is an irreducible affine root system with
ordered basis $\Delta=(a_0,\ldots,a_n)$ 
and if $A$ is the associated affine Cartan
matrix, then 
$\mathcal{W}(A)\simeq W(R)$ by $s_i\mapsto s_{a_i}$ ($0\leq i\leq n)$.
Again we write $H(W(R),k)$ for the associated affine Hecke algebra
$H(A,k)$.

\subsection{Realizations of the affine Hecke algebra}\label{realizations}
We use the notations on affine root systems as introduced in the appendix.
The construction is motivated by Cherednik's polynomial representation
\cite[Thm. 3.2.1]{C} 
of the affine Hecke algebra and its extension to the nonreduced case
by Noumi \cite{N}.

Let $R\subset\widehat{E}$ 
be an irreducible affine root system on the affine Euclidean
space $E$ (possibly nonreduced) with affine Weyl group $W$,
and fix an ordered basis $\Delta=(a_0,a_1,\ldots,a_n)$ of $R$. 
Write $A=A(R,\Delta)$ for the associated affine Cartan matrix. 

Consider the lattice $\mathbb{Z}R$ in $\widehat{E}$. It is a full
$W$-stable lattice with $\mathbb{Z}$-basis the simple affine roots. 
Denote by
$F$ the quotient field $\textup{Quot}\bigl(\mathbb{C}[\mathbb{Z}R]\bigr)$ 
of the complex group algebra $\mathbb{C}[\mathbb{Z}R]$ of $\mathbb{Z}R$.
It is convenient to write $e^\lambda$ ($\lambda\in\mathbb{Z}R$) for the natural
complex linear basis of $\mathbb{C}[\mathbb{Z}R]$. 
The multiplicative structure 
of $F$ is determined by
\[
e^0=1,\quad e^{\lambda+\mu}=e^\lambda e^\mu,\qquad \lambda,\mu\in\mathbb{Z}R.
\]

The affine Weyl group $W$ canonically acts by field automorphisms on $F$.
On the basis elements $e^\lambda$ the $W$-action reads 
$w(e^\lambda)=e^{w\lambda}$ ($w\in W$, $\lambda\in\mathbb{Z}R$).
Since $W$ acts by algebra automorphisms on $F$, we can form the
semidirect product algebra $W\ltimes F$.

Let $k: R\rightarrow \mathbb{C}^*:=\mathbb{C}\setminus\{0\}$, 
$a\mapsto k_a$ be a $W$-equivariant map,
i.e. $k_{wa}=k_a$ for all $w\in W$ and $a\in R$.
If $a\in R$ but $2a\not\in R$ then we set $k_{2a}:=k_a$.
Note that $k^{ind}:=k|_{R^{ind}}$ is a $W$-equivariant map on $R^{ind}$,
which is determined by its values $k_i:=k_{a_i}$
($0\leq i\leq n$) on the simple affine roots. 
It satisfies $k_i=k_j$ if 
$s_i$ is conjugate to $s_j$ in $W$.
Hence we can form the associated affine Hecke algebra $H(W(R),k^{ind})$.
\begin{thm}\label{beta}
With the above notations and conventions, 
there exists a unique algebra monomorphism
\[
\beta=\beta_{R,\Delta,k}: H(W(R),k^{ind})\hookrightarrow W\ltimes F
\]
satisfying
\begin{equation}\label{betaTi}
\beta(T_i)=k_{i}s_i+\frac{k_{i}-k_{i}^{-1}+
(k_{2a_i}-k_{2a_i}^{-1})e^{a_i}}{1-e^{2a_i}}(1-s_i)
\end{equation}
for $0\leq i\leq n$.
\end{thm}
The proof uses the Bernstein-Zelevinsky presentation of the
affine Hecke algebra, which we present in a slightly more general
context in Subsection \ref{Bernstein}.
\begin{rema}
If $2a_i\not\in R$ then \eqref{betaTi} simplifies to
\[
\beta(T_i)=k_{i}s_i+\frac{k_{i}-k_{i}^{-1}}{1-e^{a_i}}(1-s_i).
\]
\end{rema}

The notion of similarity of pairs $(R,\Delta)$ (see the appendix)
can be extended to triples $(R,\Delta,k)$ in the obvious way. 
The algebra homomorphisms
$\beta$ that are associated to different representatives of the 
similarity class of $(R,\Delta,k)$ 
are then equivalent in a natural sense. 
Starting from the next subsection we therefore 
will focus on the explicit representatives
of the similarity classes as described in the appendix in
Subsection \ref{special}.

Recall from the appendix that the classification of irreducible affine root
systems leads to a subdivision of irreducible affine root systems in 
three types, namely untwisted, twisted and mixed type.
It is easy to show that the algebra map $\beta_{R,\Delta,k}$ 
in case $R$ is of mixed type (possibly nonreduced) can alternatively
be written as $\beta_{(R^\prime,\Delta^\prime,k^\prime)}$ with appropriately
chosen triple $(R^\prime,\Delta^\prime,k^\prime)$ and with $R^\prime$
of untwisted or of twisted type. In the Cherednik-Macdonald
theory, the mixed type can therefore safely be ignored, as we in fact
will do.

\subsection{Initial data}
It turns out that labelling Macdonald-Koornwinder polynomials by
irreducible affine root systems is not the most convenient way to
proceed. In this subsection we give more convenient initial data
and we explain how it relates 
to affine root systems. For basic notations and facts on affine root systems
we refer to the appendix. 
\begin{defi}
The set $\mathcal{D}$ of initial data consists of quintuples
\[
D=(R_0,\Delta_0,\bullet,\Lambda,\Lambda^d)
\]
with
\begin{enumerate}
\item $R_0$ is 
a finite set of nonzero vectors in an Euclidean space $Z$
forming a finite, irreducible, reduced crystallographic 
root system within
the real span $V$ of $R_0$,
\item $\Delta_0=(\alpha_1,\ldots,\alpha_n)$ is an ordered basis of $R_0$,
\item $\bullet=u$ or $\bullet=t$,
\item $\Lambda$ and $\Lambda^d$ are full lattices in $Z$,
satisfying 
\begin{equation}\label{latticeconditions1}
\begin{split}
&\mathbb{Z}R_0\subseteq\Lambda,\qquad 
\bigl(\Lambda,\mathbb{Z}R_0^\vee\bigr)\subseteq\mathbb{Z},\\
&\mathbb{Z}R_0^{d}\subseteq\Lambda^d,\qquad
\bigl(\Lambda^d,\mathbb{Z}R_0^{d}{}^{\vee}\bigr)\subseteq\mathbb{Z}
\end{split}
\end{equation}
where $R_0^d=\{\alpha^d:=\mu_\alpha^\bullet\alpha^\vee\}_{\alpha\in R_0}$ and
\begin{equation*}
\begin{split}
\mu^u_\alpha&:=1,\qquad\qquad \alpha\in R_0,\\
\mu^t_\alpha&:=|\alpha|^2/2,\qquad \alpha\in R_0.
\end{split}
\end{equation*}
\end{enumerate}
\end{defi}
Note that $R_0^d=R_0^\vee$ if $\bullet=u$ and $=R_0$
if $\bullet=t$.

We view the vector space $\widehat{V}$ of real valued 
affine linear functions on $V$ as the subspace of $\widehat{Z}$ consisting of
affine linear functions on $Z$ which are constant on the 
orthocomplement $V^\perp$ of $V$ in $Z$.
We write $c$ for the constant function one on $V$ as well as on $Z$.
In a similar fashion we view the orthogonal group
$O(V)$ as subgroup of $O(Z)$ and 
$O_c(\widehat{V})$ as subgroup of $O_c(\widehat{Z})$, where $O_c(\widehat{V})$
is the subgroup of linear automorphisms of $\widehat{V}$ preserving $c$
and preserving the natural semi-positive definite form on $\widehat{V}$
(see the appendix for further details).

Fix $D=(R_0,\Delta_0,\bullet,\Lambda,\Lambda^d)
\in\mathcal{D}$. We associate to $D$ a triple $(R,\Delta,W)$
of an affine root system $R=R(D)$, an ordered basis $\Delta=\Delta(D)$
of $R$ and an extended affine Weyl group $W=W(D)$ as follows.
We first define a reduced irreducible affine root system $R^\bullet\subset
\widehat{V}$ 
to $D$ as
\begin{equation}\label{affinerootsystem}
R^\bullet=\{m\mu_\alpha^{\bullet}c+
\alpha\}_{m\in\mathbb{Z}, \alpha\in R_0}
\end{equation}
for $\bullet\in\{u,t\}$. In other words,
$R^u:=\mathcal{S}(R_0)$ and $R^t:=\mathcal{S}(R_0^\vee)^\vee$ 
in the notations of the
appendix (see \eqref{special}). 
Let $\varphi\in R_0$ 
(respectively $\theta\in R_0$) be the highest root
(respectively highest short root) of $R_0$ with respect to
the ordered basis $\Delta_0$ of $R_0$.
Then the ordered basis $\Delta=\Delta(D)$ of $R^\bullet$ is set to be
\[
\Delta:=(a_0,a_1,\ldots,a_n)
\]
with $a_i:=\alpha_i$ for $1\leq i\leq n$ and 
\begin{equation}\label{simpleaffine}
a_0:=
\begin{cases} 
c-\varphi\quad\quad\, \hbox{ if }\,\, &\bullet=u,\\
\frac{|\theta|^2}{2}c-\theta\quad \hbox{ if }\,\, &\bullet=t.
\end{cases}
\end{equation}
\begin{rema}
Suppose that $(R^\prime,\Delta^\prime)$ is a pair consisting 
of a reduced irreducible affine root system
$R^\prime$ and an ordered basis $\Delta^\prime$ of $R^\prime$. 
If $R^\prime$ is similar to $R^\bullet$
then there exists a similarity transformation realizing 
$R^\prime\simeq R^\bullet$ 
and mapping $\Delta^\prime$ to $\Delta^\bullet$ as unordered sets.
\end{rema}
The affine root system $R$ is the following extension of 
$R^\bullet$. 
Define the subset $S=S(D)$ by
\begin{equation}\label{S}
S:=\{i\in\{0,\ldots,n\} \, | \, \bigl(\Lambda,a_i^\vee\bigr)=
2\mathbb{Z} \}.
\end{equation}
Let $W^\bullet$ be the affine Weyl group of $R^\bullet$.
Then we set
\begin{equation}\label{affinerootsystemnonred}
R=R(D):=R^\bullet\cup\bigcup_{i\in S}W^\bullet(2a_i),
\end{equation}
which is an irreducible affine root system since 
$\mathbb{Z}R_0\subseteq\Lambda$. Note that $\Delta$ is also an
ordered basis of $R$.

\begin{rema}
Note that $R$ is an irreducible affine root system of untwisted or twisted
type, but never of mixed type (see the appendix for the terminology).
But the irreducible affine root systems of mixed type are
affine root subsystems of the affine root system of type $C^\vee C$,
which is the nonreduced extension of
the affine root system $R^t$ with $R_0$ of type $\textup{B}$ (see
Subsection \ref{CcheckC} for a detailed description of the affine root
system of type $C^\vee C$).
Accordingly, special cases of the Koornwinder polynomials are
naturally attached to 
affine root systems of mixed type, 
see Subsection \ref{realizations} and Remark \ref{Bnext}.
\end{rema}
\begin{rema}
The nonreduced extension of the affine root system $R^u$ with
$R_0$ of type $\textup{B}_2$ is not an affine root subsystem of
the affine root system of type $C^\vee C_2$. 
It can actually be better viewed as the rank two
case of the family $R^u$ with $R_0$ of type $\textup{C}_n$
since, in the corresponding
affine Dynkin diagram (see Subsection \ref{Dynkin}), 
the vertex labelled by the affine simple
root $a_0$ is double bonded with the finite Dynkin diagram of $R_0$. 
The nonreduced extension of $R^u$ 
with $R_0$ of type $\textup{C}_2$ 
was missing in Macdonald's \cite{M0} classification list. It was
added in \cite[(1.3.17)]{M}, but the associated theory of Macdonald-Koornwinder
polynomials was not developed. In the present setup it is a special case
of the general theory. We will describe this particular 
case in detail in Subsection \ref{Except2}.
\end{rema}

Finally we define the extended
affine Weyl group $W=W(D)$ (please consult the appendix for the notations). 
Write $s_i:=s_{a_i}$ ($0\leq i\leq n$) for the simple reflections
of $W^\bullet$. Note that $s_i=s_{\alpha_i}$ ($1\leq i\leq n$) are the
simple reflections of the finite Weyl group $W_0$ of $R_0$. Furthermore,
$s_0=\tau(\varphi^\vee)s_{\varphi}$ if $\bullet=u$ and
$s_0=\tau(\theta)s_\theta$ if $\bullet=t$, where $\tau(v)$ stands for the
translation by $v$ (see the appendix). 
Consequently
\[
W^\bullet\simeq W_0\ltimes \tau(\mathbb{Z} R_0^{d}).
\]
We omit $\tau$ from the notations if no confusion can arise.
The extended affine Weyl group $W=W(D)$ is now defined as
\[
W:=W_0\ltimes\Lambda^d.
\] 
It contains the
affine Weyl group $W^\bullet$ of $R^\bullet$ as normal subgroup, and
$W/W^\bullet\simeq \Lambda^d/\mathbb{Z}R_0^{d}$. 

The affine root system $R^\bullet\subset \widehat{Z}$
is $W$-stable since
\begin{equation}\label{translationaction2}
\tau(\xi)(m\mu_{\beta}^\bullet c+\beta)=\bigl(m-\bigl(\xi,\beta^{d\vee}\bigr)
\bigr)\mu_{\beta}^\bullet c+\beta,\qquad m\in\mathbb{Z},\, \beta\in R_0
\end{equation}
and $(\xi,\beta^{d\vee})\in\mathbb{Z}$ for $\xi\in\Lambda^d$ and
$\beta\in R_0$. Moreover, the  
affine root system $R$ is $W$-invariant. 

We now proceed by giving key examples of initial data.
Recall that for a finite root system $R_0\subset V$, 
\[
P(R_0):=\{\lambda\in V \,\, | \,\, (\lambda,\alpha^\vee)\in\mathbb{Z}
\quad \forall\,\alpha\in R_0\}
\]
is the weight lattice of $R_0$. If $\Delta_0=(\alpha_1,\ldots,\alpha_n)$
is an ordered basis of $R_0$ then we write $\varpi_i\in P(R_0)$ 
($1\leq i\leq n$) for the corresponding fundamental weights, which
are characterized by $(\varpi_i,\alpha_j^\vee)=\delta_{i,j}$. 
\begin{eg}\label{lattices}
{\bf (i)} Take an arbitrary finite reduced irreducible root system
$R_0$ in $V=Z$ with ordered basis $\Delta_0$. Choose $\bullet\in\{u,t\}$
and let
$\Lambda,\Lambda^d$ be lattices in $V$ satisfying
\[
\mathbb{Z}R_0\subseteq\Lambda\subseteq P(R_0),\qquad
\mathbb{Z}R_0^{d}\subseteq\Lambda^d\subseteq P(R_0^{d}).
\] 
Then $(R_0,\Delta_0,\bullet,\Lambda,\Lambda^d)\in\mathcal{D}$.
Note that if $\Lambda=P(R_0)$ then $S=\emptyset$ hence $R=R^\bullet$
is reduced (Cherednik's \cite{C} theory corresponds to the special case
$(\Lambda,\Lambda^d)=(P(R_0),P(R_0^d))$).\\
{\bf (ii)} Take $Z=\mathbb{R}^{n+1}$ with standard orthonormal basis
$\{\epsilon_i\}_{i=1}^{n+1}$ and 
$R_0=\{\epsilon_i-\epsilon_j\}_{1\leq i\not=j\leq n+1}$ for the
realization of the finite root system of type $n$ in $Z$. Then
$V=(\epsilon_1+\cdots+\epsilon_{n+1})^\perp$. 
As ordered basis take 
\[
\Delta_0=(\alpha_1,\ldots,\alpha_n)=(\epsilon_1-\epsilon_2,\ldots,
\epsilon_n-\epsilon_{n+1}).
\]
Then $\bigl(R_0,\Delta_0,u,\mathbb{Z}^{n+1},
\mathbb{Z}^{n+1}\bigr)\in\mathcal{D}$.
Note that $\theta=\varphi=\epsilon_1-\epsilon_{n+1}$,
hence the simple affine root $a_0$ of $R$ is 
$a_0=c-\epsilon_1+\epsilon_{n+1}$.
This example is naturally related to the $\textup{GL}_{n+1}$
type Macdonald polynomials, see Subsection \ref{GLn}.
\end{eg}

Given a quintuple $D=(R_0,\Delta_0,\bullet,\Lambda,\Lambda^d)$
we have the dual root system $R_0^d$ with dual ordered basis
$\Delta_0^d:=(\alpha_1^d,\ldots,\alpha_n^d)$. This
extends to an involution $D\mapsto D^d$ on $\mathcal{D}$ with
\begin{equation}\label{dualquintuples}
D^d:=(R_0^d,\Delta_0^d,\bullet,\Lambda^d,\Lambda)
\end{equation}
for $D=(R_0,\Delta,\bullet,\Lambda,\Lambda^d)\in\mathcal{D}$.
We call $D^d$ the initial data dual to $D$.

We write $\mu=\mu(D)$ and $\mu^d=\mu(D^d)$ for the function
$\mu^\bullet$ on $R_0$ and $R_0^d$ respectively.
Let $\varpi_i^d\in P(R_0^d)$ ($1\leq i\leq n$) 
be the fundamental weights with respect to $\Delta_0^d$. 

For a given $D=(R_0,\Delta_0,\bullet,\Lambda,\Lambda^d)\in\mathcal{D}$
we thus have a dual triple $(R_0^d,\Delta_0^d,\bullet,\Lambda^d,\Lambda)
\in\mathcal{D}$ and hence an associated triple $(R^d,\Delta^d,W^d)$.
Concretely, the highest root $\varphi^d$ and the highest short
root $\theta^d$ of $R_0^d$ with respect to $\Delta_0^d$ are given by
\begin{equation*}
\varphi^d=
\begin{cases}
\theta^\vee\quad &\hbox{if }\, \bullet=u,\\
\varphi\quad &\hbox{if }\, \bullet=t,
\end{cases}
\end{equation*}
and
\begin{equation*}
\theta^d=
\begin{cases}
\varphi^\vee\quad &\hbox{if }\, \bullet=u,\\
\theta\quad &\hbox{if }\, \bullet=t.
\end{cases}
\end{equation*}
Hence
\[
\Delta^d=(a_0^d,a_1^d,\ldots,a_n^d)
\]
with $a_i^d=\alpha_i^d$ ($1\leq i\leq n$) and 
$a_0^d=\mu_\theta(c-\theta^\vee)$. The dual affine root system 
$R^d=R(D^d)$ is
\[
R^d=R^{d\bullet}\cup\bigcup_{i\in S^d}W^{d\bullet}(2a_i^d)
\] 
with $S^d=\{i\in\{0,\ldots,n\} \,|\, (\Lambda^d,a_i^d{}^\vee)=2\mathbb{Z}\}$,
with $R^{d\bullet}=\{m\mu_{\alpha^d}+\alpha^d\}_{m\in\mathbb{Z},\alpha\in R_0}$
and with $W^{d\bullet}\simeq W_0\ltimes\tau(\mathbb{Z}R_0)$ the affine Weyl
group of $R^{d\bullet}$. The dual extended affine Weyl group is
$W^d=W_0\ltimes\Lambda$. The simple reflections $s_i^d:=s_{a_i^d}\in 
W^{d\bullet}$ ($0\leq i\leq n$) are $s_i^d=s_i$ for $1\leq i\leq n$ and
$s_0^d=\tau(\theta)s_\theta$.

\begin{eg}\label{surprise}
The correpondence $R\leftrightarrow R^d$ can turn nonreduced affine
root systems into reduced ones. We give here an example of untwisted type.
An example for twisted type will be given in Example \ref{surprise2}.

Take $n\geq 3$ and $R_0\subseteq V=Z:=\mathbb{R}^n$ of type
$\textup{B}_n$, realized as $R_0=\{\pm\epsilon_i\}\cup
\{\pm\epsilon\pm\epsilon_j\}_{i<j}$ (all sign combinations possible),
with $\{\epsilon_i\}$ the standard orthonormal basis of $V$.
Take as ordered basis of $R_0$,
\[
\Delta_0=(\alpha_1,\ldots,\alpha_{n-1},\alpha_n)=(\epsilon_1-\epsilon_2,
\ldots,\epsilon_{n-1}-\epsilon_n,\epsilon_n).
\]
The highest root is $\varphi=\epsilon_1+\epsilon_2\in R_0$.
We then have
\[
\mathbb{Z}R_0^\vee\subset P(R_0^\vee)=\mathbb{Z}^n=\mathbb{Z}R_0\subset P(R_0)
\]
with both sublattices $\mathbb{Z}R_0^\vee\subset\mathbb{Z}^n$ and 
$\mathbb{Z}^n\subset P(R_0)$ of index two. Taking 
$\Lambda=\mathbb{Z}^n=\Lambda^d$ we get the initial data
$D=(R_0,\Delta_0,u,\mathbb{Z}^n,\mathbb{Z}^n)\in\mathcal{D}$.
Then $S=S(D)=\{n\}$ and $R=R(D)$ is given by 
\[
R=\{\pm\epsilon_i+mc\}_{1\leq i\leq n, m\in\mathbb{Z}}\cup
\{\pm\epsilon_i\pm\epsilon_j+mc\}_{1\leq i<j\leq n, m\in\mathbb{Z}}\cup
\{\pm 2\epsilon_i+2mc\}_{1\leq i\leq n, m\in\mathbb{Z}}.
\]
It is nonreduced and of untwisted type $\textup{B}_n$. We have written
it here as the disjoint union of 
the three $W=W_0\ltimes\mathbb{Z}^n$-orbits of $R$. 
Note that $a_0$ lies in the orbit 
$\{\pm\epsilon_i\pm\epsilon_j+mc\}_{1\leq i<j\leq n, m\in\mathbb{Z}}$.

Dually, $R^d=R(D^d)$ is the reduced affine root system of untwisted type
$\textup{C}_n$. Concretely,
\[
R^d=\{\pm\epsilon_i\pm\epsilon_j+mc\}_{1\leq i<j\leq n, m\in\mathbb{Z}}
\cup\{\pm 2\epsilon_i+2mc\}_{1\leq i\leq n, m\in\mathbb{Z}}\cup
\{\pm 2\epsilon_i+(2m+1)c\}_{1\leq i\leq n, m\in\mathbb{Z}},
\]
written here as the disjoint union of the three 
$W^d=W_0\ltimes\mathbb{Z}^n$-orbits of $R^d$.
\end{eg}
This example shows that basic features of the affine root system
can alter under dualization. It turns out though that the
number of orbits with respect to the action
of the extended affine Weyl group is unaltered. 
To establish this fact it is convenient to use
the concept of a multiplicity function on $R$.

\begin{defi}
Set $\mathcal{M}=\mathcal{M}(D)$ for the complex algebraic group
of $W$-invariant functions $\kappa: R\rightarrow 
\mathbb{C}^*$, and $\nu=\nu(D)$ for the complex dimension
of $\mathcal{M}$.
\end{defi}
Note that 
$\nu=\nu(D)$ equals the number of $W$-orbits of $R$.
The value of $\kappa\in\mathcal{M}$ at an affine root $a\in R$
is denoted by $\kappa_a$. We call $\kappa\in\mathcal{M}$ a multiplicity
function. We write $\kappa^\bullet:=\kappa|_{R^\bullet}$ for its restriction
to $R^\bullet$. For a multiplicity function $\kappa\in\mathcal{M}$
we set $\kappa_{2a}:=\kappa_a$ if $a\in R$ is unmultiplyable (i.e.
$2a\not\in R$). 

First we need a more precise description of the sets $S=S(D)$ and
$S^d=S(D^d)$, see \eqref{S}. It is obtained using
the classification of affine root systems (see Subsection \ref{special}
and \cite{M0}). 
\begin{lem}\label{technical}
Let $D=(R_0,\Delta_0,\bullet,\Lambda,\Lambda^d)\in\mathcal{D}$.
Set $S_0=S\cap \{1,\ldots,n\}$.
\begin{enumerate}
\item[{\bf (a)}] If $\bullet=u$ then $\#S\leq 2$. If $\#S=2$ then
$R_0$ is of type $\textup{A}_1$. If $\#S=1$ then $R_0$ is of 
type $\textup{B}_n$
($n\geq 2$) and $S=S_0=\{j\}$ with
$\alpha_j\in\Delta_0$ the unique simple short root.
\item[{\bf (b)}] If $\bullet=t$ then $\#S=0$ or $\#S=2$. If
$\#S=2$ then $R_0$ is either of type $\textup{A}_1$ or of type
$\textup{B}_n$ ($n\geq 2$) and $S=\{0,j\}$ with $\alpha_j\in\Delta_0$
the unique short simple root.
\end{enumerate}
Note that in both the untwisted and the twisted case,
$\alpha_j^d\in W_0(D(a_0^d))$ if $S_0=\{j\}$.
\end{lem}

The following lemma should be compared with \cite[\S 5.7]{Ha}.
\begin{lem}\label{dualmult}
Let $D\in\mathcal{D}$ and $\kappa\in\mathcal{M}(D)$.
Let $\alpha_j\in\Delta_0$ 
(respectively $\alpha_{j^d}^d\in\Delta_0^d$) be a simple
short root. The assignments
\begin{equation*}
\begin{split}
&\kappa^d_{a_0^d}:=\kappa_{2\alpha_j},\\
&\kappa^d_{\alpha_i^d}:=\kappa_{\alpha_i}\,\,\,\,\quad\hbox{ if }\,\,
i\in\{1,\ldots,n\},\\
&\kappa_{2a_0^d}^d:=\kappa_{2a_0}\quad\, \hbox{ if }\,\, 0\in S^d,\\
&\kappa_{2\alpha_{j^d}^d}^d:=\kappa_{a_0}\quad\, \hbox{ if }\,\,
j^d\in S^d
\end{split}
\end{equation*}
uniquely extend to a multiplicity function
$\kappa^d\in\mathcal{M}^d:=\mathcal{M}(D^d)$. For fixed $D\in\mathcal{D}$
the map $\kappa\mapsto \kappa^d$ defines an isomorphism
$\phi_D: \mathcal{M}(D)\overset{\sim}{\longrightarrow}\mathcal{M}(D^d)$
of complex tori, with inverse $\phi_{D^d}$.
\end{lem}
\begin{rema}
Let $\kappa\in\mathcal{M}(D)$. 
Recall the convention that $\kappa_{2a}=\kappa_a$ for $a\in R$ such that
$2a\not\in R$. Then for all $\alpha\in R_0$,
\[
\kappa_{2\alpha}=\kappa_{\mu_{\alpha^d}c+\alpha^d}^d.
\]
\end{rema}
\begin{cor}\label{equalorbits}
Let $D\in\mathcal{D}$.
The number $\nu$ of $W$-orbits of $R$ is equal to the number 
$\nu^d$ of
$W^d$-orbits of $R^d$.
\end{cor}

\begin{rema}
Returning to Example \ref{surprise}, the 
correspondence from Lemma \ref{dualmult} links the orbit
$\{\pm\epsilon_i+mc\}$ of $R$ to the orbit
$\{2\epsilon_i+2mc\}$ of $R^d$,
the orbit $\{\pm\epsilon_i\pm\epsilon_j+mc\}$ of $R$ to
the orbit $\{\pm\epsilon_i\pm\epsilon_j+mc\}$ of $R^d$ and
the orbit $\{\pm 2\epsilon_i+2mc\}$ of $R$ to 
the orbit $\{\pm 2\epsilon_i+(2m+1)c\}$ of $R^d$.
\end{rema}

\subsection{The basic representation}\label{extension}

We fix throughout this subsection a quintuple $D=(R_0,\Delta_0,\bullet,
\Lambda,\Lambda^d)\in\mathcal{D}$ of initial
data. Recall that it gives rise to a triple $(R,\Delta,W)$ of
an irreducible affine root system $R$ containing $R^\bullet$, an ordered
basis $\Delta$ of $R$ as well as of $R^\bullet$, 
and an extended affine Weyl group $W=W_0\ltimes\Lambda^d$. In addition
we fix a multiplicity function $\kappa\in\mathcal{M}(D)$ and we
write $\kappa^\bullet:=\kappa|_{R^\bullet}$. It is a $W$-equivariant
map $R^\bullet\rightarrow\mathbb{C}^*$. Write 
$\kappa_i:=\kappa_{a_i}^\bullet$ for $0\leq i\leq n$. Note that 
$\kappa_i=\kappa_j$ if $s_i$ is conjugate to $s_j$ in $W=\Omega\ltimes 
W^\bullet$.

We write $R^{\pm}$ and $R^{\bullet\pm}$
for the positive respectively negative affine roots of $R$ and $R^\bullet$ with
respect to $\Delta$.
Since the affine root system $R^\bullet$ is $W$-stable, we can define the
length $l(w)=l_D(w)$
of $w\in W$ by
\[
l(w):=\#\bigl(R^{\bullet+}\cap w^{-1}R^{\bullet-}\bigr).
\]
If $w\in W^\bullet=W(R^\bullet)$ 
then $l(w)$ equals the number of simple reflections $s_i$
($0\leq i\leq n$) in a reduced 
expression of $w$. We have $W=\Omega\ltimes W^\bullet$ with
$\Omega=\Omega(D)$ the subgroup
\[\Omega:=\{w\in W \,\, | \,\, l(w)=0\}
\]
of $W$. Then $\Omega\simeq \Lambda^d/\mathbb{Z}R_0^d$. 
The abelian group $\Omega$ 
permutes the simple affine roots $a_i$ ($0\leq i\leq n$), which thus gives
rise to an action of $\Omega$ on the index set $\{0,\ldots,n\}$.
Consequently
the action of $\Omega$ on $W^\bullet$ by conjugation permutes the
set $\{s_i\}_{i=0}^n$ of simple reflections,
$ws_iw^{-1}=s_{w(i)}$ for $w\in\Omega$ and $0\leq i\leq n$ 
(cf., e.g., \cite[\S 2.5]{M}). 
A detailed description
of the group $\Omega$ in terms of a complete set of representatives
of $\Lambda^d/\mathbb{Z}R_0^d$ 
will be given in Subsection \ref{Macdonalddiff}.

Extended versions of the affine braid group and of the affine 
Hecke algebra are defined as follows. 
Let $A=A(D):=A(R^\bullet,\Delta)$ be the affine Cartan matrix
associated to $(R^\bullet,\Delta)$. Recall that the affine braid group 
$\mathcal{B}^\bullet:=\mathcal{B}(A)$ is isomorphic to the abstract
group generated by $T_w$ ($w\in W^\bullet$) with defining relations
$T_vT_w=T_{vw}$ for all $v,w\in W^\bullet$ satisfying $l(vw)=l(v)+l(w)$. 
\begin{defi}\label{ea}
{\bf (i)} The extended affine braid group 
$\mathcal{B}=\mathcal{B}(D)$ is the 
group generated by $T_w$ ($w\in W$) with defining
relations $T_vT_w=T_{vw}$ for all $v,w\in W$ satisfying
$l(vw)=l(v)+l(w)$.\\
{\bf (ii)} 
The extended affine Hecke algebra $H(\kappa^\bullet)=
H(D,\kappa^\bullet)$ is the quotient of
$\mathbb{C}[\mathcal{B}]$ by the two-sided ideal generated
by $(T_i-\kappa_i)(T_i+\kappa_i^{-1})$ ($0\leq i\leq n$).
\end{defi}
Similarly to the semidirect product decomposition
$W\simeq\Omega\ltimes W^\bullet$ we have
$\mathcal{B}\simeq\Omega\ltimes\mathcal{B}^\bullet$ and
\[
H(\kappa^\bullet)\simeq \Omega\ltimes H(W^\bullet,\kappa^\bullet),
\]
where the action of $\Omega$ on $\mathcal{B}$ by group automorphisms
(respectively on $H(W^\bullet,\kappa^\bullet)$ by algebra automorphisms)
is determined by $w\cdot T_i=T_{w(i)}$ for $w\in\Omega$ and $0\leq i\leq n$.
For $\omega\in\Omega$ we will denote the element
$T_\omega$ in the extended affine Hecke algebra 
$H(\kappa^\bullet)$ simply by $\omega$.

The algebra homomorphism 
\[\beta_{(R,\Delta,\kappa)}: H(W^\bullet,\kappa^\bullet)
\hookrightarrow W^\bullet\ltimes F\subseteq
W\ltimes F
\]
from Subsection \ref{realizations} extends to an injective algebra
map
\[
\beta_{D,\kappa}: H(\kappa^\bullet)\hookrightarrow W\ltimes F
\] 
by $\beta_{D,\kappa}(T_\omega)=\omega$ ($\omega\in\Omega$).
We will now
show that it gives rise to an action of $H(\kappa^\bullet)$
as $q$-difference reflection operators on a complex torus $T_\Lambda$.
It is called
the basic representation of 
$H(\kappa^\bullet)$. It is fundamental for the development 
of the Cherednik-Macdonald theory.

The complex torus $T_\Lambda$ (of rank $\textup{dim}_{\mathbb{R}}(Z)$) is the
algebraic group of complex characters of the lattice $\Lambda$,
\[
T_\Lambda:=\textup{Hom}_{\mathbb{Z}}\bigl(\Lambda,\mathbb{C}^*\bigr).
\]
The algebra $\mathbb{C}[T_\Lambda]$ of regular functions on $T_\Lambda$
is isomorphic to the group algebra $\mathbb{C}[\Lambda]$, where
the standard basis element $e^\lambda$ ($\lambda\in\Lambda$) of 
$\mathbb{C}[\Lambda]$ is
viewed as the regular function $t\mapsto t(\lambda)$ on $T_\Lambda$.
We write $t^\lambda$ for the value of $e^\lambda$ at $t\in T_\Lambda$.
Since $\Lambda$ is $W_0$-stable, $W_0$ acts on $T_\Lambda$, 
giving in turn rise to an action of $W_0$
on $\mathbb{C}[T_\Lambda]$ by algebra automorphisms. 
Then $w(e^\lambda)=e^{w\lambda}$ ($w\in W_0$ and $\lambda\in\Lambda$).
We now first extend it to an action of the extended
affine Weyl group $W$ on $\mathbb{C}[T_\Lambda]$ depending on
a fixed parameter $q\in\mathbb{R}_{>0}\setminus\{1\}$.

Set
\[
q_\alpha:=q^{\mu_\alpha},\qquad \alpha\in R_0
\]
and define 
$q^\xi\in T_\Lambda$ ($\xi\in\Lambda^d$) to be
the character $\lambda\mapsto q^{(\lambda,\xi)}$
of $\Lambda$. The action of $W_0$ on $T_\Lambda$
extends to a left $W$-action
$(w,t)\mapsto w_qt$ on $T_\Lambda$ by 
\[
\tau(\xi)_qt:=q^\xi t,\qquad \xi\in\Lambda^d,\,\, t\in T_\Lambda.
\]
Then $(w_qp)(t):=p(w_q^{-1}t)$ ($w\in W$, $p\in\mathbb{C}[T_\Lambda]$)
is a $W$-action by algebra automorphisms on $\mathbb{C}[T_\Lambda]$.
In particular,
\[
\tau(\xi)_q(e^\lambda)=q^{-(\lambda,\xi)}e^\lambda,\qquad
\xi\in\Lambda^d,\,\,\lambda\in\Lambda.
\]
It extends to a $W$-action by field automorphisms
on the quotient field $\mathbb{C}(T_\Lambda)$ of
$\mathbb{C}[T_\Lambda]$. 
It is useful to introduce the notation
\[
t_q^{rc+\lambda}:=q^{r}t^\lambda,
\qquad r\in\mathbb{R},\, \lambda\in\Lambda
\]
for $t\in T_\Lambda$. Then $(w_q^{-1}t)_q^{rc+\lambda}=
t_q^{w(rc+\lambda)}$
for $w\in W$, $t\in T_\Lambda$, $r\in\mathbb{R}$
and $\lambda\in\Lambda$.

We write $W\ltimes_q\mathbb{C}(T_\Lambda)$
for the resulting semidirect product algebra. It canonically acts
on $\mathbb{C}(T_\Lambda)$ by $q$-difference reflection operators. 
We thus have a sequence of algebra maps
\[
H(\kappa^\bullet)\rightarrow
W\ltimes F\rightarrow W\ltimes_q\mathbb{C}(T_\Lambda)
\rightarrow\textup{End}_{\mathbb{C}}\bigl(\mathbb{C}(T_\Lambda)\bigr),
\]
where the first map is $\beta_{D,\kappa}$ and the
second map sends $e^{\mu_\alpha mc+\alpha}$ to
$q_\alpha^{m}e^{\alpha}$ for $m\mu_\alpha c+\alpha\in R^\bullet$.
It gives the following result, which is closely related
to \cite[5.13]{Ha} in the present generality.

\begin{thm}
Let $D=(R_0,\Delta_0,\bullet,\Lambda,\Lambda^d)\in\mathcal{D}$
and $\kappa\in\mathcal{M}(D)$.
For $a\in R$ set $\kappa_{2a}:=\kappa_a$ if $2a\not\in R$.
There exists a unique algebra monomorphism
\[\pi_{\kappa,q}=\pi_{D;\kappa,q}:
H(D,\kappa^\bullet)\hookrightarrow \textup{End}_{\mathbb{C}}\bigl(
\mathbb{C}[T_\Lambda]\bigr)
\]
satisfying
\begin{equation*}
\begin{split}
\bigl(\pi_{\kappa,q}(T_i)p\bigr)(t)&=\kappa_{a_i}(s_{i,q}p)(t)+
\left(\frac{\kappa_{a_i}-\kappa_{a_i}^{-1}+
(\kappa_{2a_i}-\kappa_{2a_i}^{-1})t_q^{a_i}}{1-t_q^{2a_i}}\right)
(p(t)-(s_{i,q}p)(t)),\\
\bigl(\pi_{\kappa,q}(\omega)p\bigr)(t)&=(\omega_qp)(t)
\end{split}
\end{equation*}
for $0\leq i\leq n$, $\omega\in\Omega$, 
$p\in\mathbb{C}[T_\Lambda]$ and $t\in T_\Lambda$.
\end{thm}
If $2a_i\not\in R$ then, by the convention $\kappa_{2a_i}=\kappa_{a_i}$, 
the first formula reduces to
\[
\bigl(\pi_{\kappa,q}(T_i)p)\bigr)(t)=
\kappa_{a_i}(s_{i,q}p)(t)+
\Bigl(\frac{\kappa_{a_i}-\kappa_{a_i}^{-1}}{1-t_q^{a_i}}\Bigr)
(p(t)-(s_{i,q}p)(t)).
\]

The theorem is due to Cherednik (see \cite{C} and references therein)
in the $\textup{GL}_{n+1}$ case (see Example \ref{lattices}{\bf (ii)})
and when $D=(R_0,\Delta_0,\bullet,P(R_0),P(R_0^d))$ 
with $R_0$ an arbitrary reduced irreducible root system.
The theorem is due to Noumi \cite{N} 
for $D=(R_0,\Delta_0,t,\mathbb{Z}R_0,\mathbb{Z}R_0)$ 
with $R_0$ of type $\textup{A}_1$ or of
type $\textup{B}_n$ ($n\geq 2$). This case is special due to
its large degree of freedom ($\nu(D)=4$ if $n=1$ and $\nu(D)=5$
if $n\geq 2$). We will
describe this case in detail in Subsection \ref{CcheckC}.

\begin{rema}
{\bf (i)} {}From Theorem \ref{beta} one first obtains $\pi_{\kappa,q}$
as an algebra map from $H(\kappa^\bullet)$ to 
$\textup{End}_{\mathbb{C}}\bigl(\mathbb{C}(T_\Lambda)\bigr)$. 
The image is contained in
the subalgebra of endomorphisms 
preserving $\mathbb{C}[T_\Lambda]$ since, for $\lambda\in\Lambda$,
we have $(\lambda,a_i^\vee)\in\mathbb{Z}$ if $2a_i\not\in R$
and $(\lambda,a_i^\vee)\in 2\mathbb{Z}$ if $2a_i\in R$.\\
{\bf (ii)} Let $s=s(D)$ be the number of $W$-orbits of $R\setminus R^\bullet$.
Extending a $W$-equivariant map $\kappa^\bullet: R^\bullet
\rightarrow\mathbb{C}^*$ to a multiplicity function $\kappa\in\mathcal{M}$
on $R$ amounts to 
choosing $s$ nonzero complex parameters. 
Hence, 
the maps $\pi_{\kappa,q}$ define a family of algebra monomorphisms of 
the extended affine Hecke algebra $H(\kappa^\bullet)$ into 
$\textup{End}_{\mathbb{C}}\bigl(\mathbb{C}[T_\Lambda]\bigr)$,
parametrized by $s+1$ parameters $\kappa|_{R\setminus R^\bullet}$
and $q$.
\end{rema} 

\section{Monic Macdonald-Koornwinder polynomials}\label{3}
In this section we introduce the monic nonsymmetric and symmetric
Macdonald-Koorn\-win\-der polynomials. The terminology Macdonald
polynomials is employed in the literature for the cases that $R$ is reduced
(i.e. the cases $D=(R_0,\Delta_0,\bullet,P(R_0),P(R_0^d))$
and the $\textup{GL}_{n+1}$ case). The Koornwinder polynomials
correspond to the initial data
$D=(R_0,\Delta_0,t,\mathbb{Z}R_0,\mathbb{Z}R_0)$
with $R_0$ of type $\textup{A}_1$ of of type $\textup{B}_n$ ($n\geq 2$),
in which case $R=R(D)$ is nonreduced and of type $C^\vee C_n$.
To have uniform terminology
we will speak of Macdonald-Koornwinder polynomials when discussing the theory
for arbitrary initial data.

The monic nonsymmetric
Macdonald-Koornwinder polynomials will be introduced as the
common eigenfunctions in $\mathbb{C}[T_\Lambda]$ of a family of commuting
$q$-difference reflection operators. The operators are obtained
as images under the basic representation $\pi_{\kappa,q}$ of elements from
a large commutative subalgebra of the extended affine Hecke algebra
$H(\kappa^\bullet)$ constructed by Bernstein and Zelevinsky. 
A Hecke algebra symmetrizer turns the monic nonsymmetric 
Macdonald-Koornwinder
polynomials into monic symmetric 
Macdonald-Koornwinder polynomials, which are $W_0$-invariant regular
functions on $T_\Lambda$ solving a suitable spectral problem of a
commuting family of $q$-difference operators, called
Macdonald operators.
We in addition
determine in this section
the (bi)orthogonality relations of the polynomials.

Throughout this section we fix 
\begin{enumerate}
\item a quintuple $D=(R_0,\Delta_0,\bullet,\Lambda,\Lambda^d)$
of initial data,
\item a deformation parameter $q\in\mathbb{R}_{>0}\setminus\{1\}$,
\item a multiplicity function $\kappa\in\mathcal{M}(D)$ 
\end{enumerate}
and we freely use the resulting notations from the previous section.
\subsection{Bernstein-Zelevinsky presentation}\label{Bernstein}
An expression $w=\omega s_{i_1}s_{i_2}\cdots s_{i_r}\in W$
with $\omega\in\Omega$ and $0\leq i_j\leq n$ 
is called a reduced expression if $r=l(w)$. 

For a given reduced expression $w=\omega s_{i_1}s_{i_2}\cdots s_{i_r}$,
\[
T_w:=\omega T_{i_1}T_{i_2}\cdots T_{i_r}\in H(\kappa^\bullet)
\]
is well defined, and $\{T_w\}_{w\in W}$ is a complex linear basis
of $H(\kappa^\bullet)$.

The cones
\[\Lambda^{d\pm}:=\{\xi\in\Lambda^d \,\, | \,\, 
\bigl(\xi,\alpha^{d\vee}\bigr)\geq 0\quad \forall\, \alpha\in R_0^\pm 
\}
\]
form fundamental domains for the $W_0$-action on $\Lambda^d$.
Any $\xi\in\Lambda^d$ can be written as
$\xi=\mu-\nu$ with $\mu,\nu\in\Lambda^{d+}$ (and similarly for
$\Lambda^{d-}$). Furthermore,
if $\xi,\xi^\prime\in\Lambda^{d+}$ then 
$l(\tau(\xi+\xi^\prime))=l(\tau(\xi))+l(\tau(\xi^\prime))$.
It follows that there exists a unique group homomorphism 
$\Lambda^d\rightarrow\mathcal{B}$,
denoted by $\xi\mapsto Y^\xi$, such that
\[
Y^\xi=T_{\tau(\xi)}\qquad \forall\,\xi\in\Lambda^{d+}.
\]
On the level of the extended affine Hecke algebra it gives rise
to an algebra homomorphism
\begin{equation}\label{Ymap}
\mathbb{C}[T_{\Lambda^d}]\rightarrow H(\kappa^\bullet)
\end{equation}
which we denote by $p\mapsto p(Y)$,
where $T_{\Lambda^d}$ is the complex torus
$\textup{Hom}(\Lambda^d,\mathbb{C}^*)$, and
$p(Y)=\sum_\xi c_\xi Y^\xi$ if 
$p(t)=\sum_\xi c_\xi t^\xi$. 
The image of the map \eqref{Ymap}
is denoted by $\mathbb{C}_Y[T_{\Lambda^d}]$. 

As in Subsection \ref{extension}, the lattice $\Lambda^d$ is
$W_0$-stable, giving rise to a $W_0$-action on $T_{\Lambda^d}$, hence
a $W_0$-action on $\mathbb{C}[T_{\Lambda^d}]$ by algebra automorphisms.

Denote by $H_0$ the subalgebra of $H(\kappa^\bullet)$
generated by $T_i$ ($1\leq i\leq n$). 
We have a natural surjective algebra map 
$H(W_0,\kappa|_{R_0})\rightarrow H_0$, sending the algebraic generator $T_i$
of $H(W_0,\kappa|_{R_0})$ to $T_i\in H_0$.

The analog of the semidirect product decomposition
$W=W_0\ltimes\Lambda^d$ for the extended
affine Hecke algebra $H(\kappa^\bullet)$ is the following
Bernstein-Zelevinsky presentation of $H(\kappa^\bullet)$ (see \cite{Lu}).
\begin{thm}
\begin{enumerate}
\item The algebra maps $\mathbb{C}[T_{\Lambda^d}]\rightarrow
\mathbb{C}_Y[T_{\Lambda^d}]$ 
and $H(W_0,\kappa|_{R_0})\rightarrow H_0$ are isomorphisms,
\item multiplication defines a linear isomorphism
$H_0\otimes\mathbb{C}_Y[T_{\Lambda^d}]\simeq H(\kappa^\bullet)$,
\item for $i\in\{1,\ldots,n\}$ such that
$\bigl(\Lambda^d,\alpha_i^{d\vee}\bigr)=\mathbb{Z}$ 
we have
\begin{equation}\label{Lusztig1}
p(Y)T_i-T_i(s_ip)(Y)=\left(\frac{\kappa_i-\kappa_i^{-1}}
{1-Y^{-\alpha_i^{d}}}\right)(p(Y)-(s_ip)(Y))
\end{equation}
in $H(\kappa^\bullet)$ for all $p\in\mathbb{C}[T_{\Lambda^d}]$,
\item for $i\in\{1,\ldots,n\}$ such that 
$\bigl(\Lambda^d,\alpha_i^{d\vee}\bigr)=2\mathbb{Z}$ 
we have
\begin{equation}\label{Lusztig2}
p(Y)T_i-T_i(s_ip)(Y)=\left(\frac{\kappa_i-\kappa_i^{-1}+
(\kappa_0-\kappa_0^{-1})Y^{-\alpha_i^{d}}}
{1-Y^{-2\alpha_i^{d}}}\right)(p(Y)-(s_ip)(Y))
\end{equation}
in $H(\kappa^\bullet)$ for all $p\in\mathbb{C}[T_{\Lambda^d}]$.
\end{enumerate}
These properties characterize $H(\kappa^\bullet)$ as
a unital complex associative algebra.
\end{thm}
With the notion of the dual multiplicity parameter $\kappa^d$
(see Lemma \ref{dualmult}), the cross relations
\eqref{Lusztig1} and \eqref{Lusztig2} in the affine Hecke algebra
$H(\kappa^\bullet)$ can be
uniformly written as
\begin{equation}\label{Lusztiguniform}
p(Y)T_i-T_i(s_ip)(Y)=\left(\frac{\kappa_{\alpha_i^d}^d-
(\kappa_{\alpha_i^d}^d)^{-1}
+(\kappa_{2\alpha_i^d}^d-(\kappa_{2\alpha_i^d}^d)^{-1})Y^{-\alpha_i^d}}
{1-Y^{-2\alpha_i^d}}\right)(p(Y)-(s_ip)(Y))
\end{equation}
for $p\in\mathbb{C}[T_{\Lambda^d}]$ and $1\leq i\leq n$. 

It follows from the theorem that
the center $Z(H(\kappa^\bullet))$ of the extended affine Hecke algebra
$H(\kappa^\bullet)$ equals $\mathbb{C}_Y[T_{\Lambda^d}]^{W_0}$.
\subsection{Monic nonsymmetric Macdonald-Koornwinder
polynomials}\label{MnMp}
The results in this section are from \cite{MAst,Cns,Sa,Ha}.
For detailed proofs see, e.g., \cite[\S 2.8, \S 4.6, \S 5.2]{M}.

We put in this subsection the 
following conditions on $q$ 
and $\kappa\in\mathcal{M}(D)$,
\begin{equation}\label{smallparameter}
0<q<1,\qquad 0<\kappa_a<1\,\,\,
\forall a\in R,
\end{equation}
or
\begin{equation}\label{largeparameter}
q>1,\qquad \kappa_a>1\,\,\,
\forall a\in R.
\end{equation}
Set
\begin{equation}\label{eta}
\eta(x)=
\begin{cases}
1\quad &\hbox{ if } x>0,\\
-1\quad &\hbox{ if } x\leq 0
\end{cases}
\end{equation}
and define a $W_0$-equivariant map $\upsilon: 
R_0^d\rightarrow \mathbb{R}_{>0}$
(depending on $\kappa^\bullet$) by
\begin{equation}\label{kappa}
\upsilon_{\alpha^d}:=
\kappa_{\alpha}^{\frac{1}{2}}\kappa_{\mu_\alpha c+\alpha}^{\frac{1}{2}},
\qquad \alpha\in R_0.
\end{equation}
\begin{defi}\label{gammapoint}
For $\lambda\in\Lambda$ define 
$\gamma_{\lambda,q}=\gamma_{\lambda,q}(D;\kappa^\bullet)\in T_{\Lambda^d}$
by
\[
\gamma_{\lambda,q}:=q^{\lambda}\prod_{\alpha\in R_0^{+}}
\upsilon_{\alpha^d}^{\eta((\lambda,\alpha^\vee))\alpha^{d\vee}}.
\]
In other words, for all $\xi\in\Lambda^d$,
\[
\gamma_{\lambda,q}^\xi=q^{(\lambda,\xi)}\prod_{\alpha\in R_0^+}
\upsilon_{\alpha^d}^{\eta((\lambda,\alpha^\vee))(\xi,\alpha^{d\vee})}.
\]
\end{defi}
As a special case,
\[
\gamma_{\lambda,q}=q^{\lambda}\prod_{\alpha\in R_0^{+}}
\upsilon_{\alpha^d}^{-\alpha^{d\vee}}
\qquad \forall\,\lambda\in\Lambda^-,
\]
where 
$\Lambda^{\pm}=\{\lambda\in\Lambda\,\, | \,\, 
\bigl(\lambda,\alpha^\vee\bigr)\geq 0
\quad \forall\, \alpha\in R_0^\pm\}$.

Write $l^d=l_{D^d}$ for the length function on the dual extended
affine Weyl group $W^d$ and 
$\Omega^d=\Omega(D^d)$ for the subgroup of elements
of $W^d$ of length zero with respect to $l^d$.
We have a $q$-dependent $W^d$-action 
on $T_{\Lambda^d}$ extending the $W_0$-action by
$\tau(\lambda)_q\gamma=q^\lambda\gamma$ for all $\lambda\in\Lambda$ 
and $\gamma\in T_{\Lambda^d}$.
Then $\gamma_{\lambda,q}=\tau(\lambda)_q\gamma_{0,q}$ 
in $T_{\Lambda^d}$ if 
$\lambda\in\Lambda^-$. This generalizes as follows. 

\begin{lem}\label{XYZ}
We have for $\lambda\in\Lambda$,
\[
\gamma_{\lambda,q}=u^d(\lambda)_q\gamma_{0,q}
\]
in $T_{\Lambda^d}$, where 
$u^d(\lambda)\in W^d$ is the element of minimal
length with respect to $l^d$ in the coset $\tau(\lambda)W_0$.
\end{lem}
The conditions \eqref{smallparameter} or \eqref{largeparameter}
on the parameters, together with Lemma \ref{XYZ}, imply
\begin{lem}
The map $\Lambda\rightarrow T_{\Lambda^d}$, defined by $\lambda\mapsto
\gamma_{\lambda,q}$, is injective.
\end{lem}
For later purposes it is convenient to record the following
compatibility between the $q$-dependent
$W^d$-action on $\gamma_{\lambda,q}\in T_{\Lambda^d}$
and the $W^d$-action $(w\tau(\lambda),\lambda^\prime)\mapsto
w(\lambda+\lambda^\prime)$ on $\Lambda$ ($w\in W_0$ and 
$\lambda,\lambda^\prime\in\Lambda$). 

\begin{prop}\label{Wgammaaction}
Let $\lambda\in \Lambda$. Then
\begin{enumerate}
\item[{\bf (a)}] If $\omega\in\Omega^d$ then 
$\omega_q\gamma_{\lambda,q}=\gamma_{\omega\lambda,q}$.
\item[{\bf (b)}] If $0\leq i\leq n$ and $s_i^d\lambda\not=\lambda$
then $s_{i,q}^d\gamma_{\lambda,q}=\gamma_{s_i^d\lambda,q}$.
\item[{\bf (c)}] If $0\leq i\leq n$ and $s_i^d\lambda=\lambda$ then
$s_{i,q}^d\gamma_{\lambda,q}
=\gamma_{\lambda,q}\upsilon_{D(a_i^d)}^{2D(a_i^d)^{\vee}}$.
\end{enumerate}
\end{prop}
\noindent
{\bf Warning:} $D(a_i^d)=(Da_i)^d$ holds true for $1\leq i\leq n$
and for $i=0$ if $\bullet=t$ (then both sides equal $-\theta$). 
It is not correct when $i=0$,
$\bullet=u$ and $R_0$ has two root lengths, since then
$(Da_0)^d=-\varphi^\vee$ and $D(a_0^d)=-\theta^\vee$.

For $\lambda,\mu\in\Lambda^+$ we write $\lambda\leq\mu$ if
$\mu-\lambda$ can be written as a sum of positive roots 
$\alpha\in R_0^{+}$. We also write $\leq$ for the Bruhat order
of $W_0$ with respect to the Coxeter generators $s_i$ ($1\leq i\leq n$). 
For $\lambda\in\Lambda$ let $\lambda_{\pm}$ 
be the unique element in $\Lambda^{\pm}\cap W_0\lambda$ and
write $v(\lambda)\in W_0$
for the element of shortest length such that $v(\lambda)\lambda=\lambda_-$.
Then $\tau(\lambda)=u^d(\lambda)v(\lambda)$ in $W^d$
for $\lambda\in\Lambda$. 
\begin{defi}
Let $\lambda,\mu\in\Lambda$. We write
$\lambda\preceq\mu$ if $\lambda_+<\mu_+$ or if $\lambda_+=\mu_+$
and $v(\lambda)\geq v(\mu)$. 
\end{defi}
Note that $\preceq$ is a partial order on $\Lambda$. Furthermore, 
if $\lambda\in\Lambda^-$ then $\mu\preceq\lambda$ for all
$\mu\in W_0\lambda$.
For each $\lambda\in\Lambda$ the set of elements $\mu\in \Lambda$ satisfying 
$\mu\preceq\lambda$ thus is contained in the finite set
\[
\{\mu\in \Lambda \, | \, \mu\preceq\lambda_-\}=
\bigcup_{\mu_+\in \Lambda^+: \mu_+\leq\lambda_+}W_0\mu_+,
\]
which is the smallest saturated subset $\textup{Sat}(\lambda_+)$
of $\Lambda$ containing $\lambda_+$ 
(a subset $X\subseteq \Lambda$ is saturated if for each 
$\alpha\in R_0^+$ and $\lambda\in \Lambda$ we have
$\lambda-r\alpha\in X$ for all integers $r$ between zero and
$\bigl(\lambda,\alpha^\vee\bigr)$, including both zero and
$\bigl(\lambda,\alpha^\vee\bigr)$). 

We write $p=d_{\lambda}e^\lambda+\textup{ l.o.t. }$ 
for an element $p=\sum_{\mu\in\Lambda}d_\mu e^\mu\in\mathbb{C}[T_\Lambda]$
satisfying $d_\mu=0$ if $\mu\not\preceq\lambda$. If in addition 
$d_\lambda\not=0$ then we call $p$ of degree $\lambda$.
\begin{prop}\label{triangularity}
For $r\in\mathbb{C}[T_{\Lambda^d}]$ and $\lambda\in \Lambda$ we have
\[
\pi_{\kappa,q}(r(Y))e^\lambda=
r(\gamma_{\lambda,q}^{-1})e^\lambda+\textup{ l.o.t.}
\]
in $\mathbb{C}[T_\Lambda]$.
\end{prop}
 
\begin{cor}
For each $\lambda\in\Lambda$ there exists a unique 
$P_\lambda=P_\lambda(D;\kappa,q)\in\mathbb{C}[T_{\Lambda}]$ satisfying
\[
\pi_{\kappa,q}(r(Y))P_\lambda=r(\gamma_{\lambda,q}^{-1})P_\lambda,\qquad
\forall\, r\in\mathbb{C}[T_{\Lambda^d}]
\]
and satisfying $P_\lambda=e^\lambda+\textup{l.o.t.}$
\end{cor}
\begin{defi}
$P_\lambda=P_\lambda(D;\kappa,q)\in\mathbb{C}[T_\Lambda]$ 
is called the monic nonsymmetric Macdonald-Koornwinder polynomial
of degree $\lambda\in\Lambda$.
\end{defi}
For $D=(R_0,\Delta_0,\bullet,P(R_0),P(R_0^d))$ 
the definition of the nonsymmetric
Mac\-do\-nald-Koorn\-win\-der 
polynomial is due to Macdonald \cite{MAst} in the untwisted case
($\bullet=u$) and due to Cherednik \cite{Cns}
in the general case. For $D=(R_0,\Delta_0,t,\mathbb{Z}R_0,\mathbb{Z}R_0)$
with $R_0$ of type $\textup{A}_1$ or of type $\textup{B}_n$ ($n\geq 2$)
the nonsymmetric Macdonald-Koornwinder 
polynomials are Sahi's \cite{Sa} nonsymmetric
Koornwinder polynomials. In the present generality (with more flexible 
choices of lattices $\Lambda$ and $\Lambda^d$) the definition is close
to Haiman's \cite[\S 6]{Ha} definition. The same references apply for the
biorthogonality relations of the nonsymmetric Macdonald-Koornwinder 
polynomials discussed in the next subsection.

The $\textup{GL}_{n+1}$
nonsymmetric Macdonald polynomials (corresponding to
Example \ref{lattices}{\bf (ii)}) are often studied separately, see, 
e.g., \cite{Kn,HH}.

\subsection{Biorthogonality}\label{sectionbiortho}
We assume in this subsection that $\kappa\in\mathcal{M}$
and $q$ satisfy
\begin{equation}\label{smallparameterorth}
\begin{split}
&0<q<1,\\
&0<\kappa_a<1\qquad\quad \forall\, a\in R,\\
&0<\kappa_a\kappa_{2a}^{\pm 1}\leq 1\qquad \forall\, a\in R^\bullet
\end{split}
\end{equation}
and write for $\lambda\in\Lambda$,
\[
P_\lambda:=P_\lambda(D,\kappa,q),\qquad
P_\lambda^\circ:=P_\lambda(D,\kappa^{-1},q^{-1}),
\]
where $\kappa^{-1}\in\mathcal{M}(D)$ is the
multiplicity function $a\mapsto\kappa_a^{-1}$.

For $a\in R^\bullet$ define 
$c_a=c_a^{\kappa,q}(\cdot;D)\in\mathbb{C}(T_\Lambda)$ by
\begin{equation}\label{ca}
c_a(t)=
\frac{(1-\kappa_a\kappa_{2a}t_q^a)(1+\kappa_a\kappa_{2a}^{-1}t_q^a)}
{(1-t_q^{2a})}.
\end{equation}
Then for $w\in W$ and $a\in R^\bullet$,
\[
c_a(w_q^{-1}t)=c_{wa}(t).
\]
In addition,
\[
\pi_{\kappa,q}(T_i)=\kappa_i+\kappa_i^{-1}c_{a_i}(s_{i,q}-1)
\]
for $0\leq i\leq n$.

Since $0<q<1$, the infinite product 
\begin{equation*}
v:=\prod_{a\in R^{\bullet+}}\frac{1}{c_a}
\end{equation*}
defines a meromorphic function on $T_\Lambda$. In terms of products
\begin{equation}\label{qshifted}
\bigl(x_1,\ldots,x_m;q\bigr)_r:=\prod_{i=1}^m\prod_{j=0}^{r-1}(1-q^jx_i),
\qquad r\in\mathbb{Z}_{\geq 0}\cup\{\infty\}
\end{equation}
of $q$-shifted factorials it becomes
\begin{equation*}
\begin{split}
v(t)&=\prod_{\alpha\in R_0^{+}}\frac{1-t^{2\alpha}}
{(1-\kappa_\alpha\kappa_{2\alpha}t^\alpha)
(1+\kappa_\alpha\kappa_{2\alpha}^{-1}t^\alpha)}\\
&\times\prod_{\beta\in R_0}
\frac{\bigl(q_\beta^{2}t^{2\beta};q_\beta^{2}\bigr)_{\infty}}
{\bigl(q_\beta^{2}\kappa_\beta\kappa_{2\beta}t^\beta,
-q_\beta^{2}\kappa_\beta\kappa_{2\beta}^{-1}t^\beta,
q_\beta\kappa_{\mu_\beta c+\beta}\kappa_{2\mu_\beta c+2\beta}t^\beta,
-q_\beta\kappa_{\mu_\beta c+\beta}
\kappa_{2\mu_\beta c+2\beta}^{-1}t^\beta;q_\beta^{2}
\bigr)_{\infty}}.
\end{split}
\end{equation*}
\begin{rema}
If $\beta\in R_0$ satisfies
$\bigl(\Lambda^d,\beta^{d\vee}\bigr)=\mathbb{Z}$ 
then $\kappa_{\mu_\beta c+\beta}=\kappa_\beta$
and $\kappa_{2\mu_\beta c+2\beta}=\kappa_{2\beta}$. In that case
the $\beta$-factor in the second line of the above formula
of $v$ simplifies to
\[
\frac{\bigl(q_\beta^{2}t^{2\beta};q_\beta^{2}\bigr)_{\infty}}
{\bigl(q_\beta\kappa_\beta\kappa_{2\beta}t^\beta,
-q_\beta\kappa_\beta\kappa_{2\beta}^{-1}t^\beta;
q_\beta\bigr)_{\infty}}.
\]
If in addition $\kappa_{2\beta}=\kappa_\beta$ (for instance, 
if $2\beta\not\in R$)
then the $\beta$-factor simplifies further to
\[
\frac{\bigl(q_\beta t^{\beta};q_\beta\bigr)_{\infty}}
{\bigl(q_\beta\kappa_\beta^2t^{\beta};q_\beta\bigr)_{\infty}}.
\]
\end{rema}
By the conditions \eqref{smallparameterorth} on the parameters, 
$v$ is a continuous
function on the compact torus 
$T_\Lambda^u=\textup{Hom}(\Lambda,S^1)\subset T_\Lambda$,
where $S^1=\{z\in\mathbb{C}\, | \, |z|=1\}$. Write $dt$ for the
normalized Haar measure on $T_\Lambda^u$.

\begin{defi}
Define a sesquilinear form
$\langle\cdot,\cdot\rangle: \mathbb{C}[T_\Lambda]\times
\mathbb{C}[T_\Lambda]\rightarrow \mathbb{C}$ by
\[
\langle p,r\rangle:=\int_{T_\Lambda^u}p(t)\overline{r(t)}v(t)dt.
\]
\end{defi}

\begin{prop}\label{adjoint}
Let $p,r\in\mathbb{C}[T_\Lambda]$ and $w\in W$. Then
$\langle \pi_{\kappa,q}(T_w)p,r\rangle=\langle p,
\pi_{\kappa^{-1},q^{-1}}(T_w^{-1})r\rangle$.
\end{prop}

The biorthogonality of the nonsymmetric Macdonald-Koornwinder polynomials 
readily follows from Proposition \ref{adjoint}.
\begin{thm}\label{biortho}
If $\lambda,\mu\in \Lambda$ and $\lambda\not=\mu$ then
$\langle P_\lambda,P_\mu^\circ\rangle=0$.
\end{thm}

\subsection{Macdonald operators}\label{Macdonalddiff}
Special cases of what now are known as the 
Macdonald $q$-difference operators were explicitly written down
by Macdonald \cite{Mpol} when introducing the symmetric
Macdonald polynomials. Earlier Ruijsenaars \cite{R} 
had introduced these commuting $q$-difference operators for 
$R_0$ of type $A$ 
as the quantum Hamiltonian of a relativistic version of the quantum
trigonometric Calogero-Moser system (see Subsection \ref{GLn}). 
Koornwinder \cite{Ko} introduced a multivariable extension of 
the second-order Askey-Wilson $q$-difference operator to define
the symmetric Koornwinder polynomials. This case corresponds to
$D=(R_0,\Delta_0,t,\mathbb{Z}R_0,\mathbb{Z}R_0)$ with $R_0$
of type $\textup{A}_1$ or of type $\textup{B}_n$ with $n\geq 2$
(see Subsection \ref{CcheckC}). 

The construction of the whole family of Macdonald $q$-difference
operators using affine Hecke algebras is due to Cherednik \cite{CAnn}
in case $D=(R_0,\Delta_0,\bullet,P(R_0),P(R_0^d))$ and due to Noumi \cite{N}
in case $D=(R_0,\Delta_0,t,\mathbb{Z}R_0,\mathbb{Z}R_0)$ with $R_0$
of type $\textup{A}_1$ or of type $\textup{B}_n$ ($n\geq 2$). We explain
this construction here, see \cite[\S 4.4]{M} for a treatment close
to the present one. 

In this subsection we assume that $\kappa\in\mathcal{M}$ and $q$
satisfy \eqref{smallparameter} or \eqref{largeparameter}.
Consider the linear map 
\[
\textup{Res}_q: W\ltimes_q\mathbb{C}(T_\Lambda)\rightarrow
\tau(\Lambda^d)\ltimes_q\mathbb{C}(T_\Lambda)
\]
defined by
\[
\textup{Res}_q(\sum_{w\in W_0}D_ww):=\sum_{w\in W_0}D_w,
\]
where $D_w\in\tau(\Lambda^d)\ltimes_q\mathbb{C}(T_\Lambda)$ ($w\in W_0$).
Note that
\begin{equation}\label{Resproperty}
L(p)=\bigl(\textup{Res}_q(L)\bigr)(p)
\qquad \forall\, p\in\mathbb{C}(T_\Lambda)^{W_0}.
\end{equation}
\begin{lem}
Let $\beta_{\kappa,q}(H_0)^\prime$ be the commutant of the subalgebra
$\beta_{\kappa,q}(H_0)$ in $W\ltimes_q\mathbb{C}(T_\Lambda)$.
Then $\textup{Res}_q$ restricts to an algebra homomorphism
\[
\textup{Res}_q: \beta_{\kappa,q}(H_0)^\prime\rightarrow 
(\tau(\Lambda^d)\ltimes_q\mathbb{C}(T_\Lambda))^{W_0}
\]
where $(\tau(\Lambda^d)\ltimes_q\mathbb{C}(T_\Lambda))^{W_0}$
is the subalgebra of $W\ltimes_q\mathbb{C}(T_\Lambda)$
consisting of $W_0$-invariant $q$-difference operators.
\end{lem}
Since $Z(H(\kappa^\bullet))=\mathbb{C}_Y[T_{\Lambda^d}]^{W_0}$,
the lemma implies that the $W_0$-invariant $q$-difference operators
\[
D_p:=\textup{Res}_q\bigl(\beta_{\kappa,q}(p(Y))\bigr)
\in\bigl(\tau(\Lambda^d)\ltimes_q\mathbb{C}(T_\Lambda)\bigr)^{W_0},
\qquad p\in\mathbb{C}[T_{\Lambda^d}]^{W_0}
\]
pairwise commute. The operator $D_p$ is called
the Macdonald $q$-difference operator associated to 
$p\in\mathbb{C}[T_{\Lambda^d}]^{W_0}$.

Define the orbit sums $m_\xi^d\in\mathbb{C}[T_{\Lambda^d}]^{W_0}$ 
($\xi\in\Lambda^{d}$) by
\[
m_\xi^d(t):=\sum_{\eta\in W_0\xi}t^\eta.
\]
Then $\{m_\xi^d\}_{\xi\in\Lambda^{d-}}$ is a linear basis of
$\mathbb{C}[T_{\Lambda^d}]^{W_0}$. We write $D_\xi$ for $D_{m_\xi^d}$.

Set, for $w\in W$,
\[
c_w:=\prod_{a\in R^{\bullet+}\cap w^{-1}R^{\bullet-}}c_a\in\mathbb{C}(T_\Lambda)
\]
and $\kappa_w:=\prod_{a\in R^\bullet+\cap w^{-1}R^{\bullet-}}\kappa_a$.
Write $W_{0,\xi}$ for the stabilizer subgroup of $\xi$ in $W_0$. 
\begin{rema}
If $\xi\in \Lambda^{d-}$ and $w\in W_{0,\xi}$ then
$w(c_{\tau(-\xi)})=c_{\tau(-\xi)}$ in $\mathbb{C}(T_\Lambda)$.
\end{rema}
\begin{prop}\label{decomp}
Let $\xi\in\Lambda^{d-}$. Then
\[
D_\xi=\kappa_{\tau(-\xi)}^{-1}\sum_{w\in W_0/W_{0,\xi}}w(c_{\tau(-\xi)})\tau(w\xi)_q+
\sum_{\eta\in\textup{Sat}(\xi_+)\setminus W_0\xi}g_\eta\tau(\eta)_q
\]
for certain $g_\eta\in\mathbb{C}(T_\Lambda)$ satisfying $g_{w\eta}=w(g_\eta)$
for all $w\in W_0$ and $\eta\in\Lambda^d$.
\end{prop}

The set of dominant miniscule weights 
in $\Lambda^d$ is defined by
\[
\Lambda_{min}^{d+}:=
\{\xi\in\Lambda^d\,\, | \,\, \bigl(\xi,\alpha^{d\vee}
\bigr)\in\{0,1\} \quad \forall \alpha\in R_0^+\}.
\]
Set $\Lambda_0^d:=\Lambda^d\cap V^\perp$. We now first
give an explicit description of the dominant miniscule weights.

Recall that $Da_0=-\varphi$ if $\bullet=u$
and $Da_0=-\theta$ if $\bullet=t$. Hence
$-(Da_0)^{d\vee}$ is the highest root of $R_0^{d\vee}$.
Consider the expansion of
$-(Da_0)^{d\vee}\in R_0^{d\vee}$ with respect to
the ordered basis $\Delta_0^{d\vee}=(\alpha_1^{d\vee},\ldots,\alpha_n^{d\vee})$ 
of $R_0^{d\vee}$,
\[
-(Da_0)^{d\vee}=\sum_{i=1}^nm_i\alpha_i^{d\vee},
\]
then $m_i\in\mathbb{Z}_{\geq 1}$ for all $i$. Set
\[
J_{\Lambda^d}^+:=\{i\in\{1,\ldots,n\} \,\, | \,\, m_i=1\,\, 
\&\,\, (\varpi_i^d+ V^\perp)\cap \Lambda^d\not=\emptyset\},
\] 
where $V^\perp$ is the orthocomplement of
$V$ in $Z$. 
\begin{prop}\label{miniscule}
{\bf (i)} $\Lambda_{min}^{d+}$ is a complete set of representatives
of $\Lambda^d/\mathbb{Z}R_0^d$.\\
{\bf (ii)} For $j\in J^+_{\Lambda^d}$ choose an element 
$\widetilde{\varpi}_j^d\in
(\varpi_j^d+V^\perp)\cap \Lambda^d$. Then
\[
\Lambda_{min}^{d+}=\Lambda_0^d\cup\bigcup_{j\in J^+_{\Lambda^d}}
\bigl(\widetilde{\varpi}_j^d+\Lambda_0^d\bigr)
\]
(disjoint union).\\
{\bf (iii)} For $\eta\in\Lambda^d$ let
$u(\eta)\in W$ be the unique element of minimal length (with respect to
$l$) in the coset $\tau(\eta)W_0$. Then
$\Omega=\{u(\xi)\,\, | \,\, \xi\in \Lambda_{min}^{d+}\}$. 
\end{prop}
Since $-(Da_0)^{d\vee}\in R_0^{d\vee+}$ is the highest root,
$-(Da_0)^d\in R_0^{d+}$ is quasi-miniscule, i.e. $(-(Da_0)^d,\alpha^{d\vee})
\in\{0,1\}$ for all $\alpha^d\in R_0^{d+}\setminus\{-(Da_0)^d\}$.

Let $w_0\in W_0$ be the longest Weyl group element.
\begin{cor}\label{explicitD}
{\bf (i)} For $j\in J_{\Lambda^d}^+$ we have
\[
D_{w_0\widetilde{\varpi}_j^d}=
\kappa_{\tau(-w_0\widetilde{\varpi}_j^d)}^{-1}
\sum_{w\in W_0/W_{0,w_0\widetilde{\varpi}_j^d}}
w(c_{\tau(-w_0\widetilde{\varpi}_j^d)})
\tau(ww_0\widetilde{\varpi}_j^d)_q.
\]
{\bf (ii)} We have
\[
D_{(Da_0)^d}=\kappa_{\tau(-(Da_0)^d)}^{-1}\sum_{w\in W_0/W_{0,(Da_0)^d}}
w(c_{\tau(-(Da_0)^d)})\bigl(\tau(w(Da_0)^d)_q-1\bigr)+
m_{(Da_0)^d}^d(\gamma_{0,q}^{-1}).
\]
\end{cor}
\subsection{Monic symmetric Macdonald-Koornwinder  polynomials}
Historically the symmetric Macdonald-Koornwinder polynomials preceed
the nonsymmetric Mac\-do\-nald-Koorn\-win\-der polynomials.
The monic symmetric Macdonald polynomials associated to initial
data of the form $D=(R_0,\Delta_0,\bullet,P(R_0),P(R_0^d))$
were defined in Macdonald's handwritten preprint from 1987.
It appeared in print in 2000, see \cite{Mpol} (see \cite{Mbook} for 
the $\textup{GL}_n$ case). Macdonald used the fact that
for a fixed $\xi\in\{w_0\widetilde{\varpi}_j^d\}_{j\in J_{\Lambda^d}^+}\cup
\{(Da_0)^d\}$, the explicit Macdonald $q$-difference operators
$D_\xi$ (see Corollary \ref{explicitD}) is a linear operator 
on $\mathbb{C}[T_{\Lambda}]^{W_0}$ which is triangular with respect to
the suitable partially ordered basis of orbit sums and has
(generically) simple spectrum. 

This approach was extended by Koornwinder \cite{Ko} to the case
corresponding to the initial data
$D=(R_0,\Delta,t,\mathbb{Z}R_0,\mathbb{Z}R_0)$ with $R_0$ of type $\textup{A}_1$
or of type $\textup{B}_n$ ($n\geq 2$), in which case 
$D_{-\theta}$ is Koornwinder's multivariable extension
of the Askey-Wilson second-order $q$-difference operator (see Subsection
\ref{CcheckC}). The corresponding symmetric Macdonald-Koornwinder polynomials 
are the Askey-Wilson \cite{AW}
polynomials if $R_0$ is of rank one and the symmetric
Koornwinder \cite{Ko} polynomials if $R_0$ is of higher rank. 

In this subsection we introduce the monic symmetric Macdonald-Koornwinder 
polynomials by symmetrizing the nonsymmetric ones, cf. \cite[\S 4]{CAnn}.
We assume throughout this subsection that $\kappa\in\mathcal{M}$ 
and $q$ satisfy \eqref{smallparameter} or \eqref{largeparameter}.
Recall the notation $\kappa_w:=\prod_{a\in R^{\bullet+}\cap w^{-1}R^{\bullet-}}
\kappa_a$. It only depends on $\kappa^\bullet=\kappa|_{R^\bullet}$.
It satisfies $\kappa_v\kappa_w=\kappa_{vw}$ 
if $l(vw)=l(v)+l(w)$.
Hence there exists a unique linear character 
$\chi_+: H(\kappa^\bullet)
\rightarrow\mathbb{C}$ satisfying
$\chi_+(T_w)=\kappa_w$ for all $w\in W$, the trivial 
linear character of $H(\kappa^\bullet)$.

Define $C_+\in H_0(\kappa^\bullet|_{R_0})\subset
H(\kappa^\bullet)$ by
\begin{equation}\label{Cplus}
C_+=\frac{1}{\sum_{w\in W_0}\kappa_w^{2}}\sum_{w\in W_0}\kappa_wT_w.
\end{equation}
The normalization is such that $\chi_+(C_+)=1$. Then 
$T_iC_+=\kappa_{i}C_+=C_+T_i$ for $1\leq i\leq n$ and
$C_+^2=C_+$. 
The following lemma follows from the explicit expression
of $\pi_{\kappa,q}(T_i)$ ($1\leq i\leq n$).
\begin{lem}
The linear endomorphism $\pi_{\kappa,q}(C_+)$ of $\mathbb{C}[T_\Lambda]$ 
is an idempotent with image $\mathbb{C}[T_\Lambda]^{W_0}$.
\end{lem}
Consider the linear basis $\{m_\lambda\}_{\lambda\in \Lambda^+}$ of 
$\mathbb{C}[T_\Lambda]^{W_0}$ given by orbit sums, $m_\lambda(t)=
\sum_{\mu\in W_0\lambda}t^{\mu}$. Recall that 
$P_\lambda=P_\lambda(D;\kappa,q)$ denotes the monic nonsymmetric 
Macdonald-Koornwinder polynomial of degree $\lambda\in\Lambda$.

\begin{lem}
If $\lambda\in \Lambda^+$ then
$\pi_{\kappa,q}(C_+)P_\lambda=
\sum_{\mu\in \Lambda^+: \mu\leq\lambda}c_{\lambda,\mu}m_\mu$ 
for certain $c_{\lambda,\mu}\in\mathbb{C}$ with
$c_{\lambda,\lambda}\not=0$. 
\end{lem}
\begin{defi}
The monic symmetric Macdonald-Koornwinder 
polynomial $P_\lambda^+=P_\lambda^+(D;\kappa,q)$
of degree $\lambda\in\Lambda^+$ is defined by
\[
P_\lambda^+:=c_{\lambda,\lambda}^{-1}\pi_{\kappa,q}(C_+)P_\lambda
\in\mathbb{C}[T_\Lambda]^{W_0}.
\]
\end{defi}
\begin{thm}
The monic symmetric Macdonald-Koornwinder polynomial
$P_\lambda^+$ of degree $\lambda\in\Lambda^+$
is the unique element in $\mathbb{C}[T_\Lambda]^{W_0}$
satisfying
\begin{enumerate}
\item[{\bf (i)}] $P_\lambda^+=\sum_{\mu\in\Lambda^+: \mu\leq\lambda}
d_{\lambda,\mu}m_{\mu}$ with $d_{\lambda,\mu}\in\mathbb{C}$ and
$d_{\lambda,\lambda}=1$,
\item[{\bf (ii)}] $D_pP_\lambda^+=
p\bigl(q^{-\lambda}\prod_{\alpha\in R_0^+}\upsilon_\alpha^{-\alpha^{d\vee}}\bigr)
P_\lambda^+$ for all $p\in\mathbb{C}[T_{\Lambda^d}]^{W_0}$.
\end{enumerate}
(Note that for $p\in\mathbb{C}[T_{\Lambda^d}]^{W_0}$ and $\lambda\in\Lambda^+$
we have
$p\bigl(q^{-\lambda}\prod_{\alpha\in R_0^+}
\upsilon_\alpha^{-\alpha^{d\vee}}\bigr)=
p(\gamma_{w_0\lambda,q}^{-1})$).
\end{thm}

Explicit expressions of the monic symmetric Macdonald-Koornwinder 
polynomials when $R_0$ is of rank one are given in Subsection \ref{GLn}
and Subsection \ref{CcheckC}.

\subsection{Orthogonality}
In this subsection we assume that
$\kappa\in\mathcal{M}$ and $q$ 
satisfy the conditions \eqref{smallparameterorth}. We thus have
the monic symmetric Macdonald-Koornwinder 
polynomials $\{P_\lambda^+\}_{\lambda\in\Lambda^+}$
with respect to the parameters $(\kappa,q)$, as well as 
the monic symmetric Macdonald-Koornwinder polynomials with
respect to the parameters $(\kappa^{-1},q^{-1})$, in which case we denote
them by $\{P_\lambda^{\circ+}\}_{\lambda\in\Lambda^+}$.

Define the $W_0$-invariant meromorphic function
\[
v_+:=\prod_{a\in R^\bullet: a(0)\geq 0}\frac{1}{c_a}
\]
on $T_\Lambda$. 
It relates to the weight function $v$ by
\begin{equation}\label{mathcalC}
v=\mathcal{C}v_+,\qquad
\mathcal{C}=\mathcal{C}(\cdot;D;\kappa,q):=\prod_{\alpha\in R_0^{-}}c_\alpha.
\end{equation}
One recovers $v_+$ from $v$ up to a multiplicative constant 
by symmetrization in view of the following property of
the rational function $\mathcal{C}\in\mathbb{C}(T_\Lambda)$
(cf. \cite[(2.8)\& (2.8 nr)]{MacPoincare}). 

\begin{lem}\label{curlC}
We have
\begin{equation}\label{formulaK}
\sum_{w\in W_0}w\mathcal{C}=\mathcal{C}(\gamma_{0,q}^d)
\end{equation}
as identity in $\mathbb{C}(T_\Lambda)$, where
$\gamma_{\xi,q}^d:=\gamma_{\xi,q}(D^d;\kappa^{d\bullet})\in
T_{\Lambda}$ ($\xi\in\Lambda^d$), cf. Definition \ref{gammapoint}.
More generally, if $\xi\in\Lambda^{d-}$ and if $(\kappa,q)$ is generic
then
\[
\sum_{w\in W_0}w\mathcal{C}=\mathcal{C}(\gamma_{\xi,q}^d)
\sum_{\eta\in W_0\xi}\prod_{\alpha\in R_0^+\cap v(\eta)R_0^-}
\frac{c_\alpha(\gamma_{\xi,q}^d)}{c_{-\alpha}(\gamma_{\xi,q}^d)}
\]
as identity in $\mathbb{C}(T_\Lambda)$, where (recall)
$v(\eta)\in W_0$ is the
element of shortest length such that $v(\eta)\eta=\eta_-$.
\end{lem}

The meromorphic function $v_+$ on $T_{\Lambda}$
reads in terms of $q$-shifted factorials,
\begin{equation}\label{weightsymmetric}
v_+(t)=
\prod_{\beta\in R_0}
\frac{\bigl(t^{2\beta};q_\beta^{2}\bigr)_{\infty}}
{\bigl(\kappa_\beta\kappa_{2\beta}t^\beta,-\kappa_\beta\kappa_{2\beta}^{-1}t^\beta,
q_\beta\kappa_{\mu_\beta c+\beta}
\kappa_{2\mu_\beta c+2\beta}t^\beta,
-q_\beta\kappa_{\mu_\beta c+\beta}
\kappa_{2\mu_\beta c+2\beta}^{-1}t^\beta;
q_\beta^{2}\bigr)_{\infty}}.
\end{equation}
It is a nonnegative real-valued continuous function on $T_\Lambda^u$
(it is nonnegative since it can be written on $T_\Lambda^u$ as
$v_+(t)=|\delta(t)|^2$ with $\delta(t)$ the expression 
\eqref{weightsymmetric} with product taken 
only over the set $R_0^{+}$ of positive roots). 
\begin{cor}
$\langle\cdot,\cdot\rangle$ restricts to a positive definite,
sesquilinear form on $\mathbb{C}[T_\Lambda]^{W_0}$. In fact, for 
$p,r\in\mathbb{C}[T_\Lambda]^{W_0}$,
\[
\langle p,r\rangle=\frac{\mathcal{C}(\gamma_{0,q}^d)}{\#W_0}\int_{T_\Lambda^u}p(t)
\overline{r(t)}v_+(t)dt.
\]
\end{cor}
Symmetrizing the results on the monic nonsymmetric
Macdonald-Koornwinder polynomials using the idempotent $C_+\in H_0$ gives
the following properties of the monic symmetric Macdonald-Koornwinder 
polynomials. 
\begin{thm} Let $\lambda\in \Lambda^+$.\\
{\bf (a)}
The symmetric Macdonald-Koornwinder 
polynomial $P_\lambda^+\in\mathbb{C}[T_\Lambda]^{W_0}$
satisfies the following characterizing properties,
\begin{enumerate}
\item[{\bf (i)}] $P_\lambda^+=\sum_{\mu\in \Lambda^+: \mu\leq\lambda}
d_{\lambda,\mu}m_{\mu}$ with $d_{\lambda,\lambda}=1$,
\item[{\bf (ii)}] $\langle P_\lambda^+,m_\mu\rangle=0$ if $\mu\in \Lambda^+$
and $\mu<\lambda$.
\end{enumerate}
{\bf (b)} $P_\lambda^{\circ+}=P_\lambda^+$.\\
{\bf (c)} $\langle P_\lambda^+,P_\mu^+\rangle=0$ if $\mu\in \Lambda^+$
and $\mu\not=\lambda$.
\end{thm}

\subsection{$\textup{GL}_{n+1}$ Macdonald polynomials}\label{GLn}
Take $n\geq 1$ and 
\[V:=(\epsilon_1+\cdots+\epsilon_{n+1})^\perp\subset
\mathbb{R}^{n+1}=:Z\]
with $\{\epsilon_i\}_{i=1}^{n+1}$ the standard orthonormal basis of
$\mathbb{R}^{n+1}$. Let
$R_0=\{\epsilon_i-\epsilon_j\}_{1\leq i\not=j\leq n+1}$ with
ordered basis 
\[
\Delta_0=(\alpha_1,\ldots,\alpha_n)=(\epsilon_1-\epsilon_2,\ldots,
\epsilon_n-\epsilon_{n+1}).
\]
As lattices take $\Lambda=\mathbb{Z}^{n+1}=\Lambda^d$.
Then indeed
$D:=(R_0,\Delta_0,u,\mathbb{Z}^{n+1},\mathbb{Z}^{n+1})\in\mathcal{D}$.

The corresponding irreducible reduced affine root system is 
$R^u=\{mc+\alpha\}_{m\in\mathbb{Z}, \alpha\in R_0}$, and the corresponding
additional simple affine root is $a_0=c-\epsilon_1+\epsilon_{n+1}$.
There is no nonreduced extension of $R^u$ involved
since $\bigl(\Lambda,a_i^\vee\bigr)=\mathbb{Z}$
for $i\in\{0,\ldots,n\}$. Hence $R=R^u$.

The fundamental weights $\varpi_j=\varpi_j^u$ ($1\leq j\leq n$)
are given by
\[
\varpi_j=-\frac{j}{n+1}(\epsilon_1+\cdots+\epsilon_{n+1})+
\epsilon_1+\epsilon_2+\cdots+\epsilon_j.
\]
Define for $1\leq j\leq n+1$ the elements 
\[
\widetilde{\varpi}_j:=\epsilon_1+\cdots+\epsilon_j.
\]
The orthocomplement $V^\perp$ of $V$ in $Z$ is 
$\mathbb{R}\widetilde{\varpi}_{n+1}$ and
$\Lambda_0^d:=\Lambda^d\cap V^\perp=\mathbb{Z}\widetilde{\varpi}_{n+1}$.
Then
\[
\widetilde{\varpi}_j\in\bigl(\varpi_j+V^\perp\bigr)\cap\Lambda^d
\qquad \forall\,
j\in\{1,\ldots,n\}.
\]
Since $-(Da_0)^{d\vee}=\epsilon_1-\epsilon_{n+1}=
\sum_{i=1}^n\alpha_i$ we conclude that
\[
J_{\Lambda^d}^+=\{1,\ldots,n\}.
\]
The miniscule dominant weights in $\Lambda^d$ are thus given by
\[
\Lambda_{min}^{d+}=\Lambda_0^d\cup\bigcup_{j=1}^n(\widetilde{\varpi}_j+
\Lambda_0^d).
\]
We have $u(\epsilon_1)=\tau(\epsilon_1)s_1\cdots s_{n-1}s_n$ in the extended
affine Weyl group $W\simeq S_{n+1}\ltimes\mathbb{Z}^{n+1}$ (where $S_{n+1}$
is the symmetric group in $n+1$ letters). Note also that $u(\epsilon_1)(a_j)=
a_{j+1}$ for $0\leq j<n$ and $u(\epsilon_1)(a_n)=a_0$. 
In addition
we have $u(\epsilon_1)^j=u(\widetilde{\varpi}_j)$ for $1\leq j\leq n$
and $u(\epsilon_1)^{n+1}=u(\widetilde{\varpi}_{n+1})=
\tau(\widetilde{\varpi}_{n+1})$ in $W$. Hence 
$\mathbb{Z}\simeq \Omega$ by $m\mapsto u(\epsilon_1)^m$.

Observe that $R=R^u$ has one $W(R^u)$-orbit, hence also one $W$-orbit.
The affine Hecke algebra $H(W(R);\kappa)$ and the extended affine Hecke
algebra $H(\kappa)=H(D;\kappa)$ thus
depends on a single nonzero complex number $\kappa$. 
$\mathbb{Z}$ acts on $H(W(R);\kappa)$ by algebra automorphisms,
where $1\in\mathbb{Z}$ acts on $T_i$ 
by mapping it to $T_{i+1}$ (reading the subscript modulo $n+1$). We
write $\mathbb{Z}\ltimes H(W(R);\kappa)$ for the associated semidirect
product algebra. 
\begin{prop}
{\bf (i)} 
$\mathbb{Z}\ltimes H(W(R);\kappa)\simeq H(\kappa)$
by mapping the generator $1\in\mathbb{Z}$ to 
$u(\widetilde{\varpi}_1)=u(\epsilon_1)$.\\
{\bf (ii)} For $1\leq i\leq n+1$ we have
\[
Y^{\epsilon_i}=T_{i-1}^{-1}\cdots T_2^{-1}T_1^{-1}u(\epsilon_1)
T_{n}T_{n-1}\cdots T_i.
\]
in $H(\kappa)$.
\end{prop}

We now turn to the explicit description of the Macdonald-Ruijsenaars
$q$-difference operators and the orthogonality measure for the 
$\textup{GL}_{n+1}$ Macdonald polynomials.

The longest Weyl group element $w_0\in W_0$ maps $\epsilon_i$ to
$\epsilon_{n+2-i}$ for $1\leq i\leq n+1$, hence 
\[
R^+\cap \tau(w_0\widetilde{\varpi}_j)R^-=\{\epsilon_r-\epsilon_s\, | \, 
1\leq r\leq n+1-j\,\,\&\,\, n+2-j\leq s\leq n+1 \}.
\]
Write $t_i=t^{\epsilon_i}$ for $t\in T_{\Lambda}$ and
$1\leq i\leq n+1$. Then, for $1\leq j\leq n$ we obtain
\[
c_{\tau(-w_0\widetilde{\varpi}_j)}(t)=
\prod_{\stackrel{1\leq j\leq n+1-i}{n+2-j\leq s\leq n+1}}
\frac{(1-\kappa^2t_rt_s^{-1})}{(1-t_rt_s^{-1})}
\]
and consequently
\begin{equation*}
\begin{split}
D_{w_0\widetilde{\varpi}_j}&=
\sum_{\stackrel{I\subset\{1,\ldots,n+1\}:}{\#I=n+1-j}}
\left(\prod_{r\in I, s\not\in I}\frac{(1-\kappa^2t_rt_s^{-1})}
{\kappa(1-t_rt_s^{-1})}
\right)\tau\bigl(\sum_{s\not\in I}\epsilon_s\bigr)_q\\
&=\sum_{\stackrel{I\subset\{1,\ldots,n+1\}:}{\#I=j}}
\left(\prod_{r\in I, s\not\in I}\frac{\kappa^{-1}t_r-\kappa t_s}{t_r-t_s}\right)
\tau\bigl(\sum_{r\in I}\epsilon_r\bigr)_q
\end{split}
\end{equation*}
for $1\leq j\leq n$ by Proposition \ref{explicitD}{\bf (i)}.
These commutative 
$q$-difference operators were introduced by Ruijsenaars \cite{R}
as the quantum Hamiltonians of a relativistic version of the trigonometric
quantum Calogero-Moser system.
These operators also go by the name (trigonometric) Ruijsenaars
$q$-difference operators, or Ruijsenaars-Macdonald $q$-difference
operators.

The corresponding monic
symmetric Macdonald polynomials $\{P_\lambda\}_{\lambda\in\Lambda^+}$
are parametrized by
\[
\Lambda^+=\{\lambda\in\mathbb{Z}^{n+1} \,\, | \,\, \lambda_1\geq\lambda_2\geq
\cdots\geq\lambda_{n+1}\}.
\]
The orthogonality weight function then becomes
\[
v_+(t)=\prod_{1\leq i\not=j\leq n+1}
\frac{\bigl(t_i/t_j;q\bigr)_{\infty}}{\bigl(\kappa^2t_i/t_j;q\bigr)_{\infty}}.
\]
The $\textup{GL}_2$ Macdonald polynomials are essentially the
continuous $q$-ultraspherical polynomials (see
\cite[Chpt. 2]{C} or \cite[\S 6.3]{M}),
\[
P_\lambda^+(t)=
t_1^{\lambda_1}t_2^{\lambda_2}{}_2\phi_1\left(\begin{matrix} \kappa^2, 
q^{\lambda_2-\lambda_1}\\ q^{1+\lambda_2-\lambda_1}/\kappa^2 \end{matrix};
q,\frac{qt_2}{\kappa^2 t_1}\right),
\] 
for $\lambda=(\lambda_1,\lambda_2)\in\mathbb{Z}^2$ with 
$\lambda_1\geq\lambda_2$, where we use 
standard notations for the basic hypergeometric ${}_r\phi_s$ series
(see, e.g., \cite{GR}).

\subsection{Koornwinder polynomials}\label{CcheckC}
Take $Z=V=\mathbb{R}^n$ with standard orthonormal basis $\{\epsilon_i\}_{i=1}^n$.
We realize the root system $R_0\subset V$ of type $\textup{B}_n$
as
\[R_0=\{\pm\epsilon_i\}_{i=1}^n\cup
\{\pm\epsilon_i\pm\epsilon_j\}_{1\leq i<j\leq n}
\]
(all sign combinations allowed). 
The $W_0$-orbits are 
$\mathcal{O}_l=\{\pm\epsilon_i\pm\epsilon_j\}_{1\leq i<j\leq n}$ and
$\mathcal{O}_s=\{\pm\epsilon_i\}_{i=1}^n$.
As an ordered basis of $R_0$ take
\[
\Delta_0=(\alpha_1,\ldots,\alpha_{n-1},\alpha_n)=
(\epsilon_1-\epsilon_2,\ldots,\epsilon_{n-1}-\epsilon_n,\epsilon_n).
\]
Note that $\mathbb{Z}R_0=\mathbb{Z}^n$.
For $n=1$ this should be interpreted as
$R_0=\{\pm\epsilon_1\}$, the root system of type $\textup{A}_1$,
with basis element $\epsilon_1$.
We consider in this subsection the initial data
$D=(R_0,\Delta_0,t,\mathbb{Z}^n,\mathbb{Z}^n)\in\mathcal{D}$.

The associated reduced affine root system is
\[
R^t=\{\frac{m}{2}c+\alpha\}_{n\in\mathbb{Z},\alpha\in\mathcal{O}_s}\cup
\{mc+\beta\}_{m\in\mathbb{Z}, \beta\in\mathcal{O}_l}
\]
(for $n=1$ it should be read as $R^t=\{\frac{m}{2}c\pm\epsilon_1\}_{m\in
\mathbb{Z}}$, the affine root system $R^t$ with $R_0$ of type $\textup{A}_1$).
The associated ordered basis of $R^t$ is
\[
\Delta=(a_0,a_1,\ldots,a_n)=(\frac{1}{2}c-\theta,\alpha_1,\ldots,\alpha_n)
\]
with $\theta=\epsilon_1=\sum_{j=1}^n\alpha_j$ the highest short root of
$R_0$ with respect to $\Delta_0$. The affine root system $R^t$ has three
$W$-orbits,
\[
R^t=\widehat{\mathcal{O}}_1\cup\widehat{\mathcal{O}}_2\cup
\widehat{\mathcal{O}}_3,
\]
where 
$\widehat{\mathcal{O}}_1=
Wa_0=\bigl(\frac{1}{2}+\mathbb{Z}\bigr)c+\mathcal{O}_s$,
$\widehat{\mathcal{O}}_2=Wa_i=\mathbb{Z}c+\mathcal{O}_l$ 
(where $1\leq i<n$) and
$\widehat{\mathcal{O}}_3=Wa_n=\mathbb{Z}c+\mathcal{O}_s$.
If $n=1$ then 
$\Delta=(a_0,a_1)=(\frac{1}{2}c-\epsilon_1,\epsilon_1)$
and $R^t$ has two $W$-orbits $\widehat{\mathcal{O}}_1$
and $\widehat{\mathcal{O}}_3$.

Note that 
$\bigl(\Lambda,a_i^\vee\bigr)=2\mathbb{Z}$ for $i=0$ and $i=n$,
hence $S=\{0,n\}=S^d$ and
\[
R=R^t\cup W(2a_0)\cup W(2a_n)=
\widehat{\mathcal{O}}_1\cup\widehat{\mathcal{O}}_2\cup
\widehat{\mathcal{O}}_3\cup\widehat{\mathcal{O}}_4\cup
\widehat{\mathcal{O}}_5
\]
with additional $W$-orbits $\widehat{\mathcal{O}}_4=
2\widehat{\mathcal{O}}_1$ and $\widehat{\mathcal{O}}_5=
2\widehat{\mathcal{O}}_3$.
If $n=1$ then $R$ has four $W$-orbits 
$\widehat{\mathcal{O}}_i$ ($i=1,2,4,5$).
Since $D^d=D$ we have $R^d=R$, $\Delta^d=\Delta$
and $W^d=W$.

Suppose that $\kappa\in\mathcal{M}(D)$ and $q\in\mathbb{C}^*$
satisfy \eqref{smallparameter} or \eqref{largeparameter}. 
Fix $1\leq i<n$. 
Then $\kappa\in\mathcal{M}(D)$ is determined by five (four in
case $n=1$) independent numbers $\kappa_{a_0}, \kappa_{a_i}, \kappa_{a_n}, 
\kappa_{2a_0}$ and $\kappa_{2a_n}$. The corresponding Askey-Wilson
\cite{AW} parameters are defined by
\begin{equation}\label{parameterKoornwinder}
(a,b,c,d,k)=
(\kappa_{a_n}\kappa_{2a_n},-\kappa_{a_n}\kappa_{2a_n}^{-1},
q^{\frac{1}{2}}\kappa_{a_0}\kappa_{2a_0},
-q^{\frac{1}{2}}\kappa_{a_0}\kappa_{2a_0}^{-1},\kappa_{a_i}^2)
\end{equation}
(the parameter $k$ is dropping out in case $n=1$).
The dual multiplicity function $\kappa^d$ on $R^d=R$
is then determined by $\kappa^d_{a_0}:=\kappa_{2a_n}$, $\kappa_{a_i}^d:=
\kappa_{a_i}$, $\kappa_{a_n}^d:=\kappa_{a_n}$, $\kappa_{2a_0}^d:=\kappa_{2a_0}$
and $\kappa_{2a_n}^d:=\kappa_{a_0}$. The corresponding Askey-Wilson parameters
are
\[
(\tilde{a},\tilde{b},\tilde{c},\tilde{d},\tilde{k})=
(\kappa_{a_n}\kappa_{a_0},-\kappa_{a_n}\kappa_{a_0}^{-1},
q^{\frac{1}{2}}\kappa_{2a_n}\kappa_{2a_0},-q^{\frac{1}{2}}\kappa_{2a_n}
\kappa_{2a_0}^{-1},\kappa_{a_i}^2).
\]
In terms of the original Askey-Wilson parameters this can be
expressed as $\tilde{k}=k$ and
\[
(\tilde{a},\tilde{b},\tilde{c},\tilde{d})=
\Bigl(\sqrt{q^{-1}abcd},\frac{ab}{\sqrt{q^{-1}abcd}},
\frac{ac}{\sqrt{q^{-1}abcd}},\frac{ad}{\sqrt{q^{-1}abcd}}\Bigr).
\]

Note that 
\[
-(Da_0)^{d\vee}=\theta^\vee=2\alpha_1^\vee+2\alpha_2^\vee+
\cdots+2\alpha_{n-1}^\vee+\alpha_n^\vee.
\]
Furthermore, denoting the 
fundamental weights of $R_0$ with respect to the ordered
basis $\Delta_0$ by $\{\varpi_i\}_{i=1}^n$, we have 
\[
\varpi_n=\frac{1}{2}(\epsilon_1+\cdots+\epsilon_n)\not\in\Lambda^d=\mathbb{Z}^n.
\]
Consequently $J_{\Lambda^d}^+=\emptyset$. Hence
the only explicit Macdonald $q$-difference operator obtainable 
from Corollary \ref{explicitD} is $D_{-\theta}=D_{-\epsilon_1}$.

Write
\[
\mathcal{D}=\kappa_{\tau(\epsilon_1)}^{-1}
\sum_{w\in W_0/W_{0,\epsilon_1}}w\bigl(c_{\tau(\epsilon_1)}\bigr)
\bigl(\tau(-w\epsilon_1)_q-1\bigr),
\]
so that $D_{-\epsilon_1}=\mathcal{D}+
m_{-\epsilon_1}^d(\gamma_{0,q}^{-1})$.
Write $t_i=t^{\epsilon_i}$ for
$t\in T_{\Lambda}$ and $1\leq i\leq n$. 
Since
\[
R^{t+}\cap \tau(-\epsilon_1)R^{t-}=\{\epsilon_1,\frac{1}{2}c+\epsilon_1\}\cup
\{\epsilon_1\pm\epsilon_j\}_{j=2}^n
\]
($=\{\epsilon_1,\frac{1}{2}c+\epsilon_1\}$ if $n=1$)
it follows that 
\begin{equation}\label{factor}
\begin{split}
\kappa_{\tau(\epsilon_1)}^{-1}
c_{\tau(\epsilon_1)}(t)&=\kappa_{\tau(\epsilon_1)}^{-1}
c_{\epsilon_1}(t)c_{\frac{1}{2}c+\epsilon_1}(t)
\prod_{j=2}^nc_{\epsilon_1-\epsilon_j}(t)c_{\epsilon_1+\epsilon_j}(t)\\
&=\frac{1}{\sqrt{q^{-1}abcd}k^{n-1}}
\frac{(1-at_1)(1-bt_1)(1-ct_1)(1-dt_1)}{(1-t_1^2)(1-qt_1^2)}\\
&\times\prod_{j=2}^n\frac{(1-kt_1t_j^{-1})(1-kt_1t_j)}
{(1-t_1t_j^{-1})(1-t_1t_j)}.
\end{split}
\end{equation}
For $n=1$ the product over $j$ is
not present in \eqref{factor}. In particular, the operator 
$\mathcal{D}$ then only depends on $q$ and on the four parameters $a,b,c,d$. 

Hence $\mathcal{D}$ is 
Koornwinder's \cite{Ko} second order $q$-difference operator
\begin{equation*}
\begin{split}
\mathcal{D}&=\frac{1}{\sqrt{q^{-1}abcd}k^{n-1}}
\sum_{i=1}^n\bigl(A_i(t)(\tau(-\epsilon_i)_q-1)+
A_i(t^{-1})(\tau(\epsilon_i)_q-1),\\
A_i(t)&=\frac{(1-at_i)(1-bt_i)(1-ct_i)(1-dt_i)}{(1-t_i^2)(1-qt_i^2)}
\prod_{j\not=i}\frac{(1-kt_it_j)(1-kt_it_j^{-1})}
{(1-t_it_j)(1-t_it_j^{-1})},
\end{split}
\end{equation*}
with the obvious adjustment for $n=1$, in which case 
$\mathcal{D}$ is the Askey-Wilson \cite{AW} second-order 
$q$-difference operator (this derivation is due to Noumi \cite{N}). 
The weight function $v_+(t)$ becomes 
\begin{equation*}
v_+(t)=\prod_{i=1}^n\frac{\bigl(t_i^2,t_i^{-2};q\bigr)_{\infty}}
{\bigl(at_i,at_i^{-1},
bt_i,bt_i^{-1},ct_i,ct_i^{-1},dt_i,dt_i^{-1};q\bigr)_{\infty}}
\prod_{1\leq r\not=s\leq n}
\frac{\bigl(t_rt_s,t_rt_s^{-1};q\bigr)_{\infty}}
{\bigl(kt_rt_s,kt_rt_s^{-1};q\bigr)_{\infty}}
\end{equation*}

The dominant elements $\Lambda^+$ in $\Lambda$ are the 
partitions of length $\leq n$,
\[
\Lambda^+=\{\lambda\in\mathbb{Z}^n\,\, | \,\, \lambda_1\geq\lambda_2\geq
\cdots\geq\lambda_n\geq 0\}.
\]
The corresponding monic symmetric Macdonald-Koornwinder polynomials 
$\{P_\lambda^+\}_{\lambda\in\Lambda^+}$ are the symmetric Koornwinder
\cite{Ko} polynomials. 

For $n=1$ the Koornwinder polynomials are the Askey-Wilson \cite{AW}
polynomials.
Concretely, 
using standard notations for basic
hypergeometric series (cf. \cite{GR}),
for $m\in\mathbb{Z}_{\geq 0}=\Lambda^+$,
\[
P_m^+(t)=\textup{cst}_m{}_4\phi_3\left(
\begin{matrix} q^{-m},q^{m-1}abcd,at,at^{-1}\\
ab,ac,ad\end{matrix}; q,q\right),
\]
where
\[
\textup{cst}_m=\frac{\bigl(ab,ac,ad;q\bigr)_m}
{a^m\bigl(q^{m-1}abcd;q\bigr)_m}.
\]
The Cherednik-Macdonald theory associated to the Askey-Wilson
polynomials was worked out in detail in \cite{NS}. 
\begin{eg}\label{surprise2}
Consider the initial data $D=(R_0,\Delta_0,t,\mathbb{Z}^n,P(R_0))
\in\mathcal{D}$ 
with the root system $R_0\subset\mathbb{R}^n$ of type $\textup{B}_n$ ($n\geq 2$)
and the ordered basis $\Delta_0$ both defined in terms of the standard
orthonormal basis $\{\epsilon_i\}_{i=1}^n$ of $\mathbb{R}^n$ 
as above. Compared to the Koornwinder setup just discussed
we have thus chosen a different lattice $\Lambda^d=P(R_0)$,
which contains $\mathbb{Z}R_0=\mathbb{Z}^n$ as index two sublattice.
Still $S(D)=\{0,n\}$, hence $R$ is the nonreduced affine
root system with five $W_0\ltimes\mathbb{Z}^n$-orbits 
$\widehat{\mathcal{O}}_i$ ($1\leq i\leq 5$) as introduced above. 
But the number $\nu=\nu(D)$ of 
$W_0\ltimes P(R_0)$-orbits is three: they are given by
$\widehat{\mathcal{O}}_2$, 
$\widehat{\mathcal{O}}_1\cup\widehat{\mathcal{O}}_3$ and 
$\widehat{\mathcal{O}}_4\cup\widehat{\mathcal{O}}_5$. 

On the other hand, $R^d=R^t$ since $S(D^d)=\emptyset$, 
in particular $R^d$ is reduced. It has the three
$W^d=W_0\ltimes\mathbb{Z}^n$-orbits $\widehat{\mathcal{O}}_i$ ($i=1,2,3$).
\end{eg}
\begin{rema}\label{Bnext}
Consider the initial data $D=(R_0,\Delta_0,u,\mathbb{Z}^n,\mathbb{Z}^n)$
with $(R_0,\Delta_0)$ of type $\textup{B}_n$ ($n\geq 3$) as above.
Compared to initial data $D^{C^\vee C_n}:=(R_0,\Delta_0,t,\mathbb{Z}^n,
\mathbb{Z}^n)$ related to Koornwinder polynomials, 
we thus have only changed the type from twisted to untwisted.  
Then $R^u$ is the affine root subsystem $\widehat{\mathcal{O}}_2\cup
\widehat{\mathcal{O}}_3$ of the affine root
system $R^{C^\vee C_n}$ of type
$C^\vee C_n$ as defined above, and
$R=R(D)$ is the nonreduced
irreducible affine root subsystem
$\widehat{\mathcal{O}}_2\cup\widehat{\mathcal{O}}_3\cup\widehat{\mathcal{O}}_5$
of $R(D^{C^\vee C_n})$. The dual affine root system $R^d=R(D^d)$ is the reduced
irreducible affine root subsystem $R^d=\widehat{\mathcal{O}}_2\cup
\widehat{\mathcal{O}}_4\cup\widehat{\mathcal{O}}_5$ (it is of untwisted
type, with underlying finite root system $R_0^\vee$ of type $\textup{C}_n$). 
The Macdonald-Koornwinder theory associated to the initial data $D$ thus is
the special case of the $C^\vee C_n$ theory 
when the multiplicity function
$\kappa\in\mathcal{M}(D^{C^\vee C_n})$ takes
the value one at the two orbits $\widehat{\mathcal{O}}_1=W(a_0)$ and
$\widehat{\mathcal{O}}_4=W(2a_0)$ of $R^{C^\vee C_n}\setminus R$.
\end{rema}
\subsection{Exceptional nonreduced
 rank two Macdonald-Koornwinder polynomials}\label{Except2}
For this special case of the theory we 
take $D=(R_0,\Delta_0,u,\mathbb{Z}R_0,\mathbb{Z}R_0^\vee)$ with
$R_0\subset Z=V=\mathbb{R}^2$ the root system of type $\textup{C}_2$
(it is more natural to view it as type $\textup{C}_2$ then of type
$\textup{B}_2$) given by $R_0=R_{0,s}\cup R_{0,l}$ with
\begin{equation*}
\begin{split}
R_{0,s}&=\{\pm(\epsilon_1+\epsilon_2),\pm(\epsilon_1-\epsilon_2)\},\\
R_{0,l}&=\{\pm 2\epsilon_1,\pm 2\epsilon_2\},
\end{split}
\end{equation*}
where $\{\epsilon_1,\epsilon_2\}$ is the standard orthonormal basis
of $\mathbb{R}^2$. As ordered basis take $\Delta_0=(\alpha_1,\alpha_2)=
(\epsilon_1-\epsilon_2,2\epsilon_2)$. Then $\varphi=2\epsilon_1$ and
$\theta=\epsilon_1+\epsilon_2$. Hence the simple affine root $a_0$
of the associated reduced affine root system $R^u=\mathbb{Z}c+R_0$
is $c-2\epsilon_1$. Then $S=\{1\}$ hence
\[
R=R^u\cup W(2\alpha_1)
\]
with $W=W^u=W_0\ltimes\mathbb{Z}R_0^\vee$. Then $R$ has four $W$-orbits,
\begin{equation*}
\begin{split}
W(a_0)&=(2\mathbb{Z}+1)c+R_{0,l},\\
W(\alpha_1)&=\mathbb{Z}c+R_{0,s},\\
W(\alpha_2)&=2\mathbb{Z}c+R_{0,l},\\
W(2\alpha_1)&=2\mathbb{Z}c+2R_{0,s}.
\end{split}
\end{equation*}
For the dual initial data $D^d=(R_0^\vee,\Delta_0^\vee,u,
\mathbb{Z}R_0^\vee,\mathbb{Z}R_0)$ we have $R^{ud}=\mathbb{Z}c+R_0^\vee$
and $R^d=R^{ud}\cup W^{d}(2\alpha_2^\vee)$ with $W^d=W^{ud}=
W_0\ltimes\mathbb{Z}R_0$. The simple affine root $a_0^d$ of
$R^{ud}$ is $a_0^d=c-\epsilon_1-\epsilon_2$.
The four $W^d$-orbits of $R^d$ are
\begin{equation*}
\begin{split}
W^d(a_0^d)&=(2\mathbb{Z}+1)c+R_{0,s}^\vee,\\
W^d(\alpha_1^\vee)&=2\mathbb{Z}c+R_{0,s}^\vee,\\
W^d(\alpha_2^\vee)&=\mathbb{Z}c+R_{0,l}^\vee,\\
W^d(2\alpha_2^\vee)&=2\mathbb{Z}c+2R_{0,l}^\vee.
\end{split}
\end{equation*}
Note that $R\simeq R^d$ but not $(R,\Delta)\simeq (R^d,\Delta^d)$. 

Let $\kappa\in\mathcal{M}(D)$. We write
\begin{equation}\label{parametersExcept2}
\begin{split}
\{a,b,c,d\}&:=\{\kappa_\theta\kappa_{2\theta},
-\kappa_\theta\kappa_{2\theta}^{-1},\kappa_\varphi^2,q\kappa_0^2\},\\
\{\widetilde{a},\widetilde{b},\widetilde{c},\widetilde{d}\}&=
\{\kappa_\varphi\kappa_0,-\kappa_\varphi\kappa_0^{-1},\kappa_\theta^2,
q\kappa_{2\theta}^2\}
\end{split}
\end{equation}
(the dual parameters $(\widetilde{a},\widetilde{b},\widetilde{c},\widetilde{d})$
are the parameters $(a,b,c,d)$ with respect to the dual initial data $D^d$).
We identify $T_{\mathbb{Z}R_0}\simeq 
(\mathbb{C}^*)^2$ by $t\mapsto (z_1,z_2):=(t^{\alpha_1},t^{\alpha_2})$.
The $W_0$-action on $T_{\mathbb{Z}R_0}$ 
then reads $s_1(z_1,z_2)=(z_1^{-1},z_1^2z_2)$ and
$s_2(z_1,z_2)=(z_1z_2,z_2^{-1})$, giving rise to a 
$W_0$-action
on $\mathbb{C}[T_{\mathbb{Z}R_0}]\simeq\mathbb{C}[z_1^{\pm 1},z_2^{\pm 1}]$
by algebra automorphisms.
The weight function $v_+(z_1,z_2)$ of the associated symmetric 
Macdonald-Koornwinder
polynomials is given by $v_+(z_1,z_2)=\delta(z_1,z_2)\delta(z_1^{-1},z_2^{-1})$
with
\[
\delta(z_1,z_2):=\frac{\bigl(z_1^2,z_1^2z_2^2;q^2\bigr)_{\infty}}
{\bigl(az_1,az_1z_2,bz_1,bz_1z_2;q\bigr)_{\infty}}
\frac{\bigl(z_1^2z_2,z_2;q\bigr)_{\infty}}
{\bigl(cz_1^2z_2,cz_2,dz_1^2z_2,dz_2;q^2\bigr)_{\infty}}.
\]
The symmetric Macdonald-Koornwinder polynomials associated to the initial data
$D$ form an
orthogonal basis of $\mathbb{C}[z_1^{\pm 1},z_1^{\pm 1}]^{W_0}$ with
respect to the sesquilinear, positive pairing
\[
\langle p,r\rangle_+:=\frac{1}{(2\pi i)^2}\int_{|z_1|=1}
\int_{|z_2|=1}p(z_1,z_2)\overline{r(z_1,z_2)}v_+(z_1,z_2)\frac{dz_1}{z_1}
\frac{dz_2}{z_2}
\]
on $\mathbb{C}[z_1^{\pm 1},z_2^{\pm 1}]^{W_0}$ (note that the
monic symmetric Macdonald-Koornwinder polynomials are parametrized by the
$\lambda=\lambda_1\epsilon_1+\lambda_2\epsilon_2\in\Lambda^+$
with the $\lambda_i\in\mathbb{Z}$ satisfying
$\lambda_1+\lambda_2\in 2\mathbb{Z}$ and $\lambda_1\geq\lambda_2\geq 0$).

The associated explicit Macdonald operator is $D_{-\epsilon_1}$. Indeed,
\[
-(Da_0)^{d\vee}=-(Da_0)^\vee=\epsilon_1=\alpha_1^\vee+\alpha_2^\vee
\]
hence
\[
J_{\mathbb{Z}R_0^\vee}=
\{i\in\{1,2\} \,\, | \,\, \varpi_i^d\in \mathbb{Z}R_0^\vee\}=\{1\}
\]
since $\varpi_1^d=\epsilon_1\in\mathbb{Z}R_0^\vee$ and
$\varpi_2^d=\frac{1}{2}(\epsilon_1+\epsilon_2)\not\in\mathbb{Z}R_0^\vee$.
By a direct computation we obtain the explicit expression
\begin{equation*}
\begin{split}
\kappa_0\kappa_\theta^2\kappa_\varphi(D_{-\epsilon_1}f)(z_1,z_2)=&
A(z_1,z_2)f(qz_1,z_2)+
A(z_1^{-1},z_2)f(q^{-1}z_1,q^2z_2)\\
+&A(z_1^{-1}z_2^{-1},z_1^2z_2)f(qz_1,q^{-2}z_2)+
A(z_1^{-1}z_2^{-1},z_2)f(q^{-1}z_1,z_2),\\
A(z_1,z_2):=&\frac{(1-az_1)(1-az_1z_2)(1-bz_1)(1-bz_1z_2)(1-cz_1^2z_2)
(1-dz_1^2z_2)}{(1-z_1^2)(1-z_1^2z_2^2)(1-z_1^2z_2)(1-qz_1^2z_2)}.
\end{split}
\end{equation*}

\section{Double affine Hecke algebras and normalized
Macdonald-Koornwinder polynomials}\label{4}
Cherednik's \cite{CKZ2,CAnn} double affine braid group and double affine Hecke
algebra are fundamental for proving properties
of the (non)symmetric Macdonald-Koornwinder polynomials such as 
evaluation formula, duality and quadratic norms
(see \cite{Cev,Cns}).
We discuss these results in this section.

The first part follows closely \cite{Ha}. 
We define a double affine braid group and a
double affine Hecke algebra depending on $q$, on 
the initial data $D$ and, in case of the double affine Hecke algebra, 
on a choice of a multiplicity function $\kappa\in\mathcal{M}(D)$. We
extend the duality $D\mapsto D^d$ on initial data to an antiisomorphism
of the associated double affine braid groups, and to an algebra 
antiisomorphism
of the associated double affine Hecke algebras (on the dual side,
the double affine Hecke algebra is taken with respect to
the dual multiplicity function $\kappa^d\in\mathcal{M}(D^d)$ as defined
in Lemma \ref{dualmult}). These antiisomorphisms are called duality 
antiisomorphisms
and find their origins in the work of Cherednik \cite{Cev,Cns}
(see also \cite{Ion,Sa,M,Ha} for further results and generalizations).

As a consequence of the duality 
antiisomorphism and of the theory of intertwiners
we first derive the evaluation
formula and duality of the Macdonald-Koornwinder polynomials following closely
\cite{Cns,CI,Sa}. These results, together with quadratic norm formulas,
were first conjectured by Macdonald (see, e.g., \cite{Mpol,Mpr}).
For $R_0$ of type $A$ the evaluation formula and duality were 
proven by Koornwinder \cite{Kpr} by different methods (see \cite{Mbook}
for a detailed account). Subsequently Cherednik \cite{Cev} established
the evaluation formula and duality when 
$(\Lambda,\Lambda^d)=(P(R_0),P(R_0^d))$. The $C^\vee C$ case was established
using Koornwinder's \cite{Kpr} methods by van Diejen \cite{vD} for a suitable
subset of multiplicity functions. Cherednik's double affine Hecke algebra
methods were extended to the $C^\vee C$ case in the work of Noumi \cite{N}
and Sahi \cite{Sa}, leading to the evaluation formula and duality for all
multiplicity functions. Our uniform approach is close
to Haiman's \cite{Ha}.

Following closely \cite{Cns,CI} we use the duality antiisomorphism
and intertwiners to establish quadratic norm formulas for the (non)symmetric 
Macdonald-Koornwinder polynomials.

We fix throughout this section initial data $D=(R_0,\Delta_0,\bullet,
\Lambda,\Lambda^d)\in\mathcal{D}$, a deformation parameter 
$q\in\mathbb{R}_{>0}\setminus\{1\}$, and a multiplicity function 
$\kappa\in\mathcal{M}(D)$.

\subsection{Double affine braid groups, Weyl groups and Hecke algebras}
Consider the $W$-stable additive
subgroup $\widehat{\Lambda}:=\Lambda+\mathbb{R}c$ of $\widehat{Z}$.
In the following definition it is convenient to write 
$X^{\widehat{\lambda}}$ ($\widehat{\lambda}\in\widehat{\Lambda}$) for the elements
of $\widehat{\Lambda}$. In particular,
$X^{\widehat{\lambda}}X^{\widehat{\mu}}=X^{\widehat{\lambda}+\widehat{\mu}}$ 
for $\widehat{\lambda},\widehat{\mu}\in\widehat{\Lambda}$
and $X^0=1$.

We write $A=A(R,\Delta)$ and 
$A_0=A(R_0,\Delta_0)$. The generators of the affine braid group
$\mathcal{B}^\bullet=
\mathcal{B}(A)$ are denoted by $T_0,T_1,\ldots,T_n$, the
generators of $\mathcal{B}_0=\mathcal{B}(A_0)$ by $T_1,\ldots,T_n$.
Recall that elements $T_w\in\mathcal{B}$ ($w\in 
W^\bullet$)
and $T_w\in\mathcal{B}_0$ ($w\in W_0$) can be defined using reduced
expressions of $w$ in the Coxeter groups $W^\bullet=W(A^\bullet)$ and
$W_0=W(A_0)$, respectively. Recall furthermore that $\mathcal{B}=
\mathcal{B}(D)\simeq\Omega\ltimes\mathcal{B}^\bullet$
denotes the extended affine braid group (cf. Definition \ref{ea}).
\begin{defi}\label{DAB}
The double affine braid group $\mathbb{B}=\mathbb{B}(D)$ 
is the group generated by the groups $\mathcal{B}$ and
$\widehat{\Lambda}$
together with the relations:
\begin{equation}\label{defBB1}
T_iX^{\widehat{\lambda}}=X^{\widehat{\lambda}}T_i
\end{equation}
if $\widehat{\lambda}\in\widehat{\Lambda}$ and $0\leq i\leq n$ such that
$(\widehat{\lambda},a_i^\vee)=0$;
\begin{equation}\label{defBB2}
T_iX^{\widehat{\lambda}}T_i=X^{s_i\widehat{\lambda}}
\end{equation}
if $\widehat{\lambda}\in\widehat{\Lambda}$ and $0\leq i\leq n$ such that
$(\widehat{\lambda},a_i^\vee)=1$;
\begin{equation}\label{defBB3}
\omega X^{\widehat{\lambda}}=X^{\omega\widehat{\lambda}}\omega
\end{equation}
if $\widehat{\lambda}\in\widehat{\Lambda}$ and $\omega\in\Omega$.
\end{defi}
Note that $X^{\mathbb{R}c}$ is contained in the center of $\mathbb{B}$.
It is not necessarily
true that any element $g\in \mathbb{B}$ can be written as
$g=bX^{\widehat{\lambda}}$ (or $g=X^{\widehat{\lambda}}b$) with $b\in\mathcal{B}$
and $\widehat{\lambda}\in\widehat{\Lambda}$.
The cases that this possibly fails to be true is
if $S(D)\not=\emptyset$
(these are the cases for which $R\not=R^\bullet$, i.e. 
for which nonreduced extensions of $R^\bullet$
play a role in the theory, cf. Subsection \ref{extension}). 
Straightening the elements of $\mathbb{B}$ is always possible
as soon as quotients to double affine Weyl groups and double
affine Hecke algebras are taken, as we shall see in a moment.

Recall from Subsection \ref{Bernstein} the construction of
commuting group elements $Y^\xi\in\mathcal{B}$ 
($\xi\in\Lambda^d$).
It leads to a presentation of $\mathcal{B}$
in terms of the group $\Lambda^d$ 
(with its elements denoted by $Y^\xi$ ($\xi\in\Lambda^d$))
and the braid group $\mathcal{B}_0=\mathcal{B}(A_0)$, see 
e.g. \cite[\S 3.3]{M}.
On the level
of double affine braid groups it implies the following alternative presentation
of $\mathbb{B}$ (recall that $Da_0$ equals
$-\varphi$ if $\bullet=u$ and $-\theta$ if $\bullet=t$).
\begin{prop}
The double affine braid group $\mathbb{B}$ is isomorphic to the group generated
by the groups $\mathcal{B}_0$, $\Lambda^d$ and $\widehat{\Lambda}$,
satisfying for $1\leq i\leq n$, $\lambda\in\Lambda$
and $\xi\in\Lambda^d$, 
\begin{enumerate}
\item[{\bf (a)}] $X^{\mathbb{R}c}$ is contained in the center,
\item[{\bf (b)}] 
\begin{enumerate}
\item[{\bf 1.}] $Y^{-\xi}T_i=T_iY^{-\xi}$ if $(\xi,\alpha_i^{d\vee})=0$,
\item[{\bf 2.}] $T_iX^\lambda =X^\lambda T_i$ if $(\lambda,\alpha_i^\vee)=0$,
\end{enumerate}
\item[{\bf (c)}]
\begin{enumerate}
\item[{\bf 1.}] $T_iY^{-\xi}T_i=Y^{-s_i\xi}$ if $(\xi,\alpha_i^{d\vee})=1$,
\item[{\bf 2.}] $T_iX^\lambda T_i=X^{s_i\lambda}$ if $(\lambda,\alpha_i^\vee)=1$,
\end{enumerate}
\item[{\bf (d)}] if $(\lambda,a_0^\vee)=0$,
\[
\bigl(Y^{-(Da_0)^d}T_{s_{Da_0}}^{-1}\bigr)X^\lambda=
X^\lambda\bigl(Y^{-(Da_0)^d}T_{s_{Da_0}}^{-1}\bigr),
\]
\item[{\bf (e)}] if $(\lambda,a_0^\vee)=1$,
\[
\bigl(Y^{-(Da_0)^d}T_{s_{Da_0}}^{-1}\bigr)X^\lambda
\bigl(Y^{-(Da_0)^d}T_{s_{Da_0}}^{-1}\bigr)=p_{Da_0}^{-1}
X^{s_{Da_0}\lambda}
\]
where $p_\alpha:=X^{\mu_\alpha c}$ for $\alpha\in R_0$.
\end{enumerate}
\end{prop}
Recall that $u(\eta)\in W\simeq W_0\ltimes\Lambda^d$ denotes the element
of minimal length in $\tau(\eta)W_0$ for all $\eta\in\Lambda^d$. 
Then $u(\eta)=\tau(\eta)v(\eta)^{-1}$
with $v(\eta)\in W_0$ the element of minimal length such that 
$v(\eta)\eta=\eta_-$. Recall furthermore that $\Omega=\{u(\xi)\, | \,
\xi\in\Lambda_{min}^{d+}\}$.

The identification of the two different sets of generators
of $\mathbb{B}$ then is as follows: 
$T_0=Y^{-(Da_0)^d}T_{s_{Da_0}}^{-1}$ and 
$u(\xi)=Y^\xi T_{v(\xi)}^{-1}$ for $\xi\in\Lambda_{min}^{d+}$. Conversely,
for $\xi=\eta_1-\eta_2\in\Lambda^d$ with $\eta_1,\eta_2\in\Lambda^{d+}$,
$Y^\xi=Y^{\eta_1}\bigl(Y^{\eta_2}\bigr)^{-1}$ and
\[
Y^{\eta_s}=T_{\tau(\eta_s)}=\omega T_{i_1}T_{i_2}\cdots T_{i_r} 
\]
if $\tau(\eta_s)=\omega s_{i_1}s_{i_2}\cdots s_{i_r}\in W$
is a reduced expression ($\omega\in\Omega$ and $0\leq i_j\leq n$).

Recall the set $S=S(D)$ given by \eqref{S}.
Write 
\begin{equation}\label{Vi}
V_i:=X^{-a_i}T_i^{-1}\in\mathbb{B}\qquad \forall\, i\in S.
\end{equation}
Recall furthermore that $R=R^\bullet\cup\bigcup_{i\in S}W^\bullet(2a_i)$.
\begin{defi}
{\bf (i)} The double affine Weyl group $\mathbb{W}=\mathbb{W}(D)$ 
is the quotient
of $\mathbb{B}$ by the normal subgroup generated by
$T_j^2$ ($j\in\{0,\ldots,n\}$) and $V_i^2$ ($i\in S$).\\
{\bf (ii)} The double affine Hecke algebra $\mathbb{H}(\kappa,q)=
\mathbb{H}(D;\kappa,q)$ is $\mathbb{C}[\mathbb{B}]/
\widehat{I}_{\kappa,q}$, where $\widehat{I}_{\kappa,q}$ is the two-sided
ideal of $\mathbb{C}[\mathbb{B}]$ generated by $X^{rc}-q^{r}$
($r\in\mathbb{R}$),
$(T_j-\kappa_j)(T_j+\kappa_j^{-1})$ ($j\in\{0,\ldots,n\}$) and
$(V_i-\kappa_{2a_i})(V_i+\kappa_{2a_i}^{-1})$ ($i\in S$).
\end{defi}
Recall that $W$ acts on $\widehat{\Lambda}$ by group automorphisms.
\begin{prop}
$\mathbb{W}\simeq W\ltimes\widehat{\Lambda}$.
\end{prop}

For the double affine Hecke algebra, note that we have canonical
algebra homomorphisms 
$H(\kappa^\bullet)\rightarrow \mathbb{H}(\kappa,q)$
and $\mathbb{C}[\widehat{\Lambda}]\rightarrow \mathbb{H}(\kappa,q)$. We write
$\widetilde{h}$, $\widetilde{X}_q^{\widehat{\lambda}}$ and $\widetilde{X}^\lambda$ 
for the images of
$h\in\widehat{H}$, $X^{\widehat{\lambda}}\in\mathbb{C}[\widehat{\Lambda}]$ 
($\widehat{\lambda}\in\widehat{\Lambda}$)
and $X^\lambda\in\mathbb{C}[\Lambda]$ ($\lambda\in\Lambda$)
in $\mathbb{H}(\kappa,q)$ respectively.
Define a linear map
\[
m: H(\kappa^\bullet)\otimes_{\mathbb{C}}\mathbb{C}[\Lambda]
\rightarrow \mathbb{H}(\kappa,q)
\]
by $m(h\otimes X^\lambda):=\widetilde{h}\widetilde{X}^\lambda$
for $h\in H(\kappa^\bullet)$ and $\lambda\in\Lambda$.

If 
$\widehat{\lambda}=\lambda+r c\in\widehat{\Lambda}$  
then we interpret $e_q^{\widehat{\lambda}}$
as the endomorphism of $\mathbb{C}[T_\Lambda]$ by
\[\bigl(e_q^{\widehat{\lambda}}p\bigr)(t):=t_q^{\widehat{\lambda}}p(t)=
q^{r}t^{\lambda}p(t),\qquad r\in\mathbb{R},\,\, \lambda\in\Lambda.
\] 
\begin{thm}
{\bf (i)} We have a unique algebra monomorphism 
$\mathbb{H}(\kappa,q)\hookrightarrow\textup{End}_{\mathbb{C}}\bigl(
\mathbb{C}[T_\Lambda]\bigr)$ defined by
\begin{equation*}
\begin{split}
\widetilde{h}&\mapsto \pi_{\kappa,q}(h)\qquad \forall\, 
h\in H(\kappa^\bullet),\\
\widetilde{X}_q^{\widehat{\lambda}}&\mapsto e_q^{\widehat{\lambda}}\qquad
\qquad\forall\, \widehat{\lambda}\in\widehat{\Lambda}.
\end{split}
\end{equation*}
{\bf (ii)} The linear map $m$ defines a complex linear isomorphism
$H(\kappa^\bullet)\otimes_{\mathbb{C}}\mathbb{C}[\Lambda]
\overset{\sim}{\longrightarrow}\mathbb{H}(\kappa,q)$.
\end{thm}

To simplify notations, we omit the tilde when writing 
the elements in $\mathbb{H}$. With this convention $X^{\widehat{\lambda}}_q=
q^{r}X^\lambda$ in $\mathbb{H}(\kappa,q)$
if $\widehat{\lambda}=\lambda+rc$.

Together with the Bernstein-Zelevinsky presentation of the extended affine Hecke
algebra $H(\kappa^\bullet)$ (see Subsection \ref{Bernstein})
we conclude that
\[
\mathbb{C}[T_{\Lambda^d}]\otimes_{\mathbb{C}}H_0(\kappa|_{R_0})\otimes_{\mathbb{C}}
\mathbb{C}[T_\Lambda]\simeq \mathbb{H}(\kappa,q)
\]
as complex vector spaces by mapping $e^\xi\otimes h\otimes e^\lambda$
to $Y^\xi h X^\lambda$ ($\xi\in\Lambda^d$, $h\in H_0(\kappa|_{R_0})$ and
$\lambda\in\Lambda$). This is the Poincar{\'e}-Birkhoff-Witt property
of the double affine Hecke algebra $\mathbb{H}(\kappa,q)$.

The dual version of the Lusztig relations \eqref{Lusztig1}
and \eqref{Lusztig2}, now also including a commutation relation
for $T_0$, is given as follows. Write $p(X)\in\mathbb{H}(\kappa,q)$
for the element corresponding to $p\in\mathbb{C}[T_\Lambda]$. In other
words, $p(X)=\sum_{\lambda\in\Lambda}c_\lambda X^\lambda$ 
if $p(t)=\sum_{\lambda\in\Lambda}c_\lambda t^\lambda$.
\begin{cor}
Let $0\leq i\leq n$ and $p\in\mathbb{C}[T_\Lambda]\subset\mathbb{H}(\kappa,q)$.
Then 
\begin{equation}\label{dualLusztig}
T_ip(X)-(s_{i,q}p)(X)T_i=
\left(\frac{\kappa_i-\kappa_i^{-1}+
(\kappa_{2a_i}-\kappa_{2a_i}^{-1})X_q^{a_i}}{1-X_q^{2a_i}}\right)
(p(X)-(s_{i,q}p)(X))
\end{equation}
in $\mathbb{H}(\kappa,q)$.
\end{cor}
\subsection{Duality antiisomorphism}\label{dualitysection}
Recall that we have associated to 
$D=(R_0,\Delta_0,\bullet,\Lambda,\Lambda^d)$ 
the dual initial data $D^d=(R_0^d,\Delta_0^d,\bullet,\Lambda^d,\Lambda)$,
see \eqref{dualquintuples}. 

We define the dual double affine braid group by $\mathbb{B}^d:=
\mathbb{B}(D^d)$.
We add a superscript $d$ to elements of $\mathbb{B}^d$ if 
confusion may arise. So we write ${}^dY^\lambda$ ($\lambda\in\Lambda$),
${}^dX^\xi$ ($\xi\in\Lambda^d$) and $T_i^d$ ($0\leq i\leq n$) in
$\mathbb{B}^d$. The group generators $T_i^d$ ($1\leq i\leq n$)
of the homomorphic image of $\mathcal{B}_0$ in $\mathbb{B}^d$
will usually be written without superscripts. 

Recall that 
$Da_0=-\varphi$
(respectively $-\theta$) if $\bullet=u$ (respectively $\bullet=t$).
We have $s_0=\tau(-(Da_0)^d)s_{Da_0}\in W^\bullet\subseteq W$ and 
$T_0=Y^{-(Da_0)^d}T_{s_{Da_0}}^{-1}\in\mathbb{B}$.
Dually, $D(a_0^d)=-\theta^d$, $s_0^d=\tau(\theta)s_\theta\in W^d$ 
and $T_0^d=Y^{\theta} T_{s_\theta}^{-1}\in \mathbb{B}^d$.
The following result, in the present generality, is from \cite[\S 4]{Ha} 
(in \cite{Ha} the duality isomorphism is constructed, which is related
to the duality antiisomorphism below via an elementary antiisomorphism).
See also \cite{Cev,Cns,Ion,Sa,M} for special cases.
\begin{thm}\label{dualB}
There exists a unique antiisomorphism $\delta: 
\mathbb{B}\rightarrow\mathbb{B}^d$
satisfying
\begin{equation*}
\begin{split}
\delta(X^{rc})&=X^{rc}\qquad\qquad\, \forall\, r\in\mathbb{R},\\
\delta\bigl(Y^\xi\bigr)&={}^dX^{-\xi}\qquad\quad\,\,\,
\forall\, \xi\in\Lambda^d,\\
\delta\bigl(T_i\bigr)&=T_i\qquad\qquad\,\,
\,\,\,\, \forall\, i\in\{1,\ldots,n\},\\
\delta\bigl(X^\lambda\bigr)&={}^dY^{-\lambda}\qquad\quad\,\,\,\,
\forall\, \lambda\in\Lambda.
\end{split}
\end{equation*}
\end{thm}
Note that $\delta^d=\delta^{-1}$, where $\delta^d$ is the duality
antiisomorphism with respect to the dual initial data $D^d$.

Recall that for $\lambda\in\Lambda$, $v(\lambda)$ is the element in $W_0$
of smallest length such that $\lambda_-=v(\lambda)\lambda$.
Then $\tau(\lambda)=u^d(\lambda)v(\lambda)$ in $W^d$,
where $u^d(\lambda)$ is the shortest element (with respect to the length 
function $l^d$ on $W^d$) of the coset $\tau(\lambda)W_0$ in $W^d$.
Then $\Omega=\{u^d(\lambda)\, | \, \lambda\in\Lambda_{min}^+\}$ and in
$\mathbb{B}^d$,
\[
{}^dY^{\lambda}=u^d(\lambda)T_{v(\lambda)}
\]
hence
\[
\delta^d(u^d(\lambda))=T_{v(\lambda)^{-1}}^{-1}X^{-\lambda}
\]
in $\mathbb{B}$. On the other hand, $\tau(\theta)=s_0^ds_\theta$ and
$l^d(\tau(\theta))=l(s_\theta)+1$, hence
\[
T_0^d={}^dY^{\theta}T_{s_\theta}^{-1}
\]
in $\mathbb{B}^d$ and
\[
\delta^{d}(T_0^d)=T_{s_\theta}^{-1}X^{-\theta}
\]
in $\mathbb{B}$.

The next result shows that the duality antiisomorphism from Theorem
\ref{dualB} descends to an antiisomorphism between double affine Hecke
algebras. In order to do so, we use the
isomorphism $\mathcal{M}\overset{\sim}{\longrightarrow}\mathcal{M}^d$
from the complex torus $\mathcal{M}=\mathcal{M}(D)$ onto 
$\mathcal{M}^d=\mathcal{M}(D^d)$ from Lemma \ref{dualmult}.
\begin{thm}
The antiisomorphism $\delta: 
\mathbb{B}\rightarrow
\mathbb{B}^d$ descends
to an antiisomorphism $\delta: 
\mathbb{H}(\kappa,q)\rightarrow
\mathbb{H}^d(\kappa^d,q):=\mathbb{H}(D^d;\kappa^d,q)$.
\end{thm}
For instance, the cross relations \eqref{dualLusztig} 
in $\mathbb{H}(\kappa,q)$ for $i\in\{1,\ldots,n\}$ match with the 
the cross relations \eqref{Lusztiguniform} in $\mathbb{H}^d(\kappa^d,q)$
through the antiisomorphism $\delta$.
Note that $\delta^d=\delta^{-1}$ also on the level of the double affine
Hecke algebra, where $\delta^d$ is the duality
antiisomorphism with respect to the dual data $(D^d,\kappa^d)$.

\subsection{Evaluation formulas}
We follow closely the arguments from \cite{Cns} (which relates to the special
case that $(\Lambda,\Lambda^d)=(P(R_0),P(R_0^d))$). See also \cite{St,M,Ha}.
We assume in this subsection that $q$ and $\kappa\in\mathcal{M}(D)$ satisfy
either \eqref{smallparameter} or \eqref{largeparameter}.
Then we have the monic nonsymmetric Macdonald-Koornwinder polynomials
$P_\lambda\in\mathbb{C}[T_\Lambda]$ ($\lambda\in\Lambda$) associated
to $(D,\kappa,q)$, as well as the dual monic nonsymmetric 
Macdonald-Koornwinder polynomials
$P_\xi^d\in\mathbb{C}[T_{\Lambda^d}]$ ($\xi\in\Lambda^d$) associated to 
$(D^d,\kappa^d,q)$. 

Write $\gamma_{\xi,q}^d=\gamma_{\xi,q}(D^d;\kappa^{d\bullet})
\in T_\Lambda$ ($\xi\in\Lambda^d$) for the spectral points with respect to
the dual initial data $(D^d,\kappa^d,q)$. Concretely, they are given by
\[
\gamma_{\xi,q}^d=q^\xi\prod_{\alpha\in R_0^+}
{}^d\upsilon_\alpha^{\eta((\xi,\alpha^{d\vee}))\alpha^\vee}
\]
with ${}^d\upsilon_\alpha:=(\kappa_{\alpha^d}^d)^{\frac{1}{2}}
(\kappa_{\mu_{\alpha^d}c+\alpha^d}^d)^{\frac{1}{2}}=
\kappa_\alpha^{\frac{1}{2}}\kappa_{2\alpha}^{\frac{1}{2}}$ for $\alpha\in R_0^+$.

Cherednik's basic representations 
\begin{equation*}
\begin{split}
\pi_{D;\kappa,q}: H(\kappa^\bullet)\hookrightarrow
\textup{End}_{\mathbb{C}}\bigl(\mathbb{C}[T_\Lambda]\bigr),\\
\pi_{D^d;\kappa^d,q}: H(\kappa^{d\bullet})\hookrightarrow
\textup{End}_{\mathbb{C}}\bigl(\mathbb{C}[T_{\Lambda^d}]\bigr)
\end{split}
\end{equation*}
extend to algebra maps
\begin{equation*}
\begin{split}
\widehat{\pi}: \mathbb{H}(\kappa,q)\hookrightarrow
\textup{End}_{\mathbb{C}}\bigl(\mathbb{C}[T_\Lambda]\bigr),\\
\widehat{\pi}^d: \mathbb{H}(\kappa^{d},q)\hookrightarrow
\textup{End}_{\mathbb{C}}\bigl(\mathbb{C}[T_{\Lambda^d}]\bigr)
\end{split}
\end{equation*}
by $\widehat{\pi}\bigl(X_q^{\widehat{\lambda}}\bigr)=
e_q^{\widehat{\lambda}}$ for $\widehat{\lambda}\in\widehat{\Lambda}$
and $\widehat{\pi}^d\bigl({}^dX_q^{\widehat{\xi}}\bigr)=
e_q^{\widehat{\xi}}$ for $\widehat{\xi}\in\widehat{\Lambda^d}$.

Note that
\[
\gamma_{0,q}^d=
\prod_{\alpha\in R_0^+}{}^d\upsilon_\alpha^{-\alpha^\vee}\in T_\Lambda,
\qquad
\gamma_{0,q}=
\prod_{\alpha\in R_0^+}\upsilon_\alpha^{-\alpha^{d\vee}}\in T_{\Lambda^d}.
\]
\begin{defi} Define evaluation maps 
$\textup{Ev}: \mathbb{H}(\kappa,q)\rightarrow\mathbb{C}$ 
and $\textup{Ev}^d: {}^d\mathbb{H}(\kappa^d,q)\rightarrow\mathbb{C}$
by
\[
\textup{Ev}(Z):=\bigl(\widehat{\pi}(Z)1\bigr)(\gamma_{0,q}^{d}),
\qquad Z\in\mathbb{H}(\kappa,q),
\]
\[
\textup{Ev}^d(Z):=\bigl(\widehat{\pi}^d(Z)1\bigr)(\gamma_{0,q}),\qquad
Z\in\mathbb{H}^d(\kappa^d,q).
\]
\end{defi}
The following lemma is crucial for the duality of the (nonsymmetric)
Macdonald-Koornwinder polynomials.
\begin{lem}
For all $Z\in\mathbb{H}(\kappa,q)$,
\[
\textup{Ev}^d(\delta(Z))=\textup{Ev}(Z).
\]
\end{lem}
Write $\mathbb{H}=\mathbb{H}(\kappa,q)$ and 
$\mathbb{H}^d=\mathbb{H}^d(\kappa^d,q)$. Define bilinear forms
\begin{equation*}
\begin{split}
B: \mathbb{H}\times\mathbb{H}^d\rightarrow\mathbb{C},\qquad
&B(Z,\widetilde{Z})=\textup{Ev}(\delta^d(\widetilde{Z})Z),\\
B^d: \mathbb{H}^d\times\mathbb{H}\rightarrow\mathbb{C},\qquad
&B^d(\widetilde{Z},Z)=\textup{Ev}^d(\delta(Z)\widetilde{Z})
\end{split}
\end{equation*}
for $Z\in\mathbb{H}$ and $\widetilde{Z}\in\mathbb{H}^d$.
\begin{cor}\label{aaaa}
$B(Z,\widetilde{Z})=B^d(\widetilde{Z},Z)$ for all $Z\in\mathbb{H}$
and all $\widetilde{Z}\in\mathbb{H}^d$.
\end{cor}
The following elementary lemma provides convenient tools to derive the
evaluation formula for nonsymmetric Macdonald-Koornwinder polynomials.
\begin{lem}\label{bbbb}
Let $p\in\mathbb{C}[T_\Lambda]$, $\widetilde{p}\in\mathbb{C}[T_{\Lambda^d}]$,
$Z,Z_1,Z_2\in\mathbb{H}$ and $\widetilde{Z},\widetilde{Z}_1,\widetilde{Z}_2
\in\mathbb{H}^d$. Then\\
{\bf (i)} $B(Z_1Z_2,\widetilde{Z})=B(Z_2,\delta(Z_1)\widetilde{Z})$.\\
{\bf (ii)} $B(ZT_i,\widetilde{Z})=\kappa_iB(Z,\widetilde{Z})$ for
$0\leq i\leq n$.\\
{\bf (iii)} $B((\widehat{\pi}(Z)(p))(X),\widetilde{Z})=
B(Zp(X),\widetilde{Z})$.
\end{lem}
Lemma \ref{bbbb} and Corollary \ref{aaaa} imply
\begin{prop}
{\bf (i)} For $\xi\in\Lambda^d$ and $p\in\mathbb{C}[T_\Lambda]$,
\[
P_\xi^d(\gamma_{0,q})p(\gamma_{\xi,q}^d)=B(p,P_\xi^d).
\]
{\bf (ii)} For $\lambda\in\Lambda$ and $\widetilde{p}\in
\mathbb{C}[T_{\Lambda^d}]$,
\[
P_\lambda(\gamma_{0,q}^d)\widetilde{p}(\gamma_{\lambda,q})=
B^d(\widetilde{p},P_\lambda).
\]
\end{prop}

\begin{cor}[duality]
For $\lambda\in\Lambda$ and $\xi\in\Lambda^d$ we have
\[
P_\xi^d(\gamma_{0,q})P_\lambda(\gamma_{\xi,q}^d)=
P_\lambda(\gamma_{0,q}^d)P_\xi^d(\gamma_{\lambda,q}).
\]
\end{cor}
The next aim is to explicitly evaluate 
$\textup{Ev}(P_\lambda)=P_\lambda(\gamma_{0,q}^d)$.

Using the results on the pairing $B$ above and using
Proposition \ref{Wgammaaction}, one first derives an important intermediate
result which shows how the action of generators of 
the double affine Hecke algebra on monic nonsymmetric Macdonald-Koornwinder 
polynomials $P_\lambda$ can be explicitly expressed in terms of an 
action on the degree $\lambda\in\Lambda$ of $P_\lambda$.

Recall that $W^d$ acts on $\Lambda$ by $(w\tau(\lambda),\lambda^\prime)\mapsto
w(\lambda+\lambda^\prime)$ ($w\in W_0$ and $\lambda,\lambda^\prime\in\Lambda$).
\begin{prop}\label{DAHAaction}
{\bf (i)} If $s_i^d\lambda=\lambda$ for $1\leq i\leq n$ then
\[
\widehat{\pi}(T_i)P_\lambda=\kappa_i^dP_\lambda.
\]
{\bf (ii)} If $s_0^d\lambda=\lambda$ then
\[
\widehat{\pi}(\delta^{d}(T_0^d))P_\lambda=\kappa_0^dP_\lambda.
\]
{\bf (iii)} Let $\lambda\in\Lambda_{min}^+$, then
\[
\widehat{\pi}\bigl(\delta^{d}(u^d(\lambda))\bigr)1=
\kappa_{v(\lambda)}^dP_{-\lambda_-}.
\]
{\bf (iv)} Suppose $1\leq i\leq n$, $\lambda\in\Lambda$ such that
$(\lambda,\alpha_i^\vee)>0$. Then
\[
\widehat{\pi}(T_i)P_\lambda=
\Bigl(\frac{\kappa_i^d-(\kappa_i^d)^{-1}+
(\kappa_{2\alpha_i^d}^d-(\kappa_{2\alpha_i^d}^d)^{-1})\gamma_{\lambda,q}^{\alpha_i^d}}
{1-\gamma_{\lambda,q}^{2\alpha_i^d}}\Bigr)P_\lambda
+(\kappa_i^d)^{-1}P_{s_i^d\lambda}.
\]
{\bf (v)} Suppose $\lambda\in\Lambda$ such that $a_0^d(\lambda)>0$. Then
\[
\widehat{\pi}\bigl(\delta^{-1}(T_0^d)\bigr)P_\lambda=
\Bigl(\frac{\kappa_0^d-(\kappa_0^d)^{-1}+
(\kappa_{2a_0^d}^d-(\kappa_{2a_0^d}^d)^{-1})q_\theta\gamma_{\lambda,q}^{-\theta^d}}
{1-q_\theta^2\gamma_{\lambda,q}^{-2\theta^d}}\Bigr)P_\lambda
+\kappa_{v(s_0^d\lambda)}^d(\kappa_{v(\lambda)}^d)^{-1}P_{s_0^d\lambda}
\]
(note that $a_0^d=\mu_{\theta}c-\theta^d=\mu_\theta(c-\theta^\vee)$
and $s_0^d\lambda=\lambda+(1-(\lambda,\theta^\vee))\theta$).
\end{prop}

The proposition gives the following recursion relations for
$\textup{Ev}(P_\lambda)=P_\lambda(\gamma_{0,q}^d)$.
\begin{cor}
{\bf (i)} If $\lambda\in\Lambda_{min}^+$ then
\[\textup{Ev}(P_\lambda)
=(\kappa_{v(\lambda)}^d)^{-1}.
\]
{\bf (ii)} If $\lambda\in\Lambda$, $1\leq i\leq n$ and $a_i^d(\lambda)>0$ then
\[
\textup{Ev}(P_{s_i^d\lambda})=
\frac{(1-\kappa_i^d\kappa_{2\alpha_i^d}^d\gamma_{\lambda,q}^{\alpha_i^d})
(1+\kappa_i^d(\kappa_{2\alpha_i^d}^d)^{-1}\gamma_{\lambda,q}^{\alpha_i^d})}
{(1-\gamma_{\lambda,q}^{2\alpha_i^d})}
\textup{Ev}(P_\lambda).
\]
{\bf (iii)} If $\lambda\in\Lambda$ and $a_0^d(\lambda)>0$ then
\[
\textup{Ev}(P_{s_0^d\lambda})=
(\kappa_0^d)^{-1}\kappa_{v(\lambda)}^d(\kappa_{v(s_0^d\lambda)}^d)^{-1}
\frac{(1-q_\theta\kappa_0^d\kappa_{2a_0^d}^d\gamma_{\lambda,q}^{-\theta^d})
(1+q_\theta\kappa_0^d(\kappa_{2a_0^d}^d)^{-1}\gamma_{\lambda,q}^{-\theta^d})}
{(1-q_\theta^2\gamma_{\lambda,q}^{-2\theta^d})}
\textup{Ev}(P_\lambda).
\]
\end{cor}
Recall the definition $c_a=c_a^{\kappa,q}(\cdot;D)\in\mathbb{C}(T_\Lambda)$
for $a\in R^\bullet$ from \eqref{ca}. The dual version is denoted
by $c_a^d=c_a^{\kappa^d,q}(\cdot;D^d)$ ($a\in R^{d\bullet}$). Concretely, for
$a\in R^{d\bullet}$, $c_a^d\in\mathbb{C}(T_{\Lambda^d})$ is given by
\[
c_a^d(t):=\frac{(1-\kappa_a^d\kappa_{2a}^dt_q^a)
(1+\kappa_a^d(\kappa_{2a}^d)^{-1}t_q^a)}{1-t_q^{2a}}.
\]
We also set for $w\in W^d$,
\[
c_w^d:=\prod_{a\in R^{d\bullet,+}\cap w^{-1}(R^{d\bullet,-})}c_a^d\in
\mathbb{C}(T_{\Lambda^d}).
\]
An induction argument gives now the explicit evaluation formula
for the nonsymmetric Macdonald-Koornwinder polynomials (see \cite{Cns,St,M}).
\begin{thm}\label{evthm}
For $\lambda\in\Lambda$ we have
\[
\textup{Ev}\bigl(P_\lambda\bigr)
=(\kappa_{\tau(\lambda)}^d)^{-1}c_{u^d(\lambda)}^d(\gamma_{0,q}).
\]
\end{thm}

\subsection{Normalized nonsymmetric Macdonald-Koornwinder polynomials
and duality}
The treatment in this subsection is close to \cite{Cns} and \cite{St}, 
which deal with the case that $(\Lambda,\Lambda^d)=(P(R_0),P(R_0^d))$
and $C^\vee C$ case, respectively.
We assume that $q$ and $\kappa\in\mathcal{M}(D)$ satisfy
\eqref{smallparameter} or \eqref{largeparameter}.
Then
\[
P_\lambda(\gamma_{0,q}^d)\not=0\quad \&\quad P_\xi^d(\gamma_{0,q})\not=0
\]
for all $\lambda\in\Lambda$ and $\xi\in\Lambda^d$ in view of the evaluation
formula (Theorem \ref{evthm}).

Recall the Macdonald-Koornwinder polynomials $P_\lambda\in\mathbb{C}[T_\Lambda]$
($\lambda\in\Lambda$) and $P_\xi^d\in\mathbb{C}[T_{\Lambda^d}]$ 
($\xi\in\Lambda^d$)
satisfy
\begin{equation*}
\begin{split}
\widehat{\pi}(p(Y))P_\lambda&=p(\gamma_{\lambda,q}^{-1})P_\lambda,\\
\widehat{\pi}^d(r({}^dY))P_\xi^d&=r(\gamma_{\xi,q}^{d\,-1})P_\xi^d
\end{split}
\end{equation*}
for all $p\in\mathbb{C}[T_{\Lambda^d}]$ and $r\in\mathbb{C}[T_\Lambda]$.
This motivates the following notation for the normalized nonsymmetric
Macdonald-Koornwinder polynomials.
\begin{defi}
The normalized nonsymmetric Macdonald-Koornwinder polynomials are defined by
\[
E(\gamma_{\lambda,q}^{-1};\cdot):=\frac{P_\lambda}{P_\lambda(\gamma_{0,q}^d)}
\in\mathbb{C}[T_\Lambda],\qquad \lambda\in\Lambda
\]
and 
\[
E^d((\gamma_{\xi,q}^d)^{-1};\cdot):=\frac{P_\xi^d}{P_\xi^d(\gamma_{0,q})}
\in\mathbb{C}[T_{\Lambda^d}],\qquad \xi\in\Lambda^d.
\]
\end{defi}
We denote by $E^\circ(\gamma_{\lambda,q};\cdot)\in\mathbb{C}[T_\Lambda]$
the normalized Macdonald-Koornwinder 
polynomial with respect to the inverted parameters
$(\kappa^{-1},q^{-1})$ (and similarly for $E^{d\circ}(\gamma_{\xi,q}^d;\cdot)$).
 
For $\lambda\in\Lambda$ we have $E(\gamma_{\lambda,q}^{-1};\gamma_{0,q}^d)=1$.
{}From the previous subsection we immediately get
the self-duality of the normalized nonsymmetric Macdonald-Koornwinder 
polynomials (cf. \cite{Cns,Sa,M,Ha}).
\begin{cor}
For all $p\in\mathbb{C}[T_\Lambda]$,
$r\in\mathbb{C}[T_{\Lambda^d}]$, $\lambda\in\Lambda$ and $\xi\in\Lambda^d$,
\begin{equation*}
\begin{split}
B(p,E^d(\gamma_{\xi,q}^{d\,-1};\cdot))&=p(\gamma_{\xi,q}^d),\\
B(E(\gamma_{\lambda,q}^{-1};\cdot),r)&=r(\gamma_{\lambda,q}).
\end{split}
\end{equation*}
In particular, $E(\gamma_{\lambda,q}^{-1};\gamma_{\xi,q}^d)=
E^d(\gamma_{\xi,q}^{d\,-1};\gamma_{\lambda,q})$ for all
$\lambda\in\Lambda$ and $\gamma\in\Lambda^d$.
\end{cor}
\subsection{Polynomial Fourier transform}
For the remainder of the text we assume that the
parameters $q$ and $\kappa\in\mathcal{M}(D)$ satisfy 
the more restrictive parameter conditions
\eqref{smallparameterorth}.

In order to compute the norms of the normalized nonsymmetric
Macdonald-Koornwinder 
polynomials it is convenient to formulate the explicit formulas
from Proposition \ref{DAHAaction} in terms of properties of 
a Fourier transform whose kernel is given by the normalized nonsymmetric
Macdonald-Koornwinder polynomial.
For this we first need to consider
the adjoint of the double affine Hecke algebra action
with respect to the sesquilinear form 
\[
\langle p_1,p_2\rangle:=\int_{T_\Lambda^u}p_1(t)\overline{p_2(t)}v(t)dt,
\qquad p_1,p_2\in\mathbb{C}[T_\Lambda]
\]
and express it in terms
of an explicit antilinear antiisomorphism of the double affine
Hecke algebra.
\begin{lem}\label{antilem}
There exists a unique antilinear antialgebra isomorphism
${}^\ddagger: \mathbb{H}(\kappa,q)\rightarrow\mathbb{H}(\kappa^{-1},q^{-1})$
satisfying $T_w^\ddagger=T_w^{-1}$ ($w\in W$)
and $(X^\lambda)^\ddagger=X^{-\lambda}$ ($\lambda\in\Lambda$). In addition,
\[
\langle \widehat{\pi}_{\kappa,q}(h)p_1,p_2\rangle=\langle p_1,
\widehat{\pi}_{\kappa^{-1},q^{-1}}(h^\ddagger)p_2\rangle
\]
for all $h\in\mathbb{H}$.
\end{lem}
Define 
\begin{equation*}
\begin{split}
\mathcal{S}&:=\{\gamma_{\lambda,q} \,\, | \,\, \lambda\in\Lambda\}
\subset T_{\Lambda^d},\\
\mathcal{S}^d&:=\{\gamma_{\xi,q}^d \,\, | \,\, \xi\in\Lambda^d\}\subset
T_\Lambda,
\end{split}
\end{equation*}
and write $F(\mathcal{S})$ (respectively
$F(\mathcal{S}^d)$) for the space of finitely supported complex valued
functions on $\mathcal{S}$ (respectively $\mathcal{S}^d$). The following
lemma follows easily from the results in Subsection \ref{realizations}
and from Proposition \ref{Wgammaaction}. 
\begin{lem} 
There exists a unique algebra homomorphism 
$\widehat{\rho}^d: \mathbb{H}^d\rightarrow 
\textup{End}_{\mathbb{C}}(F(\mathcal{S}))$ satisfying
for $0\leq i\leq n$,
\begin{equation*}
\begin{split}
(\widehat{\rho}^d(T_i^d)g)(\gamma_{\lambda,q})&=\kappa_i^dg(\gamma_{\lambda,q})+
(\kappa_i^d)^{-1}c_{a_i^d}^d(\gamma_{\lambda,q})
(g(\gamma_{s_i^d\lambda,q})-g(\gamma_{\lambda,q}))\quad 
\hbox{ if } s_i^d\lambda\not=\lambda,\\
(\widehat{\rho}^d(T_i^d)g)(\gamma_{\lambda,q})&=\kappa_i^dg(\gamma_{\lambda,q})
\qquad\qquad\qquad\qquad\qquad\qquad\qquad\qquad\,\,\,
\hbox{ if } s_i^d\lambda=\lambda,
\end{split}
\end{equation*}
and satisfying for $\omega\in\Omega^d$,
\[
(\widehat{\rho}^d(\omega)g)(\gamma_{\lambda,q})=g(\gamma_{\omega^{-1}\lambda,q})
\]
and for $\xi\in\Lambda^d$,
\[
(\widehat{\rho}({}^dX^\xi)g)(\gamma_{\lambda,q})=
\gamma_{\lambda,q}^\xi g(\gamma_{\lambda,q}).
\]
\end{lem}

\begin{defi}
We define the complex linear map
\[
\mathcal{F}: \mathbb{C}[T_\Lambda]\rightarrow F(\mathcal{S})
\]
by
\[(\mathcal{F}p)(\gamma):=\langle p,E^\circ(\gamma;\cdot)\rangle,\qquad
p\in\mathbb{C}[T_\Lambda],\,\, \gamma\in\mathcal{S}.
\]
\end{defi}

For generators $X\in\mathbb{H}$ and $p\in\mathbb{C}[T_\Lambda]$
we can re-express
$\mathcal{F}(\widehat{\pi}(X)p)$ as an explicit linear operator
acting on $\mathcal{F}p\in F(\mathcal{S})$ using Lemma \ref{antilem},
Proposition \ref{DAHAaction} and Theorem \ref{evthm}. 
It gives the following result
(the first part of the theorem follows easily from the duality antiisomorphism).
\begin{thm}\label{kpk}
{\bf (i)} The following formulas define an algebra isomorphism
$\Phi: \mathbb{H}\rightarrow\mathbb{H}^d$,
\begin{equation*}
\begin{split}
\Phi(T_i)&=T_i^d,\qquad\qquad\quad 1\leq i\leq n,\\
\Phi(T_{s_\theta}^{-1}X^{-\theta})&=T_0^d,\\
\Phi(T_{v(\lambda)^{-1}}^{-1}X^{-\lambda})&=u^d(\lambda)^{-1},\qquad\,\,\,\,\, 
\lambda\in\Lambda_{min}^+,\\
\Phi(Y^\xi)&={}^dX^{-\xi},\qquad\qquad\,\,\, \xi\in\Lambda^d.
\end{split}
\end{equation*}
{\bf (ii)} We have for all $h\in\mathbb{H}$,
\[
\mathcal{F}\circ\widehat{\pi}(h)=\widehat{\rho}^d(\Phi(h))\circ\mathcal{F}.
\]
\end{thm}

\subsection{Intertwiners and norm formulas}
Define for $0\leq i\leq n$,
\[
I_i^d:=[T_i^d,{}^dX_q^{a_i^d}]\in\mathbb{H}^d,
\]
and $I_\omega^d:=\omega\in\mathbb{H}^d$ ($\omega\in\Omega^d$).
The following theorem extends results from \cite{CI,Sa,St}.
\begin{thm}
For a reduced expression $w=\omega s_{i_1}^ds_{i_2}^d\cdots s_{i_r}^d\in 
W^d$ ($\omega\in\Omega^d$, $0\leq i_j\leq n$), the expression
\[
I_w^d:=I_{\omega}^dI_{i_1}^dI_{i_2}^d\cdots I_{i_r}^d\in\mathbb{H}^d
\]
is well defined (independent of the choice of reduced expression).
In addition,
\[
\widehat{\pi}^d(I_w^d):=r_w^d\cdot w_q\in
\textup{End}_{\mathbb{C}}(\mathbb{C}[T_{\Lambda^d}]),
\]
where $r_w^d:=\prod_{a\in {}^dR^{\bullet\,+}\cap w({}^dR^{\bullet\,-})}r_a^d
\in\mathbb{C}[T_{\Lambda^d}]$ and 
$r_a^d\in\mathbb{C}[T_{\Lambda^d}]$ is defined by
\[
r_a^d(t):=(\kappa_a^d)^{-1}(1-\kappa_a^d\kappa_{2a}^dt_q^a)
(1+\kappa_a^d(\kappa_{2a}^d)^{-1}t_q^a)
\]
for $t\in T_{\Lambda^d}$. In addition, in $\mathbb{H}^d$,
\[
I_w^dI_{w^{-1}}^d=r_w^d({}^dX)(w_qr_{w^{-1}}^d)({}^dX)
\]
and
\[
I_w^dp({}^dX)=(w_qp)({}^dX)I_w^d
\]
for all $p\in\mathbb{C}[T_{\Lambda^d}]$.
\end{thm}
The $I_w^d\in\mathbb{H}^d$ ($w\in W^d$) are called the dual intertwiners.
The intertwiners are
defined by
\[
\mathcal{I}_w:=\delta^{d}(I_w^d)\in\mathbb{H},\qquad w\in W^d
\]
Then in $\mathbb{H}$,
\[
\mathcal{I}_{w^{-1}}\mathcal{I}_w=r_w^d(Y^{-1})(w_qr_{w^{-1}}^d)(Y^{-1})
\]
and
\[
p(Y^{-1})\mathcal{I}_w=\mathcal{I}_w(w_qp)(Y^{-1})
\]
for all $p\in\mathbb{C}[T_{\Lambda^d}]$. Proposition \ref{DAHAaction}
and Theorem \ref{evthm} give the following result.
\begin{prop}\label{intertwineraction}
{\bf (i)} If $\lambda\in\Lambda$ and $0\leq i\leq n$ satisfy
$s_i^d\lambda\not=\lambda$ then
\[
\widehat{\pi}(\mathcal{I}_i)E(\gamma_{\lambda,q}^{-1};\cdot)=
\gamma_{\lambda,q}^{-a_i^d}r_{a_i^d}^d(\gamma_{\lambda,q})
E(\gamma_{s_i^d\lambda,q}^{-1};\cdot).
\]
{\bf (ii)} If $\lambda\in\Lambda$ and $0\leq i\leq n$ satisfy
$s_i^d\lambda=\lambda$ then
\[
\widehat{\pi}(\mathcal{I}_i)E(\gamma_{\lambda,q}^{-1};\cdot)=0.
\]
{\bf (iii)} If $\omega\in\Omega^d$ and $\lambda\in\Lambda$ then
\[
\widehat{\pi}(\mathcal{I}_{\omega})E(\gamma_{\lambda,q}^{-1};\cdot)=
E(\gamma_{\omega^{-1}\lambda,q}^{-1};\cdot).
\]
\end{prop}
This proposition shows that intertwiners can be used to create
the nonsymmetric Macdonald-Koornwinder polynomial from the the constant
polynomial $E(\gamma_{0,q}^{-1};\cdot)\equiv 1$ (cf., e.g.,
\cite[\S 5.10]{M}). We explore now this observation to express
the norms of the nonsymmetric Macdonald-Koornwinder in terms of 
the nonzero constant term 
\[
\langle 1,1\rangle=\int_{T_\Lambda^u}v(t)dt=
\frac{\mathcal{C}(\gamma_{0,q}^d)}{\#W_0}\int_{T_\Lambda^u}v_+(t)dt.
\]
The constant term $\langle 1,1\rangle$
is a $q$-analog and root system generalization of the Selberg integral,
whose explicit evaluation was conjectured by Macdonald \cite{Mpol}
in case $(\Lambda,\Lambda^d)=(P(R_0),P(R_0^d))$. 
By various methods it was evaluated in special cases (for references
we refer to the detailed discussions in \cite{Mpol,CAnn}). A uniform proof
in case $(\Lambda,\Lambda^d)=(P(R_0),P(R_0^d))$ 
using shift operators was given in \cite[Thm. 0.1]{CAnn} (see \cite{St}
for the $C^\vee C$ case).

Write for $\lambda\in\Lambda$,
\begin{equation}\label{quadraticnorm}
N(\lambda):=\frac{\langle E(\gamma_{\lambda,q}^{-1};\cdot),
E^\circ(\gamma_{\lambda,q};\cdot)\rangle}{\langle 1,1\rangle}.
\end{equation}
We write
\begin{equation*}
\begin{split}
c_w^{d}&:=\prod_{a\in R^{d\bullet,+}\cap w^{-1}(R^{d\bullet,-})}
c_a^{d}\in\mathbb{C}[T_{\Lambda^d}],\\
c_{-w}^d&:=\prod_{a\in R^{d\bullet,+}\cap w^{-1}(R^{d\bullet,-})}
c_{-a}^d\in\mathbb{C}[T_{\Lambda^d}]
\end{split}
\end{equation*}
for $w\in W^d$ (warning: $c_{-w}^d(t)=c_w^d(t^{-1})$ is only
valid if $w\in W_0$). Part {\bf (i)} of the following theorem
should be compared with \cite[Prop. 3.4.1]{C} and \cite[(5.2.11)]{M}
(it originates from \cite{Cns} in case $(\Lambda,\Lambda^d)=
(P(R_0),P(R_0^d))$ and \cite{Sa,St} in the $C^\vee C$ case). The first
part of the theorem follows from Proposition \ref{intertwineraction}
and the fact that the 
intertwiners behave nicely with respect to the antiinvolution
$^\ddagger$ (see Lemma \ref{antilem}):
for all $0\leq i\leq n$,
\[
\mathcal{I}_i^\ddagger=\mathcal{I}_i
\]
in $\mathbb{H}(\kappa^{-1},q^{-1})$ and for all $\omega\in\Omega^d$,
\[
\mathcal{I}_{\omega}^\ddagger=\mathcal{I}_{\omega^{-1}}
\]
in $\mathbb{H}(\kappa^{-1},q^{-1})$. The second
part of the theorem is immediate from the first and from the biorthogonality
of the nonsymmetric Macdonald-Koornwinder 
polynomials (see Theorem \ref{biortho}).
\begin{thm}
{\bf (i)} For $\lambda\in\Lambda$,
\[
N(\lambda)=
\frac{c_{-u^d(\lambda)}^d(\gamma_{0,q})}{c_{u^d(\lambda)}^d(\gamma_{0,q})},
\]
which is nonzero for all $\lambda\in\Lambda$.\\
{\bf (ii)} The transform 
$\mathcal{F}: \mathbb{C}[T_\Lambda]\rightarrow F(\mathcal{S})$ is a linear
bijection with inverse $\mathcal{G}: F(\mathcal{S})\rightarrow
\mathbb{C}[T_\Lambda]$ given by
\[
(\mathcal{G}f)(t):=\frac{1}{\langle 1,1\rangle}
\sum_{\lambda\in\Lambda}N(\lambda)^{-1}
f(\gamma_{\lambda,q})E(\gamma_{\lambda,q}^{-1};t)
\]
for $f\in F(\mathcal{S})$ and $t\in T_\Lambda$.
\end{thm}
Combined with the evaluation formula (Theorem \ref{evthm}) we get
\begin{cor}
For all $\lambda\in\Lambda$,
\[
\frac{\langle P_\lambda,P_\lambda^\circ\rangle}{\langle 1,1\rangle}=
\bigl(\kappa_{u^d(\lambda)}^d\bigr)^2c_{u^d(\lambda)}^d(\gamma_{0,q})
c_{-u^d(\lambda)}^d(\gamma_{0,q}).
\]
\end{cor}

\subsection{Normalized symmetric Macdonald-Koornwinder polynomials}
We keep the same assumptions on $q$ and on $\kappa\in\mathcal{M}(D)$ as in
the previous subsection. The results in this subsection are from 
\cite{Cns} in case $(\Lambda,\Lambda^d)=(P(R_0),P(R_0^d))$ and from
\cite{Sa,St} in the $C^\vee C$ case.

For $\lambda\in\Lambda^-$ define
\[
E^+(\gamma_{\lambda,q}^{-1};\cdot)\in\mathbb{C}[T_\Lambda]^{W_0}
\]
by 
\[
E^+(\gamma_{\lambda,q}^{-1};\cdot):=\widehat{\pi}_{\kappa,q}(C_+)
E(\gamma_{\lambda,q}^{-1};\cdot),
\]
where (recall)
\[
C_+:=\frac{1}{\sum_{w\in W_0}\kappa_w^2}\sum_{w\in W_0}\kappa_wT_w.
\]
We call $E^+(\gamma_{\lambda,q}^{-1};\cdot)$ the normalized symmetric
Macdonald-Koornwinder polynomial of degree $\lambda\in\Lambda^-$. 
\begin{lem}
{\bf (i)} For all $w\in W_0$ and $\lambda\in\Lambda^-$,
\begin{equation}\label{welldefEplus}
E^+(\gamma_{\lambda,q}^{-1};\cdot)=\widehat{\pi}(C_+)E(\gamma_{w\lambda,q}^{-1};
\cdot).
\end{equation}
{\bf (ii)} $E^+(\gamma_{\lambda,q}^{-1};\gamma_{0,q}^d)=1$ for all
$\lambda\in\Lambda^-$.\\
{\bf (iii)} $\{E^+(\gamma_{\lambda,q}^{-1};\cdot)\}_{\lambda\in\Lambda^-}$
is a basis of $\mathbb{C}[T_\Lambda]^{W_0}$ satisfying, for all
$p\in\mathbb{C}[T_{\Lambda^d}]^{W_0}$,
\[
D_p(E^+(\gamma_{\lambda,q}^{-1};\cdot))=p(\gamma_{\lambda,q}^{-1})
E^+(\gamma_{\lambda,q}^{-1};\cdot)=
p(q^{-\lambda}\gamma_{0,q}^{-1})E^+(\gamma_{\lambda,q}^{-1};\cdot).
\]
{\bf (iv)} For all $\lambda\in\Lambda^-$, 
\[
E^+(\gamma_{\lambda,q}^{-1};\cdot)=\frac{P_{\lambda_+}^+(\cdot)}
{P_{\lambda_+}^+(\gamma_{0,q}^d)}.
\]
{\bf (iv)} For all $\lambda\in\Lambda^-$ and $\xi\in\Lambda^{d-}$,
\[
E^+(\gamma_{\lambda,q}^{-1};\gamma_{\xi,q}^d)=
E^{+\,d}(\gamma_{\xi,q}^{d\,-1};\gamma_{\lambda,q}).
\]
\end{lem}

As before we 
write superindex $\circ$ to indicate that the parameters $(\kappa,q)$
are inverted. The nonsymmetric and symmetric Macdonald-Koornwinder 
polynomials with inverted parameters $(\kappa^{-1},q^{-1})$
can be explicitly expressed in terms of 
the ones with parameters $(\kappa,q)$. The result is as follows
(see \cite[\S 3.3.2]{C}
for $(\Lambda,\Lambda^d)=(P(R_0),P(R_0^d))$ and \cite[(2.5)]{StPr}
in the twisted case). 
\begin{prop}
{\bf (i)} For all $\lambda\in\Lambda$,
\[
E^\circ(\gamma_{\lambda,q};t^{-1})=\kappa_{w_0}^{-1}
\bigl(\widehat{\pi}_{\kappa,q}(T_{w_0})E(\gamma_{-w_0\lambda,q}^{-1};\cdot)\bigr)(t)
\]
in $\mathbb{C}[T_\Lambda]$,
where $w_0\in W_0$ is the longest Weyl group element and $t\in T_\Lambda$.\\
{\bf (ii)} For all $\lambda\in\Lambda^-$ we have
\begin{equation*}
\begin{split}
E^{+\,\circ}(\gamma_{\lambda,q};t^{-1})&=E^+(\gamma_{-w_0\lambda,q}^{-1};t),\\
E^{+\,\circ}(\gamma_{\lambda,q};t)&=E^+(\gamma_{\lambda,q}^{-1};t)
\end{split}
\end{equation*}
in $\mathbb{C}[T_\Lambda]^{W_0}$.
\end{prop}
Using intertwiners or Proposition \ref{DAHAaction} it is now
possible to expand the normalized symmetric Macdonald-Koornwinder polynomials
in nonsymmetric ones. This in turn leads to an explicit expression of 
the quadratic norms of the symmetric Macdonald-Koornwinder polynomials in
terms of the nonsymmetric ones. Recall the rational function 
$\mathcal{C}(\cdot)=\mathcal{C}(\cdot;D;\kappa,q)\in\mathbb{C}(T_\Lambda)$
from \eqref{mathcalC}. We write
$\mathcal{C}^d(\cdot)=\mathcal{C}(\cdot;D^d;\kappa^d,q)\in
\mathbb{C}(T_{\Lambda^d})$ for its dual version.
\begin{thm}
{\bf (i)} For $\lambda\in\Lambda^-$ we have
\begin{equation}\label{expansionP}
P_{\lambda_+}^+(t)=\sum_{\mu\in W_0\lambda}\left(
\prod_{\alpha\in R_0^+\cap v(\mu)R_0^-}c_\alpha^d(\gamma_{\lambda,q})
\right)P_\mu(t).
\end{equation}
{\bf (ii)} For $\lambda\in\Lambda^-$ we have
\[
E^+(\gamma_{\lambda,q}^{-1};t)=
\sum_{\mu\in W_0\lambda}\frac{\mathcal{C}^d(\gamma_{\mu,q})}
{\mathcal{C}^d(\gamma_{0,q})}E(\gamma_{\mu,q}^{-1};t).
\]
{\bf (iii)} For $\lambda,\mu\in\Lambda^-$ we have
\[
\frac{\langle E^+(\gamma_{\lambda,q}^{-1};\cdot),
E^+(\gamma_{\mu,q}^{-1};\cdot)\rangle_+}{\langle 1,1\rangle_+}
=\delta_{\lambda,\mu}
\frac{\mathcal{C}^d(\gamma_{\lambda,q})N(\lambda)}
{\mathcal{C}^d(\gamma_{0,q})}
\]
where $\langle p,r\rangle_+:=\int_{T_\Lambda^u}p(t)\overline{r(t)}v_+(t)dt$.
\end{thm}

As a consequence we get the following explicit evaluation formulas
and quadratic norm formulas
for the monic symmetric Macdonald-Koornwinder polynomials.
Set 
\[
N^+(\lambda):=\frac{\langle E^+(\gamma_{\lambda,q}^{-1};\cdot),
E^+(\gamma_{\lambda,q}^{-1};\cdot)\rangle_+}{\langle 1,1\rangle_+},
\qquad \lambda\in\Lambda^-.
\]

\begin{cor}
{\bf (i)} For $\lambda\in\Lambda^-$ we have
\[
P_{\lambda_+}^+(\gamma_{0,q}^d)=
\frac{\mathcal{C}^d(\gamma_{0,q})c_{\tau(\lambda)}^d(\gamma_{0,q})}
{\mathcal{C}^d(\gamma_{\lambda,q})\kappa_{\tau(\lambda)}^d}.
\]
{\bf (ii)}
For $\lambda\in\Lambda^-$,
\[
N^+(\lambda)=\frac{\mathcal{C}^d(\gamma_{\lambda,q})
c_{-\tau(\lambda)}^d(\gamma_{0,q})}
{\mathcal{C}^d(\gamma_{0,q})c_{\tau(\lambda)}^d(\gamma_{0,q})}.
\]
\end{cor}

\begin{rema}
For $\lambda\in\Lambda^-$,
\[
\sum_{\mu\in W_0\lambda}N(\mu)^{-1}=N^+(\lambda)^{-1}.
\]
\end{rema}

\section{Explicit evaluation and norm formulas}\label{5}
We rewrite in this subsection the explicit 
evaluation formulas and the
quadratic norm espressions for the symmetric Macdonald-Koornwinder 
polynomials in terms of $q$-shifted factorials \eqref{qshifted}. The
explicit formulas for the $\textup{GL}_{n+1}$ symmetric
Macdonald polynomials (see Subsection \ref{GLn}) can be immediately obtained
as special cases of the explicit formulas below. We keep the conditions
\eqref{smallparameterorth} on the parameters $(\kappa,q)$.

\subsection{The twisted cases}
In case we have initial data $D=(R_0,\Delta_0,\bullet,\Lambda,\Lambda^d)$
with $\bullet=t$ the evaluation and norm formulas take the following
explicit form.
\begin{cor}
{\bf (i)} Suppose $\bullet=t$ and $S=\emptyset=S^d$. Then
$R=R^\bullet=R^d$, $\kappa=\kappa^d$ and 
$\kappa_{m\mu_\alpha c+\alpha}=\kappa_\alpha$ 
for all $m\in\mathbb{Z}$ and $\alpha\in R_0$.
Then for all $\lambda\in\Lambda^-$,
\begin{equation*}
\begin{split}
P_{\lambda_+}^+(\gamma_{0,q})&=\gamma_{0,q}^{-\lambda}
\prod_{\alpha\in R_0^+}
\frac{\bigl(\kappa_\alpha^2\gamma_{0,q}^{-\alpha};q_\alpha\bigr)_{
-(\lambda,\alpha^\vee)}}
{\bigl(\gamma_{0,q}^{-\alpha};q_\alpha\bigr)_{-(\lambda,\alpha^\vee)}},\\
N^+(\lambda)&=\gamma_{0,q}^{2\lambda}\prod_{\alpha\in R_0^+}
\left(
\frac{1-\gamma_{0,q}^{-\alpha}}{1-q_\alpha^{-(\lambda,\alpha^\vee)}
\gamma_{0,q}^{-\alpha}}\right)
\frac{\bigl(q_\alpha \kappa_\alpha^{-2}\gamma_{0,q}^{-\alpha};
q_\alpha\bigr)_{-(\lambda,\alpha^\vee)}}
{\bigl(\kappa_\alpha^2\gamma_{0,q}^{-\alpha};q_\alpha
\bigr)_{-(\lambda,\alpha^\vee)}}.
\end{split}
\end{equation*}
{\bf (ii)} Suppose $\bullet=t$ and $(\Lambda,\Lambda^d)=(\mathbb{Z}R_0,
\mathbb{Z}R_0)$ (the $C^\vee C$ case), realized concretely
as in Subsection \ref{CcheckC}. 
Recall the relabelling of the multiplicity functions $\kappa$
and $\kappa^d$ given by
$k=\kappa_{\varphi}^2=\widetilde{k}$ and 
\begin{equation*}
\begin{split}
\{a,b,c,d\}&=\{\kappa_\theta\kappa_{2\theta}, -\kappa_\theta\kappa_{2\theta}^{-1},
q_\theta\kappa_0\kappa_{2a_0}, -q_\theta\kappa_0\kappa_{2a_0}^{-1}\},\\
\{\widetilde{a},\widetilde{b},\widetilde{c},\widetilde{d}\}&=
\{\kappa_\theta\kappa_0,-\kappa_\theta\kappa_0^{-1}, 
q_\theta\kappa_{2\theta}\kappa_{2a_0},
-q_\theta\kappa_{2\theta}\kappa_{2a_0}^{-1}\}.
\end{split}
\end{equation*}
Then for $\lambda\in\Lambda^-$ we have
\begin{equation*}
\begin{split}
P_{\lambda_+}^+(\gamma_{0,q}^d)=(\gamma_{0,q}^d)^{-\lambda}
&\prod_{\alpha\in R_{0,l}^+}
\frac{\bigl(\kappa_\varphi^2\gamma_{0,q}^{-\alpha};q_\varphi
\bigr)_{-(\lambda,\alpha^\vee)}}
{\bigl(\gamma_{0,q}^{-\alpha};q_\varphi\bigr)_{-(\lambda,\alpha^\vee)}}\\
&\times\prod_{\alpha\in R_{0,s}^+}
\frac{\bigl(\widetilde{a}\gamma_{0,q}^{-\alpha},
\widetilde{b}\gamma_{0,q}^{-\alpha},
\widetilde{c}\gamma_{0,q}^{-\alpha},\widetilde{d}\gamma_{0,q}^{-\alpha};
q_\theta^2\bigr)_{-(\lambda,\alpha^\vee)/2}}
{\bigl(\gamma_{0,q}^{-2\alpha};q_\theta^2\bigr)_{-(\lambda,\alpha^\vee)}}
\end{split}
\end{equation*}
and
\begin{equation*}
\begin{split}
N^+(\lambda)&=(\gamma_{0,q}^d)^{2\lambda}
\prod_{\alpha\in R_{0,l}^+}
\left(\frac{1-\gamma_{0,q}^{-\alpha}}{1-q_\varphi^{-(\lambda,\alpha^\vee)}
\gamma_{0,q}^{-\alpha}}\right)
\frac{\bigl(
q_\varphi k_\varphi^{-2}\gamma_{0,q}^{-\alpha};q_\varphi\bigr)_{-(\lambda,\alpha^\vee)}}
{\bigl(k_\varphi^2\gamma_{0,q}^{-\alpha};q_\varphi
\bigr)_{-(\lambda,\alpha^\vee)}}\\
&\times\prod_{\alpha\in R_{0,s}^+}
\left(\frac{1-\gamma_{0,q}^{-2\alpha}}
{1-q_\theta^{-2(\lambda,\alpha^\vee)}\gamma_{0,q}^{-2\alpha}}\right)
\frac{\bigl(q_\theta^2\widetilde{a}^{-1}\gamma_{0,q}^{-\alpha},
q_\theta^2\widetilde{b}^{-1}\gamma_{0,q}^{-\alpha},
q_\theta^2\widetilde{c}^{-1}\gamma_{0,q}^{-\alpha},
q_\theta^2\widetilde{d}^{-1}\gamma_{0,q}^{-\alpha};
q_\theta^2\bigr)_{-(\lambda,\alpha^\vee)/2}}
{\bigl(\widetilde{a}\gamma_{0,q}^{-\alpha},\widetilde{b}\gamma_{0,q}^{-\alpha},
\widetilde{c}\gamma_{0,q}^{-\alpha},\widetilde{d}\gamma_{0,q}^{-\alpha};q_\theta^2
\bigr)_{-(\lambda,\alpha^\vee)/2}},
\end{split}
\end{equation*}
where $R_{0,s}\subset R_0$ (respectively $R_{0,l}\subset R_0$)
are the short (respectively long) roots in $R_0$. In addition,
$q_\varphi=q_\theta^2$.
\end{cor}
The formulas in the intermediate case $\bullet=t$ and $S=\emptyset\not=S^d$
(or $S\not=\emptyset=S^d$) are special cases of the $C^\vee C$ case
(by choosing an appropriate specialization of the multiplicity function).

\subsection{The untwisted cases}
If the initial data is of the form
$D=(R_0,\Delta_0,u,\Lambda,\Lambda^d)$ then 
there are essentially two cases to consider,
namely the case that $S=\emptyset=S^d$ and the exceptional rank two
case (see Subsection \ref{Except2}). Indeed,
the cases corresponding to a nonreduced
extension of the untwisted affine root system with underlying finite root
system $R_0$ of type $\textup{B}_n$ ($n\geq 3)$ or $\textup{BC}_n$ ($n\geq 1$)
are special cases of the twisted $C^\vee C_n$ case. 

\begin{cor}
{\bf (i)} Suppose $\bullet=u$ and $S=\emptyset=S^d$. Then
$R=\mathbb{Z}c+R_0$, $R^d=\mathbb{Z}c+R_0^\vee$ and 
$\kappa_{mc+\alpha}=\kappa_\alpha=\kappa_{mc+\alpha^\vee}^d$
for all $m\in\mathbb{Z}$ and $\alpha\in R_0$.
Then for all $\lambda\in\Lambda^-$,
\begin{equation*}
\begin{split}
P_{\lambda_+}^+(\gamma_{0,q}^d)&=\gamma_{0,q}^{-\lambda}
\prod_{\alpha\in R_0^+}
\frac{\bigl(\kappa_\alpha^2\gamma_{0,q}^{-\alpha^\vee};q\bigr)_{
-(\lambda,\alpha^\vee)}}
{\bigl(\gamma_{0,q}^{-\alpha^\vee};q\bigr)_{-(\lambda,\alpha^\vee)}},\\
N^+(\lambda)&=\gamma_{0,q}^{2\lambda}\prod_{\alpha\in R_0^+}
\left(\frac{1-\gamma_{0,q}^{-\alpha^\vee}}{1-q^{-(\lambda,\alpha^\vee)}
\gamma_{0,q}^{-\alpha^\vee}}\right)
\frac{\bigl(q\kappa_\alpha^{-2}\gamma_{0,q}^{-\alpha^\vee};
q\bigr)_{-(\lambda,\alpha^\vee)}}
{\bigl(\kappa_\alpha^2\gamma_{0,q}^{-\alpha^\vee};
q\bigr)_{-(\lambda,\alpha^\vee)}}.
\end{split}
\end{equation*}
{\bf (ii)} Suppose $D=(R_0,\Delta_0,u,\mathbb{Z}R_0,\mathbb{Z}R_0^\vee)$
with $R_0$ of type $\textup{C}_2$, realized concretely as in Subsection
\ref{Except2}. 
Recall the relabelling of the multiplicity functions $\kappa$
and $\kappa^d$ given by \eqref{parametersExcept2}.
Then for $\lambda\in\Lambda^-$ we have
\[
P_{\lambda_+}^+(\gamma_{0,q}^d)=(\gamma_{0,q}^d)^{-\lambda}
\prod_{\alpha\in R_{0,s}^+}\frac{\bigl(\widetilde{c}\gamma_{0,q}^{-\alpha^\vee},
\widetilde{d}\gamma_{0,q}^{-\alpha^\vee};q^2\bigr)_{-(\lambda,\alpha^\vee)/2}}
{\bigl(\gamma_{0,q}^{-\alpha^\vee};q\bigr)_{-(\lambda,\alpha^\vee)}}
\prod_{\beta\in R_{0,l}^+}\frac{\bigl(\widetilde{a}\gamma_{0,q}^{-\beta^\vee},
\widetilde{b}\gamma_{0,q}^{-\beta^\vee};q\bigr)_{-(\lambda,\beta^\vee)}}
{\bigl(\gamma_{0,q}^{-2\beta^\vee};q^2\bigr)_{-(\lambda,\beta^\vee)}}
\]
and
\begin{equation*}
\begin{split}
N^+(\lambda)&=(\gamma_{0,q}^d)^{2\lambda}
\prod_{\alpha\in R_{0,s}^+}
\left(\frac{1-\gamma_{0,q}^{-\alpha^\vee}}{1-q^{-(\lambda,\alpha^\vee)}
\gamma_{0,q}^{-\alpha^\vee}}\right)
\frac{\bigl(q^2\widetilde{c}^{-1}\gamma_{0,q}^{-\alpha^\vee},
q^2\widetilde{d}^{-1}\gamma_{0,q}^{-\alpha^\vee};
q^2\bigr)_{-(\lambda,\alpha^\vee)/2}}
{\bigl(\widetilde{c}\gamma_{0,q}^{-\alpha^\vee},
\widetilde{d}\gamma_{0,q}^{-\alpha^\vee};q^2\bigr)_{-(\lambda,\alpha^\vee)/2}}\\
&\times
\prod_{\beta\in R_{0,l}^+}\left(
\frac{1-\gamma_{0,q}^{-2\beta^\vee}}{1-q^{-2(\lambda,\beta^\vee)}
\gamma_{0,q}^{-2\beta^\vee}}\right)
\frac{\bigl(q\widetilde{a}^{-1}\gamma_{0,q}^{-\beta^\vee},
q\widetilde{b}^{-1}\gamma_{0,q}^{-\beta^\vee};q\bigr)_{-(\lambda,\beta^\vee)}}
{\bigl(\widetilde{a}\gamma_{0,q}^{-\beta^\vee},
\widetilde{b}\gamma_{0,q}^{-\beta^\vee};q\bigr)_{-(\lambda,\beta^\vee)}},
\end{split}
\end{equation*}
where $R_{0,s}\subset R_0$ (respectively $R_{0,l}\subset R_0$)
are the short (respectively long) roots in $R_0$.
\end{cor}

\section{Appendix on affine root systems}\label{1}
\subsection{Affine root systems}
The main reference for this subsection is Macdonald \cite{M0}.

Let $E$ be a real affine space with the associated space of translations
$V$ of dimension $n\geq 1$. Fix a real scalar product $\bigl(\cdot,\cdot\bigr)$
on $V$ and set $|v|^2=(v,v)$ for all $v\in V$. It turns $E$ into a metric
space, called an affine Euclidean space \cite[\S 3]{BT}.

A map $\Psi: E\rightarrow E$ is called an affine linear endomorphism of $E$
if there exists a linear endomorphism $d\Psi$ of $V$ such that
$\Psi(e+v)=\Psi(e)+d\Psi(v)$ for all $e\in E$ and all $v\in V$.
Set $\textup{O}(E)$ for 
the group of affine linear isometric automorphisms of $E$.
Let $\tau_E: V\rightarrow \textup{O}(E)$ be the group
monomorphism defined by $\tau_E(v)(e)=e+v$. 
 Then we have a short exact
sequence of groups
\[
\tau_E(V)\hookrightarrow
\textup{O}(E)\overset{d}{\twoheadrightarrow}
\textup{O}(V).
\]
For a subgroup $W\subseteq\textup{O}(E)$
let $L_W\subseteq V$ be the additive subgroup such that
$\tau_E(L_W)=\textup{Ker}(d|_W)$.

Set $\widehat{E}$ for the real $n+1$-dimensional 
vector space
of affine linear functions $a:E\rightarrow\mathbb{R}$. 
Let $c\in\widehat{E}$ be the constant function one.
The gradient of an affine 
linear function $a\in\widehat{E}$
is the unique vector $Da\in V$
such that $a(e+v)=a(e)+(Da,v)$ for all $e\in E$ and
$v\in V$. The gradient map
$D: \widehat{E}\rightarrow V$ is linear, surjective, with kernel 
consisting of the constant functions on $E$. 

For $a,b\in\widehat{E}$ set
\[
\bigl(a,b\bigr):=\bigl(Da,Db\bigr),
\]
which defines a semi-positive definite symmetric bilinear form on
$\widehat{E}$. The radical consists of the constant functions on $E$. We write
$|a|^2:=(a,a)$ for $a\in\widehat{E}$. 
Let $O(\widehat{E})$ be the form-preserving linear
automorphisms of $\widehat{E}$ and $O_c(\widehat{E})$ its subgroup 
of automorphisms fixing the constant functions. 

The contragredient action of $g\in O(E)$ on $\widehat{E}$, 
given by $(ga)(e):=a(g^{-1}e)$ 
($a\in\widehat{E}$, $e\in E$), realizes a group isomorphism
$O(E)\simeq O_c(\widehat{E})$. Note that
$\tau_E(v)a=a-(Da,v)c$ for $a\in \widehat{E}$ and $v\in V$.

A vector $a\in\widehat{E}$ is called nonisotropic if $Da\not=0$.
For a nonisotropic vector $a\in\widehat{E}$ let
$s_a: E\rightarrow E$ be the
orthogonal reflection in the affine hyperplane $a^{-1}(0)$ of $E$.
It is explicitly given by
\[
s_a(e)=e-a(e)Da^\vee,\qquad e\in E,
\]
where $v^\vee:=2v/|v|^2\in V$ is the covector of $v\in V\setminus\{0\}$. 
Viewed as element of $O_c(\widehat{E})$ it reads
\[
s_a(b)=b-(a^\vee,b)a,\qquad b\in \widehat{E},
\]
where $a^\vee:=2a/|a|^2\in\widehat{E}$ is the covector of $a$. 
For a subset $R$ of nonisotropic vectors in $\widehat{E}$ 
let $W(R)$ be the subgroup of 
$O(E)\simeq O_c(\widehat{E})$ generated by the orthogonal reflections
$s_a$ ($a\in R$). 
\begin{defi}\label{ars}
A set $R$ of nonisotropic vectors in $\widehat{E}$ is an affine root
system on $E$ if 
\begin{enumerate}
\item[(1)] $R$ spans $\widehat{E}$,
\item[(2)] $W(R)$ stabilizes $R$,
\item[(3)] $(a^\vee,b)\in\mathbb{Z}$ for all $a,b\in R$,
\item[(4)] $W(R)$ acts properly on $E$ (i.e., if 
$K_1$ and $K_2$ are two compact subsets of $E$ then
$w(K_1)\cap K_2\not=\emptyset$ for at most finitely many $w\in W(R)$),
\item[(5)] $L_{W(R)}$ spans $V$.
\end{enumerate}
\end{defi}
The elements $a\in R$ are called affine roots. The group $W(R)$ is
called the affine Weyl group of $R$. The real dimension $n$ of $V$
is the rank of $R$.

\begin{defi}
Let $R$ be an affine root system on $E$. A nonempty subset $R^\prime\subseteq R$
is called an affine root subsystem if $W(R^\prime)$ stabilizes $R^\prime$
and if $L_{W(R^\prime)}$ generates the real span $V^\prime$ of the set
$\{Da\}_{a\in R^\prime}$ of gradients of $R^\prime$ in $V$.
\end{defi}
\begin{rema}
With the notations from the previous definitions, let $E^\prime$ be the
$V^{\prime\perp}$-orbits of $E$, where $V^{\prime\perp}$ is the orthocomplement
of $V^\prime$ in $V$. It is an affine Euclidean space with $V^\prime$ 
the associated space
of translations and with norm induced by the scalar product on 
$V^\prime$. 
Let $F\subseteq\widehat{E}$ be the real span of $R^\prime$. 
Then $F\overset{\sim}{\rightarrow}
\widehat{E^\prime}$ with form-preserving linear
isomorphism $a\mapsto a^\prime$ defined by $a^\prime(e+V^{\prime\perp}):=a(e)$
for all $a\in F$ and $e\in E$. With this identification $R^\prime$
is an affine root system on $E^\prime$. Furthermore,
the corresponding affine Weyl group
$W(R^\prime)$ is isomorphic to the subgroup of $W(R)\subset
O(E)$ generated by the orthogonal reflections $s_a$ ($a\in R^\prime$).
\end{rema}

We call an affine root system $R$ irreducible if it cannot be written
as a nontrivial orthogonal disjoint union $R^\prime\cup R^{\prime\prime}$ 
(orthogonal meaning that $(a,b)=0$
for all $a\in R^\prime$ and all $b\in R^{\prime\prime}$). 
It is called reducible otherwise. In that case both $R^\prime$
and $R^{\prime\prime}$ are affine root subsystems of $R$.
Each affine root system is an orthogonal disjoint union of irreducible
affine root subsystems
(cf. \cite[\S 3]{M0}). 
\begin{rema}
Macdonald's \cite[\S 2]{M0} definition of an affine root system 
is (1)-(4) of Definition \ref{ars}. Careful analysis 
reveals that Macdonald tacitly
assumes condition (5), which only follows from the four axioms (1)-(4)
if $R$ is irreducible. Mark Reeder 
independently noted that an extra condition besides 
(1)-(4) is needed in order to avoid examples
of affine root systems given as an orthogonal
disjoint union of an affine root system $R^\prime$
and a finite crystallographic 
root system $R^{\prime\prime}$ 
(compare, on the level of affine Weyl groups,
with \cite[Chpt. V, \S 3]{B} and \cite[\S 1.3]{BT}).
Mark Reeder proposes to add to (1)-(4) the axiom that
for each $\alpha\in D(R)$ there exists at least two affine roots with
gradient $\alpha$. The resulting definition 
is equivalent to Definition \ref{ars}, as well as 
to the notion of an {\'e}chelonnage from \cite[(1.4.1)]{BT} (take here
\cite[(1.3.2)]{BT} into account).
\end{rema}

An affine root system $R$ is called reduced if $\mathbb{R}a\cap R=\{\pm a\}$
for all $a\in R$, and nonreduced otherwise.
If $R$ is nonreduced then $R=R^{ind}\cup R^{unm}$ with $R^{ind}$ 
(respectively $R^{unm}$) the reduced affine root subsystem of $R$ consisting
of indivisible (respectively unmultiplyable) affine roots.

If $R\subset \widehat{E}$ is an affine root system then $R_0:=D(R)\subset V$
is a finite crystallographic root system in $V$, called the gradient root
system of $R$. The associated Weyl group $W_0=W_0(R_0)$  is the subgroup
of $\textup{O}(V)$ generated by the orthogonal reflections
$s_\alpha\in O(V)$ in the hyperplanes $\alpha^{\perp}$ ($\alpha\in R_0$), 
which are explicitly given by $s_\alpha(v)=v-(\alpha^\vee,v)\alpha$ for $v\in V$.
The Weyl group $W_0$ coincides with
the image of $W(R)$ under the differential $d$.

We now define an appropriate equivalence relation between affine root systems,
called similarity. 
It is a slightly weaker notion of similarity compared to the one used
in \cite[\S 3]{M0}. This is to render affine root systems
that differ by a rescaling of the underlying gradient root system
similar, see Remark \ref{similaritydifference}.
\begin{defi}
We call two affine root systems $R\subset\widehat{E}$ 
and $R^\prime\subset\widehat{E^\prime}$ similar,
$R\simeq R^\prime$, if there exists a linear isomorphism $T: \widehat{E}
\overset{\sim}{\longrightarrow} \widehat{E^\prime}$ which restricts to 
a bijection of $R$ onto $R^\prime$
preserving Cartan integers,
\[
\bigl((Ta)^\vee,Tb\bigr)=\bigl(a^\vee,b\bigr),\qquad \forall\, a,b\in R.
\]
\end{defi}
Similarity respects basic notions as affine root subsystems and
irreducibility. If $R\simeq R^\prime$, realized by the linear isomorphism $T$, 
then $Ts_aT^{-1}=s_{Ta}$ for all $a\in R$. In particular,
$W(R)\simeq W(R^\prime)$. Note that
$T$ maps constant functions to constant functions.  
Replacing $T$ by $-T$ if necessary,
we may assume without loss of generality that 
$T(c)\in\mathbb{R}_{>0}c^\prime$, 
where $c\in\widehat{E}$ and $c^\prime\in\widehat{E^\prime}$ denote the
constant functions one on $E$ and $E^\prime$ respectively. With this additional
condition we call $T$ a similarity transformation between $R$ and $R^\prime$.

If $R$ is an affine root system and $\lambda\in\mathbb{R}^*:=\mathbb{R}
\setminus\{0\}$ then
$\lambda R:=\{\lambda a\}_{a\in R}$ is an affine root system similar to $R$
(the similarity transformation realizing $R\overset{\sim}{\rightarrow}\lambda R$
is scalar
multiplication by $|\lambda|$). We call $\lambda R$ a rescaling of the
affine root system $R$. If two affine root systems $R$ and $R^\prime$
are similar, then a similarity transformation $T$ between $R$ and
a rescaling of $R^\prime$ exists such that $T(c)=c^\prime$. In this case 
$T$ arises as the contragredient of an affine linear isomorphism
from $E^\prime$ onto $E$. For instance, each affine 
Weyl group element $w\in W(R)$ is a selfsimilarity 
transformation of $R$ in this way.

If the affine root systems $R\subset\widehat{E}$ and 
$R^\prime\subset\widehat{E^\prime}$ are similar, then so are
their gradients $R_0\subset V$ and $R_0^\prime\subset V^\prime$
(i.e. there exists a linear isomorphism $t$ of $V$ onto $V^\prime$
restricting to a bijection of $R_0\overset{\sim}{\rightarrow}
R_0^\prime$ and preserving 
Cartan integers).
Indeed, if $T$ is a similarity transformation between $R$ and
$R^\prime$, then the unique
linear isomorphism $t: V\overset{\sim}{\longrightarrow} 
V^\prime$ such that $D\circ T=t\circ D$ realizes
the similarity between $R_0$ and $R_0^\prime$. 

\begin{rema}\label{similaritydifference}
Let $R\subset\widehat{E}$ be an irreducible affine root system with associated
gradient $R_0\subset V$. Fix an origin $e\in E$.
For $\lambda\in\mathbb{R}^*\setminus\{\pm 1\}$ and $a\in R$ define 
$a_\lambda\in\widehat{E}$ by 
$a_\lambda(e+v):=a(e)+\lambda(Da,v)$ for all $v\in V$.
Then $R_\lambda=\{a_\lambda\}_{a\in R}\subset\widehat{E}$ 
is an affine root system similar to $R$, and $R_{\lambda,0}=\lambda R_0$.
The affine root systems $R$ and $R_\lambda$ are not similar if one
uses Macdonald's \cite[\S 3]{M0} definition of similarity.
\end{rema}

In the remainder of this section we assume that $R$ is an irreducible 
affine root system of rank $n$.
Since $W(R)$ acts properly on $E$, the set of regular elements
\[
E_{reg}:=E\setminus\bigcup_{a\in R}a^{-1}(0)
\]
decomposes as the disjoint union of open $n$-simplices, called chambers
of $R$. For a fixed chamber $C$ there exists a unique
$\mathbb{R}$-basis $\Delta=\Delta(C,R)$ of $\widehat{E}$
consisting of indivisible affine roots $a_i=a_{i,C}$ ($0\leq i\leq n$) such that
$C=\{e\in E \, | \, a_i(e)>0\,\,\, \forall i\in\{0,\ldots,n\}\}$. 
The set of roots $\Delta$ is called the basis of $R$ associated to the chamber
$C$. The affine roots $a_i$ are called simple affine roots.
Any affine root $a\in R$ can be uniquely written as
$a=\sum_{i=0}^n\lambda_ia_i$ with either all $\lambda_i\in\mathbb{Z}_{\geq 0}$
or all $\lambda_i\in\mathbb{Z}_{\leq 0}$. The subset of 
affine roots of the first type is denoted by $R^+$ and is called the
set of positive affine roots with respect to $\Delta$. Then $R=R^+\cup R^-$
(disjoint union) with $R^-:=-R^+$ the subset of negative affine roots.
The affine Weyl group $W(R)$ is a Coxeter group
with Coxeter generators the simple reflections $s_{a_i}$ 
($0\leq i\leq n$).

A rank $n$
{\it affine Cartan matrix} is a rational integral 
$(n+1)\times (n+1)$-matrix $A=(a_{ij})_{i,j=0}^n$ satisfying the four conditions
\begin{enumerate}
\item $a_{ii}=2$,
\item $a_{ij}\in\mathbb{Z}_{\leq 0}$ if $i\not=j$,
\item $a_{ij}=0$ implies $a_{ji}=0$,
\item $\det(A)=0$ and all the proper principal minors of $A$ are strictly
positive.
\item $A$ is indecomposable (i.e. the matrices obtained from $A$ by
a simultaneous permutation of its rows and columns
are not the direct sum of two nontrivial blocks). 
\end{enumerate}

The larger class of rational integral matrices satisfying 
conditions (1)-(3) are called generalized Cartan matrices.
They correspond to Kac-Moody Lie algebras, see \cite{K}. 
The Kac-Moody Lie algebras related to the 
subclass of affine Cartan matrices are the affine Lie algebras. 
The affine Cartan matrices have been classified by Kac \cite[Chpt. 4]{K}. 
 
Fix an ordered basis $\Delta=(a_0,a_1,\ldots,a_n)$ of $R$.
The matrix $A=A(R,\Delta)=(a_{ij})_{0\leq i,j\leq n}$ defined by 
\[
a_{ij}:=\bigl(a_i^\vee,a_j\bigr)
\]
is an affine Cartan matrix.
The coefficients $a_{ij}$ ($0\leq i,j\leq n$) 
are called the affine Cartan integers of $R$. 

If $R\simeq R^\prime$ with associated similarity transformation $T$,
then the $T$-image of an ordered basis $\Delta$ of $R$
is an ordered
basis $\Delta^\prime$ of $R^\prime$. 
Let $R$ and $R^\prime$ be irreducible affine root systems with ordered
bases $\Delta$ and $\Delta^\prime$ respectively.
We say that $(R,\Delta)$ 
is similar to $(R^\prime,\Delta^\prime)$,
$(R,\Delta)\simeq (R^\prime,\Delta^\prime)$,
if the irreducible affine root
systems $R$ and $R^\prime$ are similar and if 
there exists an associated similarity 
transformation $T$
mapping the ordered basis $\Delta$ of $R$
to the ordered basis $\Delta^\prime$ of $R^\prime$.
For similar pairs the associated 
affine Cartan matrices coincide.
Moreover, the affine Cartan matrix $A(R,\Delta)$
modulo simultaneous permutations of its row and columns 
does not depend on the choice of ordered basis $\Delta$. 
This leads to a map from the set of similarity classes of reduced irreducible
affine root systems to the set of affine Cartan matrices
up to simultaneous permutations of rows and columns. 
Using Kac' \cite[Chpt. 4]{K} classification of the affine
Cartan matrices and the explicit construction of reduced irreducible
affine root systems from \cite{M0} (see also the next subsection),
it follows that the map is 
surjective. It is injective by a straightforward adjustment of the proof
of \cite[Prop. 11.1]{Hu} to the present affine setup. Hence we obtain
the following classification result.

\begin{thm}\label{classification}
Reduced irreducible affine root systems up to similarity
are in bijective correspondence to affine 
Cartan matrices up to simultaneous
permutations of the rows and columns.
\end{thm}
\begin{rema}
Affine Cartan matrices also parametrize affine Lie algebras (see \cite{K}).
For a given affine Cartan matrix the
set of real roots of the associated affine Lie algebra is the
associated irreducible reduced affine root system. 
\end{rema}
Affine Cartan matrices up to simultaneous permutations of the
rows and columns (and hence similarity classes of irreducible
reduced affine root systems) 
can be naturally encoded by affine Dynkin diagrams
\cite{K,M0}. The affine Dynkin diagram associated to
an affine Cartan matrix $A=(a_{ij})_{i,j=0}^n$ is the graph with $n+1$
vertices in which we join the $i$th and $j$th node ($i\not=j$) 
by $a_{ij}a_{ji}$ edges. In addition we put an
arrow towards the $i$th node if $|a_{ij}|>1$. At the end of the
appendix we list all affine Dynkin diagrams and link it to Kac's
\cite{K} classification.

\subsection{Explicit constructions}\label{special}
The main reference for the results in this subsection is \cite{M0}.
For an irreducible finite crystallographic root system $R_0$ of type
$A,D,E$ or $BC$ there is exactly one similarity class of reduced irreducible
affine root systems whose gradient root system is similar to $R_0$.
For the other types of root systems $R_0$, there are two such similarity
classes of reduced irreducible affine root systems.
We now proceed to realize them explicitly.

Let $R_0\subset V$ be an irreducible finite crystallographic root system
(possibly nonreduced). 
The affine space $E$ is taken to be $V$ with forgotten origin. We will
write $V$ for $E$ in the sequel if no confusion is possible.

We identify the space 
$\widehat{E}$ of affine linear functions on $E$
with $V\oplus\mathbb{R}c$ as real vector space,
with $c$ the 
constant function identically equal to one on $V$ and with $V^*\simeq V$
the linear functionals on $V$ (the identification with $V$ is realized
by the scalar product on $V$). With these identifications,
\[
O_c(\widehat{E})\simeq O(E)=O(V)\ltimes\tau(V),
\]
with
$\tau(v)(e)=e+v$.
Regarding $\tau(v)$ as element of $O_c(\widehat{V})$ it is given
by
\[
\tau(v)a=-(Da,v)c+a,\qquad a\in\widehat{V}.
\]
Note that the orthogonal reflection $s_a\in O(E)$ associated to
$a=\lambda c+\alpha\in\widehat{E}$ ($\lambda\in\mathbb{R}$,
$\alpha\in V\setminus\{0\}$) decomposes as
$s_{\lambda c+\alpha}=\tau(-\lambda\alpha^\vee)s_\alpha$.

Consider the subset 
\[
\mathcal{S}(R_0):=\{mc+\alpha\}_{m\in\mathbb{Z}, \alpha\in R_0^{ind}}\cup
\{(2m+1)c+\beta\}_{m\in\mathbb{Z}, \beta\in R_0\setminus R_0^{ind}}
\]
of $\widehat{V}$, where $R_0^{ind}\subseteq R_0$ is the root subsystem
of indivisible roots. Then $\mathcal{S}(R_0)$ and $\mathcal{S}(R_0^\vee)^\vee$ 
are reduced irreducible affine root systems with gradient root system
$R_0$. We call $\mathcal{S}(R_0)$ (respectively $\mathcal{S}(R_0^\vee)^\vee$) 
the untwisted (respectively twisted) 
reduced irreducible affine root system associated
to $R_0$. Note that $\mathcal{S}(R_0)\simeq \mathcal{S}(R_0^\vee)^\vee$ 
if $R_0$ is of type $A,D,E$ or $BC$. 

\begin{prop}
The following reduced irreducible affine root systems form 
a complete set of representatives
of the similarity classes of reduced irreducible affine root systems:
\begin{enumerate}
\item[(1)] $\mathcal{S}(R_0)$ with $R_0$ running through the
similarity classes of reduced irreducible finite crystallographic root systems
(i.e. $R_0$ of type $A,B,\ldots,G$),
\item[(2)] $\mathcal{S}(R_0^\vee)^\vee$ with $R_0$ running
through the similarity classes of reduced irreducible 
finite crystallographic root systems having two 
root lengths
(i.e. $R_0$ of type $B_n$ ($n\geq 2$), $C_n$ ($n\geq 3$), $F_4$ and $G_2$),
\item[(3)] $\mathcal{S}(R_0)$ with $R_0$ a nonreduced irreducible
finite crystallographic root system (i.e. $R_0$ of type $BC_n$ ($n\geq 1$)).
\end{enumerate}
\end{prop}

In view of the above proposition we use the following terminology: 
a reduced irreducible affine root system $R$ is said to be of
{\it untwisted type}
if $R\simeq \mathcal{S}(R_0)$ with 
$R_0$ reduced, of
{\it twisted type} if $R\simeq \mathcal{S}(R_0^\vee)^\vee$ with
$R_0$ reduced, and of
{\it mixed type} if $R\simeq \mathcal{S}(R_0)$ with $R_0$ nonreduced.
Note that a reduced irreducible affine root system $R$ with 
gradient root system of type $A,D$ or $E$ is of untwisted and of twisted
type. 

Suppose that $\Delta_0=(\alpha_1,\ldots,\alpha_n)$ is an ordered basis
of $R_0$. Let $\varphi\in R_0$ (respectively $\theta\in R_0$)
be the associated highest root (respectively the highest short root).
Then
\[
\Delta:=(a_0,a_1,\ldots,a_n)=(c-\varphi,\alpha_1,\ldots,\alpha_n)
\]
is an ordered basis of $\mathcal{S}(R_0)$,
while 
\[
\Delta:=(a_0,a_1,\ldots,a_n)=
(\frac{|\theta|^2}{2}c-\theta,\alpha_1,\ldots,\alpha_n)
\]
is an ordered basis of $\mathcal{S}(R_0^\vee)^\vee$.

\subsection{Nonreduced irreducible affine root systems}
If $R$ is an irreducible affine root system with
ordered basis $\Delta$, then $\Delta$ is also an ordered basis of the 
affine root subsystem $R^{ind}$ of indivisible roots. Furthermore,
$R\simeq R^\prime$ implies $R^{ind}\simeq R^{\prime,ind}$. 
To classify nonreduced
irreducible affine root systems up to similarity, one thus only needs
to understand the possible ways to extend reduced irreducible affine root
systems to nonreduced ones.

Let $R^\prime$ be a reduced irreducible affine root system with affine
Weyl group $W=W(R^\prime)$. Choose an
ordered basis $\Delta=(a_0,a_1,\ldots,a_n)$ of $R^\prime$.
Set
\begin{equation}\label{multiplyablesimple}
S:=
\{a\in \Delta \,\, | \,\, \bigl(\mathbb{Z}R^\prime,a^\vee\bigr)=2\mathbb{Z} \}.
\end{equation}
Let $1\leq m\leq \#S+1$ 
and choose a subset $S_m$ of $S$ of cardinality $m$.
Then
\[
R^{(m)}:=R^\prime\cup\bigcup_{a\in S_m}W(2a)
\]
is an irreducible affine root system with $R^{(m),ind}\simeq R^\prime$.

Considering the possible affine Dynkin diagrams associated
to $(R^\prime,\Delta)$ (see Subsection \ref{Dynkin})
it follows that 
the set $S$ \eqref{multiplyablesimple} is of cardinality $\leq 2$.
It is of cardinality two iff 
$R^\prime\simeq \mathcal{S}(R_0^\vee)^\vee$ 
with $R_0$ of type $\textup{A}_1$ or with
$R_0$ of type $\textup{B}_n$ ($n\geq 2$).
It is of cardinality one iff
$R^\prime\simeq \mathcal{S}(R_0)$ with $R_0$ of type $\textup{B}_n$ ($n\geq 2$) 
or of type $\textup{BC}_n$ ($n\geq 1$).
Hence the similarity class of $R^{(m)}$
does not depend on the choice of subset
$S_m\subseteq S$ of cardinality $m$, and it
does not depend on the choice of ordered basis $\Delta$
of $R^\prime$. The number of $W$-orbits of $R^{(m)}$ equals the number
of $W$-orbits of $R^\prime$ plus $m$. The number of 
similarity classes of irreducible affine root systems $R$ satisfying
$R^{ind}\simeq R^\prime$ is $\#S+1$.

If $R_0$ is of type $\textup{A}_1$ or of type $\textup{B}_n$
($n\geq 2$) we thus have 
a nonreduced irreducible affine root system 
in which two $W$-orbits are added to $\mathcal{S}(R_0^\vee)^\vee$. 
It is labelled as $C^\vee C_n$ by Macdonald \cite{M0}. 
In the rank one case it has four $W$-orbits, otherwise five. 
A detailed description of this affine root system is given in
Subsection \ref{CcheckC}.

Irreducible affine
root subsystems with underlying reduced affine root system $\mathcal{S}(R_0)$ 
having finite
root system $R_0$ of type $\textup{BC}_n$ ($n\geq 1$) or of type
$\textup{B}_n$ ($n\geq 3$) can be naturally viewed as 
affine root subsystems of the affine root system of type $C^\vee C_n$.
This is not the case for the nonreduced extension of the
affine root system $\mathcal{S}(R_0)$ with
$R_0$ of type $\textup{B}_2$. It can actually be better viewed as the rank two
case of the family $\mathcal{S}(R_0)$ with $R_0$ of type $\textup{C}_n$
since, in the corresponding
affine Dynkin diagram, the vertex labelled by the affine simple
root $a_0$ is double bonded with the finite Dynkin diagram of $R_0$. 
The nonreduced extension of $\mathcal{S}(R_0)$ 
with $R_0$ of type $\textup{C}_2$ 
was missing in Macdonald's \cite{M0} classification list. It was
added in \cite[(1.3.17)]{M}.

\subsection{Affine Dynkin diagrams}\label{Dynkin}
In this subsection we list the connected affine Dynkin diagrams 
(cf. \cite[Appendix 1]{M0}) which,
as we have seen, are in one to one correspondence to similarity classes
of irreducible reduced affine root systems. Each similarity class of 
irreducible reduced affine root systems has a representative of the form
$\mathcal{S}(R_0)$ or $\mathcal{S}(R_0^\vee)^\vee$ for a unique 
irreducible finite crystallographic root system $R_0$ up to similarity,
see Subsection \ref{special}.
Recall that $\mathcal{S}(R_0^\vee)^\vee\simeq \mathcal{S}(R_0)$ if $R_0$ is
of type $A,D,E,BC$. We label the connected affine Dynkin diagram
by $\widehat{X}$ with $X$ the type 
of the associated finite root system $R_0$ if 
$X\in\{A,D,E,BC\}$. If the associated finite root system $R_0$ is
of type $X\in\{B,C,F,G\}$ then we label the connected affine Dynkin
diagram by $\widehat{X}^u$ (respectively $\widehat{X}^t$) 
if the associated irreducible reduced affine
root system is $\mathcal{S}(R_0)$ (respectively $\mathcal{S}(R_0^\vee)^\vee$).
Since $A_1\simeq B_1\simeq C_1$ and $B_2\simeq C_2$ there is some
redundancy in the notations, we pick the one
which is most convenient to fit it into an infinite family of
affine Dynkin diagrams. In the terminology of Subsection \ref{special},
the irreducible reduced affine root systems corresponding to
affine Dynkin diagrams labelled by $\widehat{X}$ with
$X\in\{A,D,E\}$ are of untwisted and of
twisted type, labelled by $\widehat{BC}$ 
of mixed type, labelled by $\widehat{X}^u$ of untwisted
type and labelled by $\widehat{X}^t$ of twisted type.
In \cite[Appendix 1]{M0} the affine Dynkin diagrams labelled $\widehat{B}_n^t$
and $\widehat{C}_n^t$ are called of type $C_n^\vee$ and $B_n^\vee$ respectively.
The remaining relations with the notations and terminologies 
in \cite[Appendix 1]{M0} are self-explanatory. 

We specify in each affine Dynkin diagram a particular vertex
(the grey vertex) which is labelled by the unique affine simple root
$a_0$ in the particular choice of ordered basis
$\Delta$ of $\mathcal{S}(R_0)$ or $\mathcal{S}(R_0^\vee)^\vee$ as specified
in Subsection \ref{special}.  

In Kac's notations (see Tables Aff 1--3 in \cite[\S 4.8]{K})
the affine Dynkin diagrams 
are labelled differently: our label $\widehat{X}$ corresponds to
$X^{(1)}$ if $X\in\{A,D,E\}$ and
$\widehat{BC}_n$ corresponds to $A_{2n}^{(2)}$ ($n\geq 1$).
Our label $\widehat{X}^u$ 
corresponds to $X^{(1)}$ if $X\in \{B,C,F,G\}$. 
Finally, our label $\widehat{B}_n^{t}$ corresponds to
$D_{n+1}^{(2)}$ ($n\geq 2$), 
$\widehat{C}_n^{t}$ corresponds to $A_{2n-1}^{(2)}$ ($n\geq 3$),
$\widehat{F}_4^{t}$ to $E_6^{(2)}$ and
$\widehat{G}_2^{t}$ to $D_4^{(3)}$.
\newpage
\begin{figure}[http]
\begin{center}
\includegraphics[width=3cm,height=1.5cm]{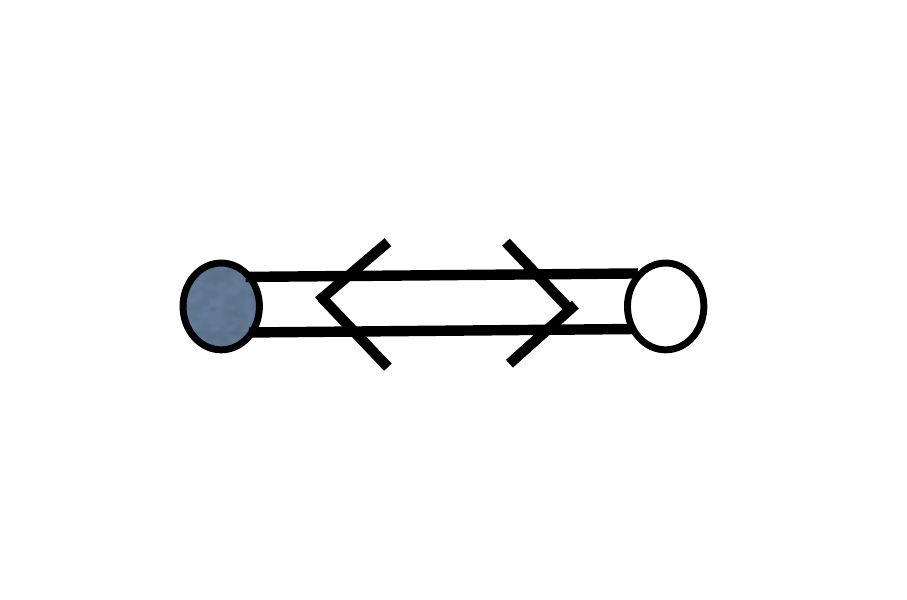}
\includegraphics[width=5cm,height=2.5cm]{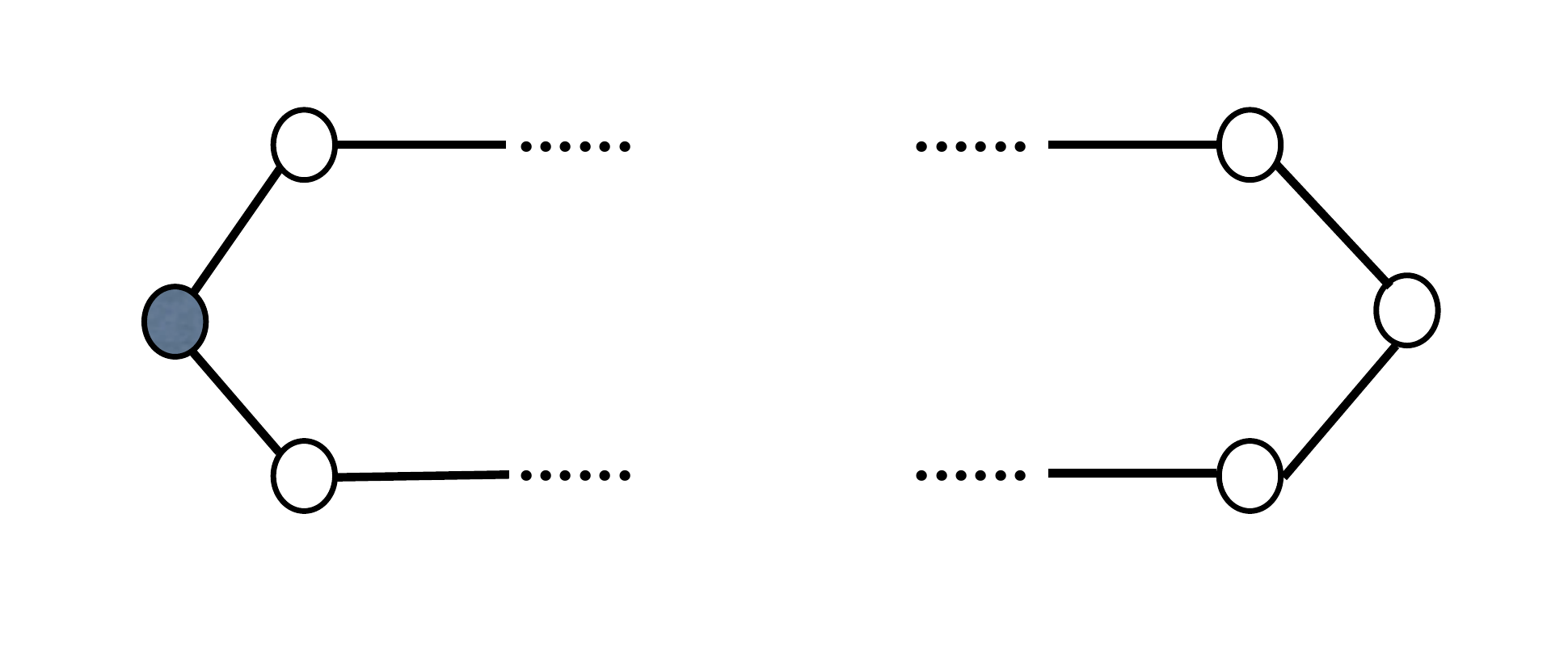}
\caption{$\widehat{A}_1$ and $\widehat{A}_n$ ($n\geq 2$)}
\end{center}
\end{figure}
\begin{figure}[http]
\begin{center}
\includegraphics[width=5cm,height=3cm]{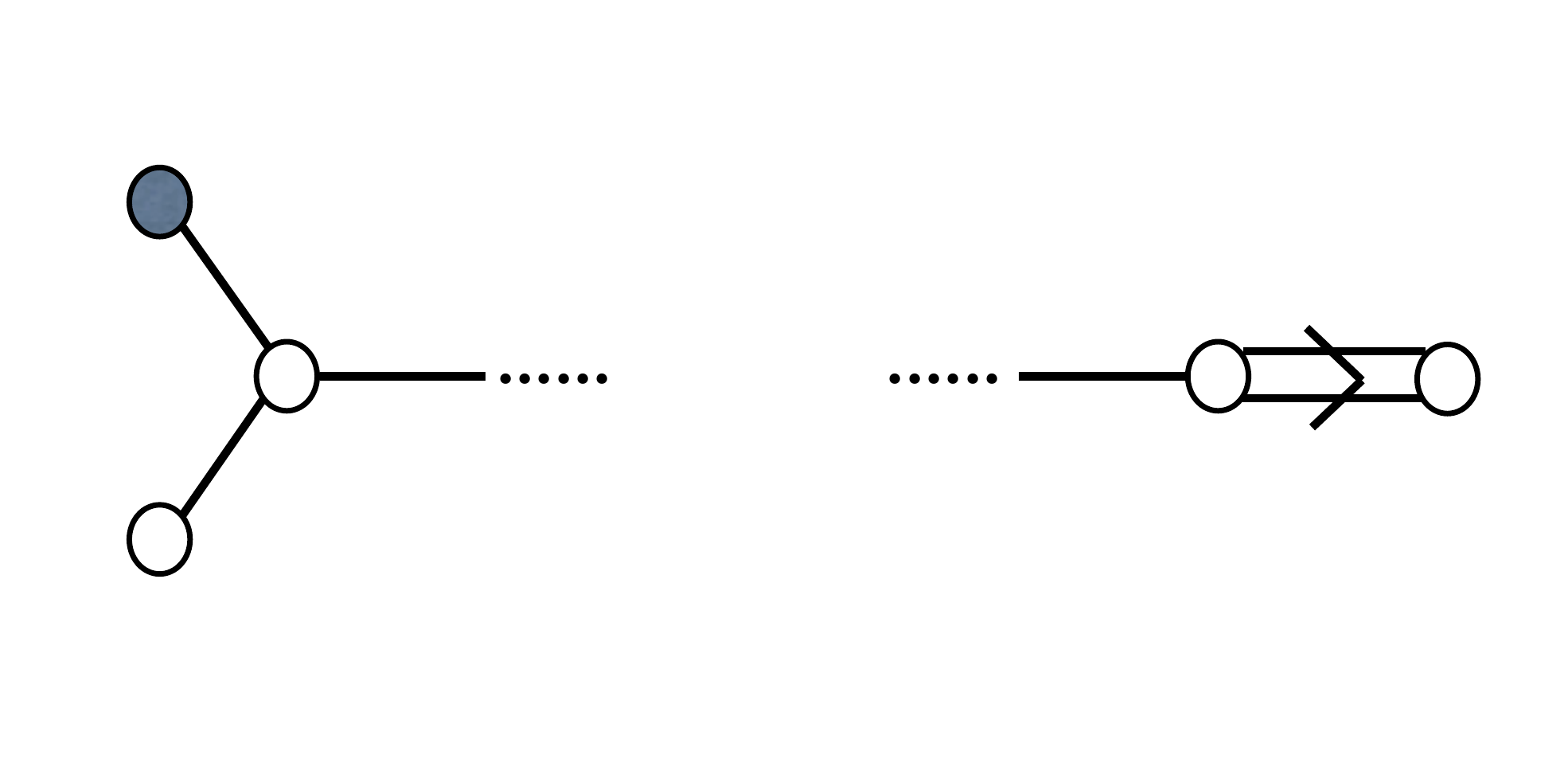}
\caption{$\widehat{B}_n^u$ ($n\geq 3$)}
\end{center}
\end{figure}
\begin{figure}[http]
\begin{center}
\includegraphics[width=5cm,height=2cm]{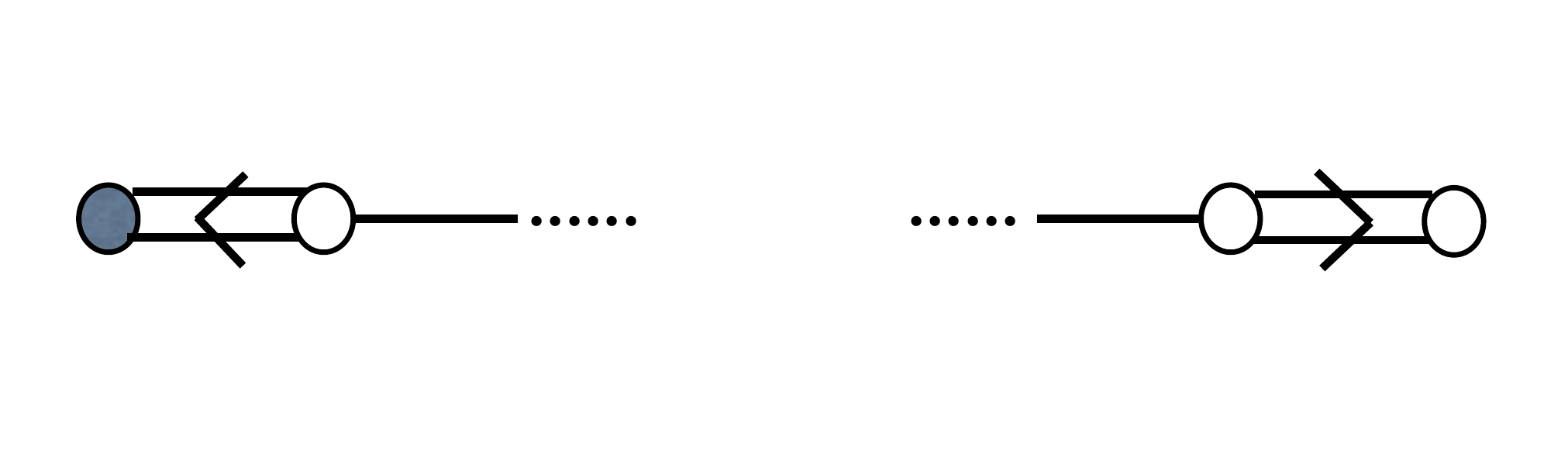}
\caption{$\widehat{B}_n^{t}$ ($n\geq 2$)}
\end{center}
\end{figure}
\begin{figure}[http]
\begin{center}
\includegraphics[width=2cm,height=1cm]{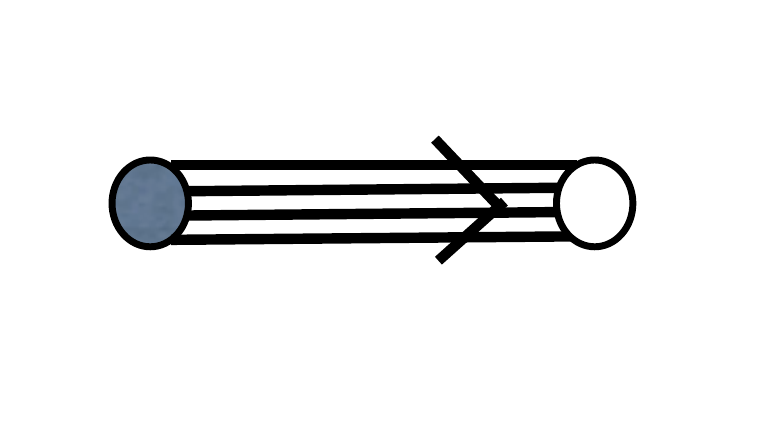}
\includegraphics[width=4cm,height=1.5cm]{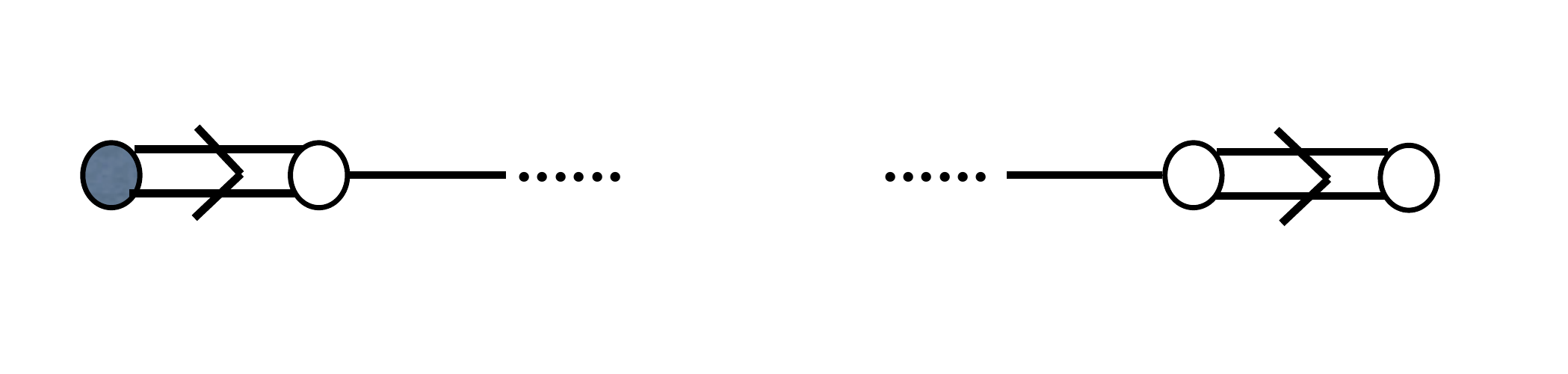}
\caption{$\widehat{BC}_{1}$ and $\widehat{BC}_{n}$ ($n\geq 2$)}
\end{center}
\end{figure}
\begin{figure}[http]
\begin{center}
\includegraphics[width=5cm,height=2cm]{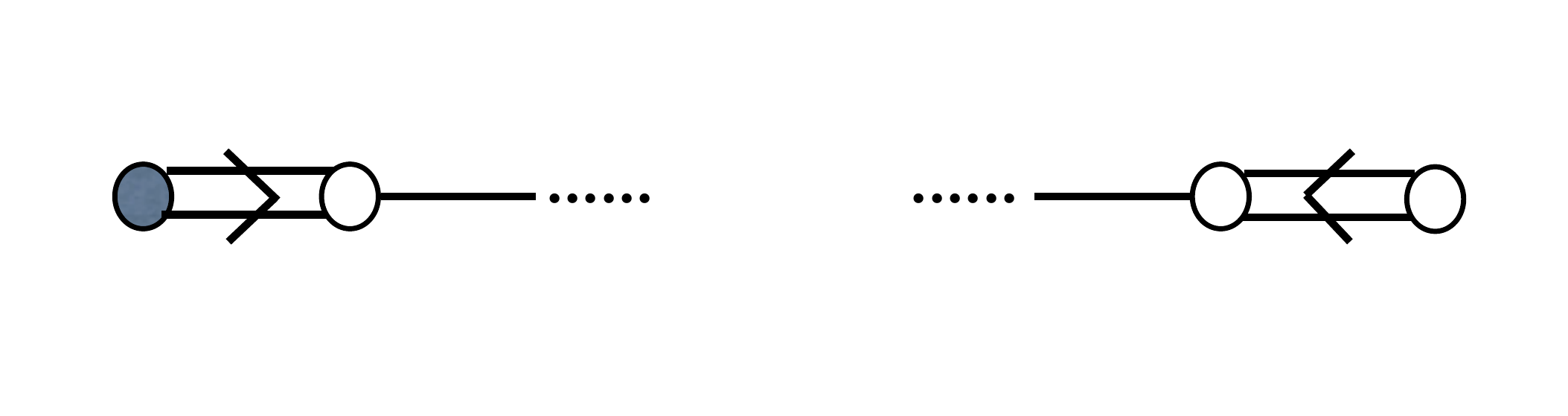}
\caption{$\widehat{C}_n^u$ ($n\geq 2$)}
\end{center}
\end{figure}
\begin{figure}[http]
\begin{center}
\includegraphics[width=5cm,height=2.5cm]{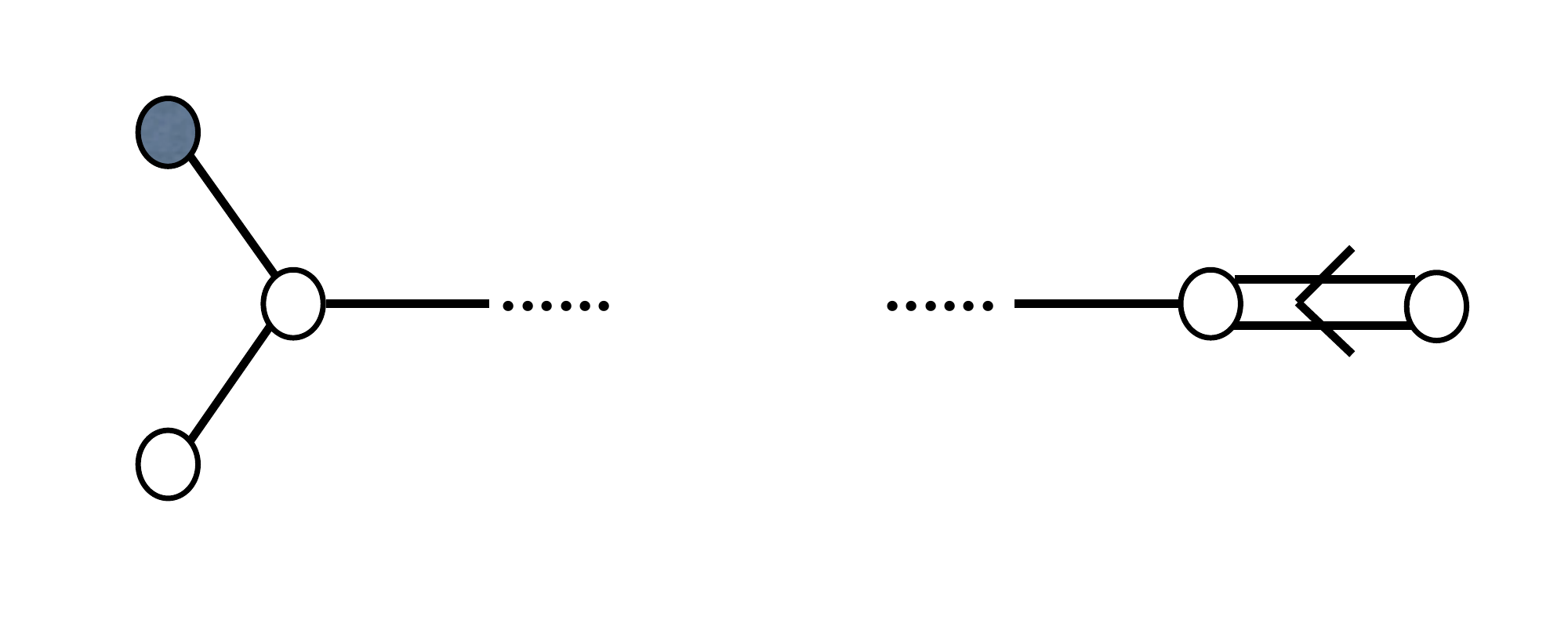}
\caption{$\widehat{C}_n^{t}$ ($n\geq 3$)}
\end{center}
\end{figure}
\begin{figure}[http]
\begin{center}
\includegraphics[width=5cm,height=2cm]{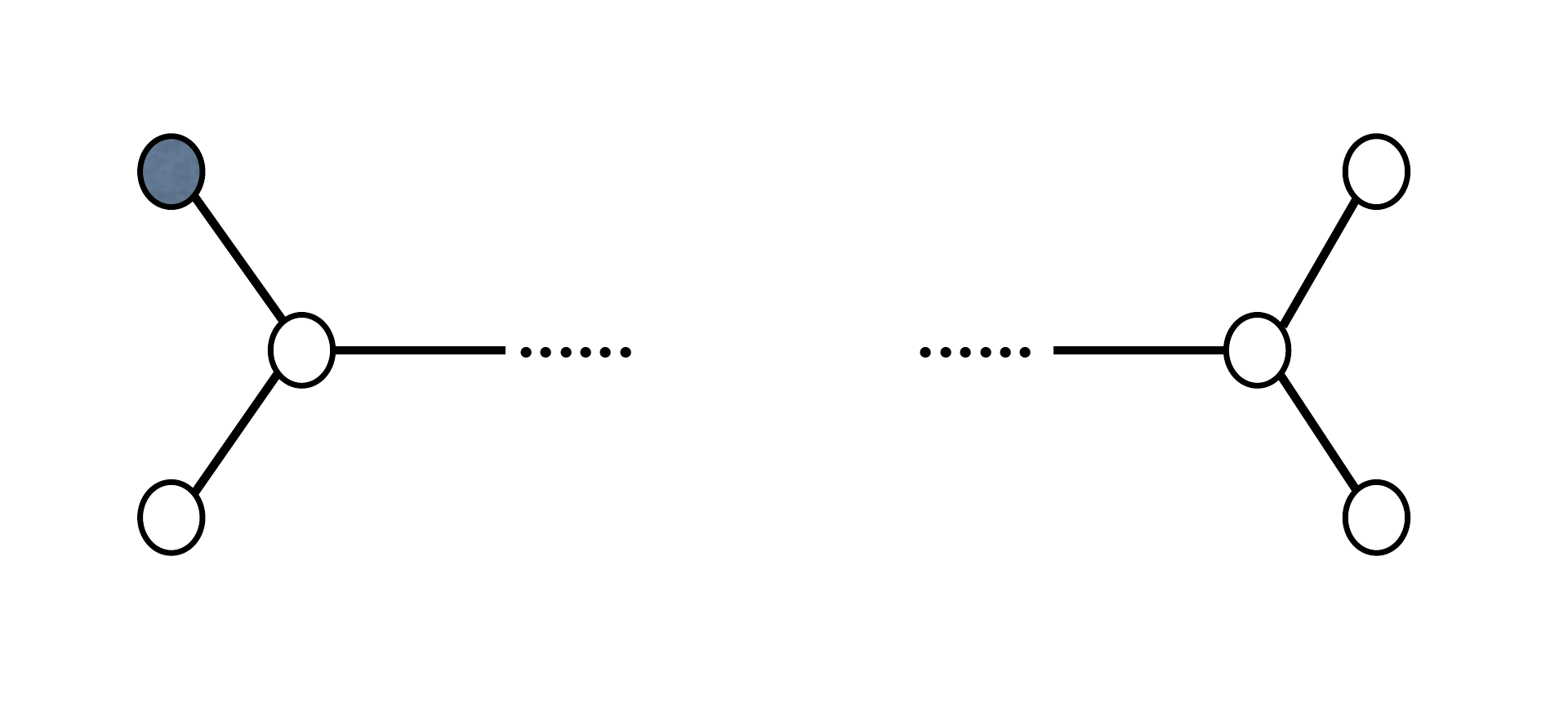}
\caption{$\widehat{D}_n$ ($n\geq 4$)}
\end{center}
\end{figure}
\begin{figure}[http]
\begin{center}
\includegraphics[width=3cm,height=2cm]{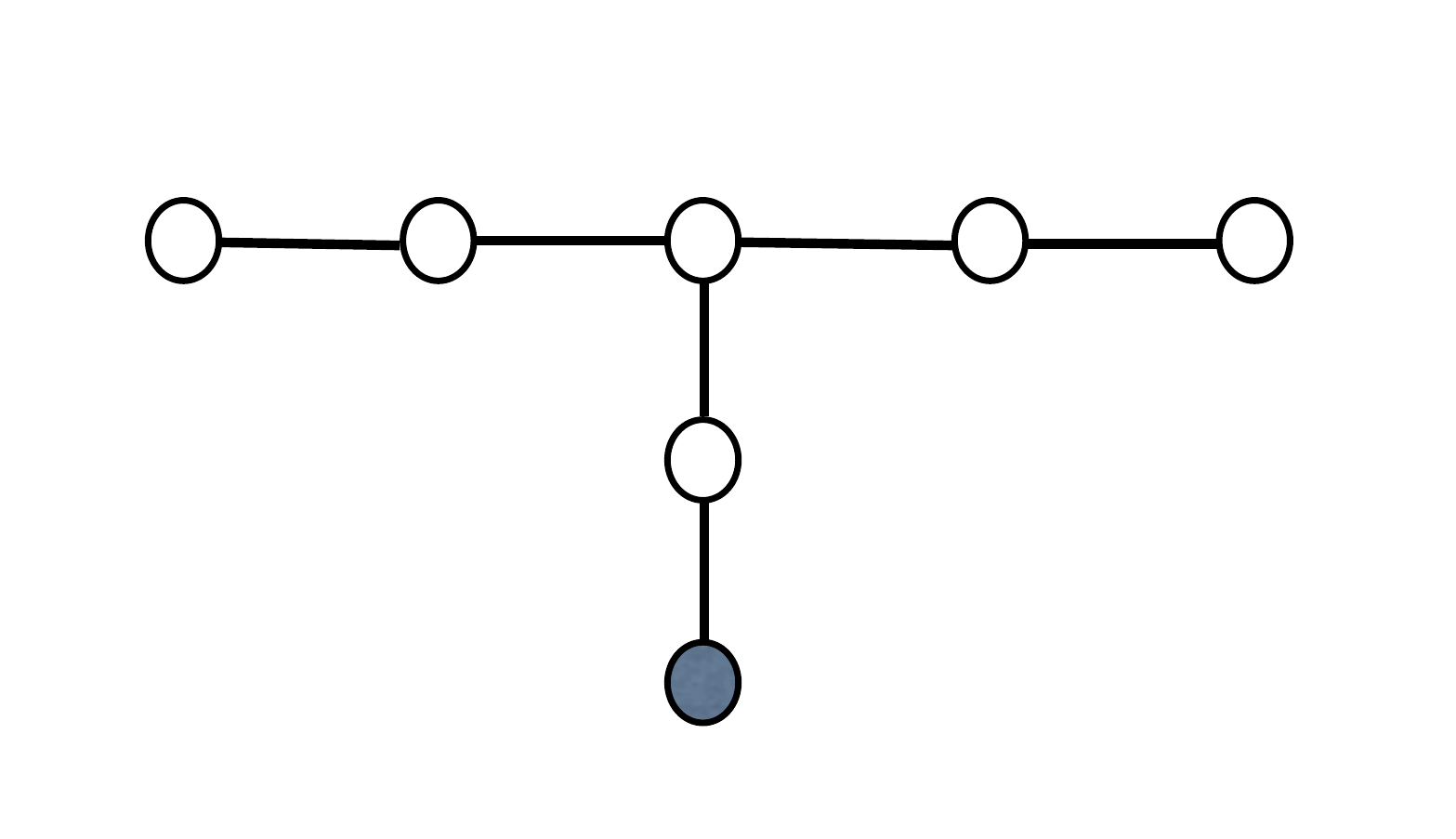}
\includegraphics[width=3cm,height=2cm]{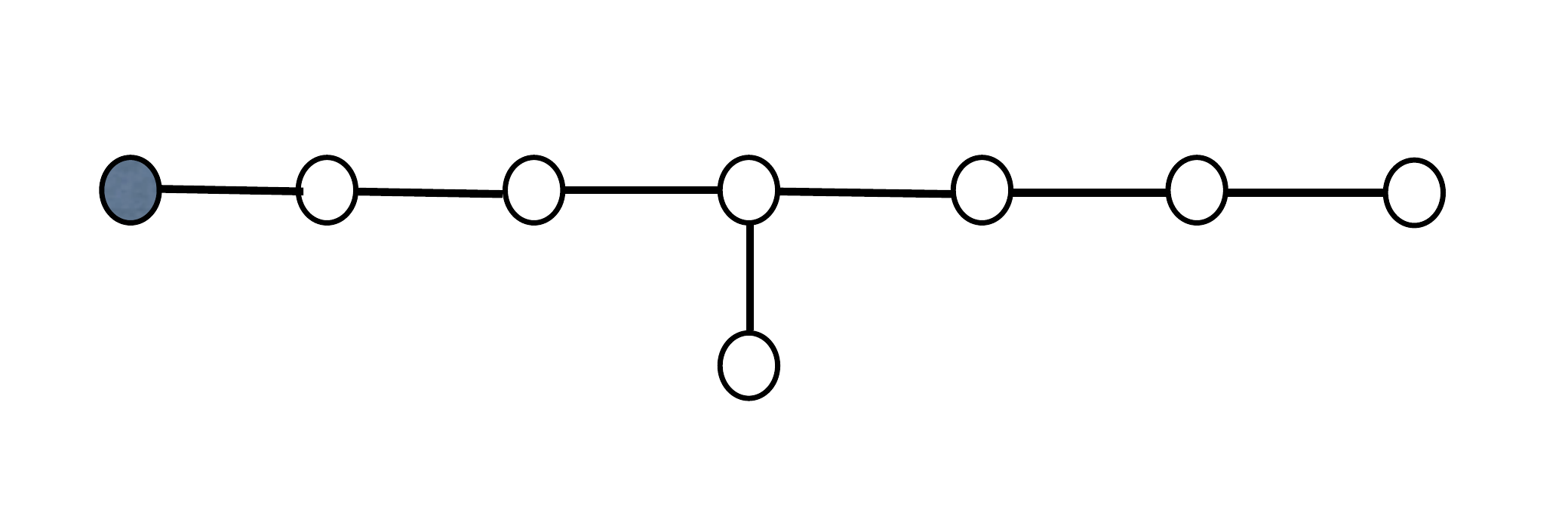}
\includegraphics[width=3cm,height=2cm]{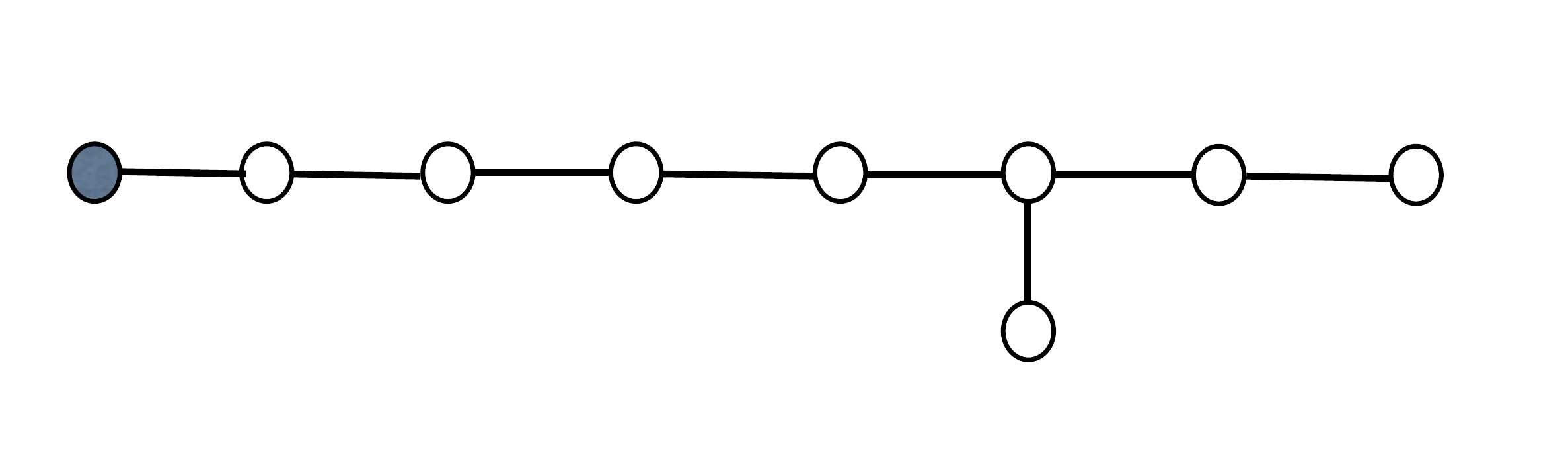}
\caption{$\widehat{E}_6$, $\widehat{E}_7$ and
$\widehat{E}_8$}
\end{center}
\end{figure}
\begin{figure}[http]
\begin{center}
\includegraphics[width=2.5cm,height=1.5cm]{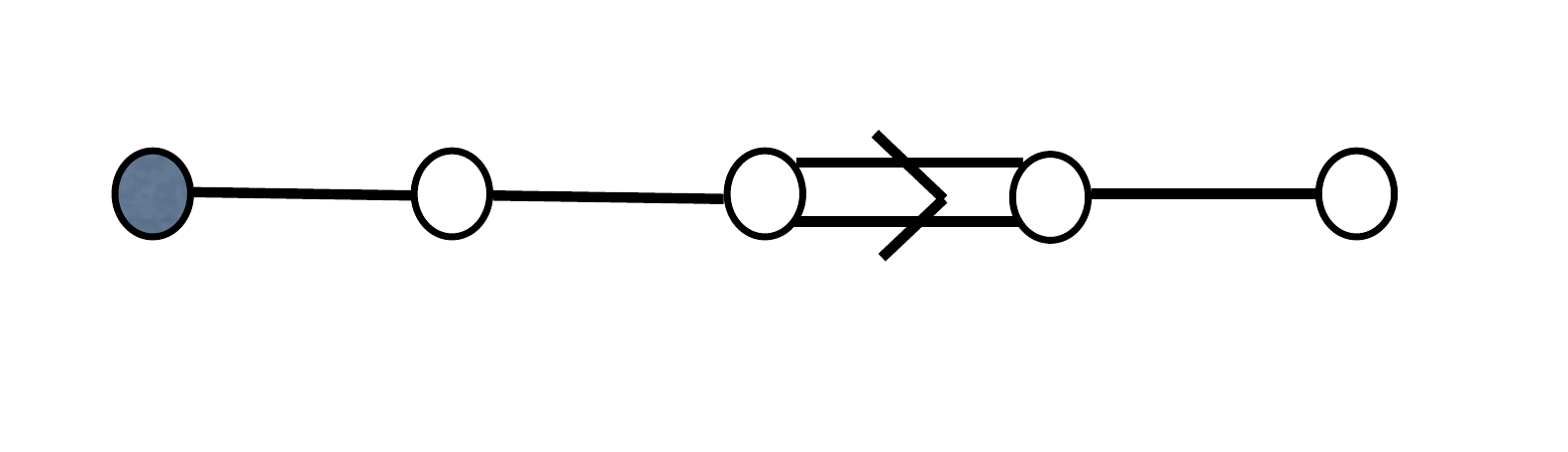}
\includegraphics[width=2.5cm,height=1.5cm]{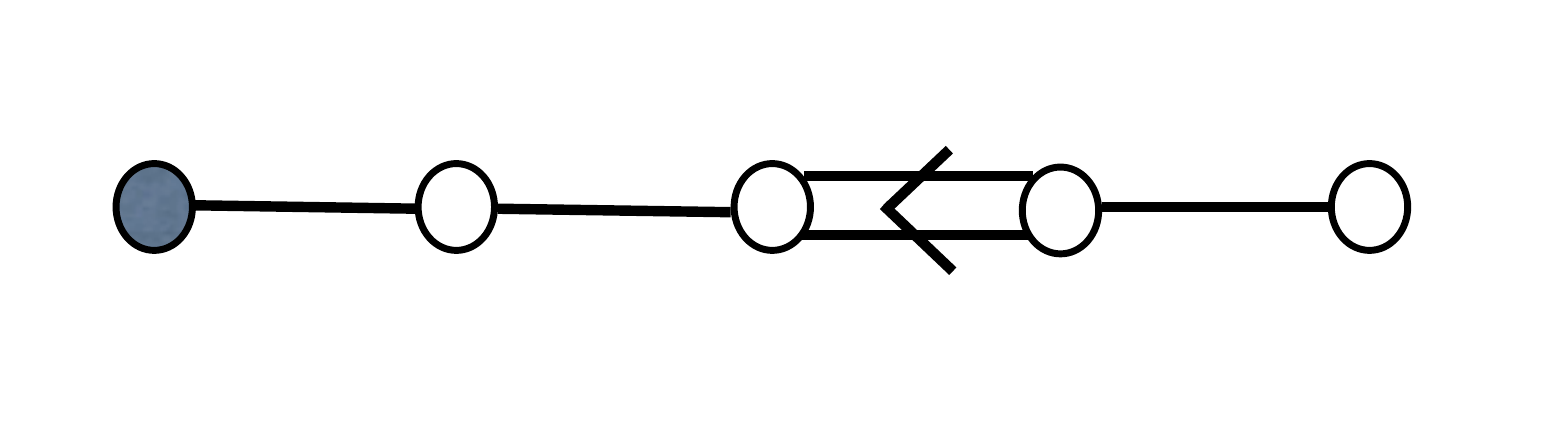}
\caption{$\widehat{F}_4^u$ and $\widehat{F}_4^{t}$}
\end{center}
\end{figure}
\begin{figure}[http]
\begin{center}
\includegraphics[width=2.5cm,height=1.5cm]{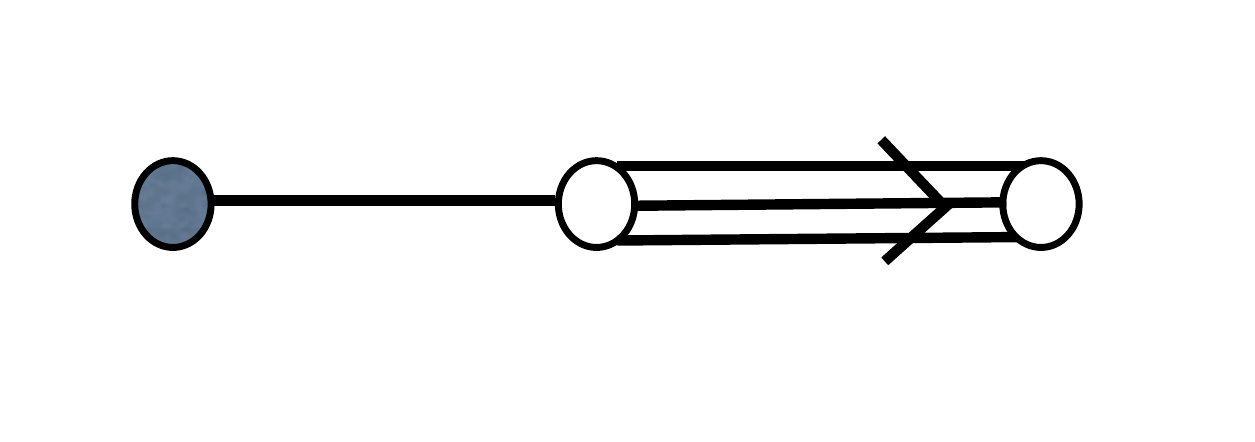}
\includegraphics[width=2.5cm,height=1.5cm]{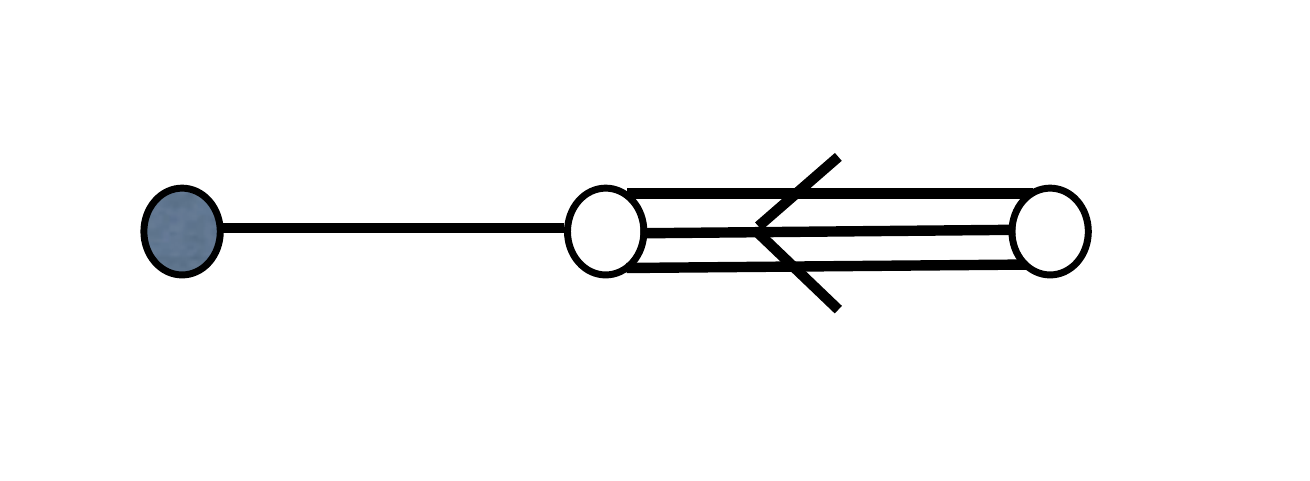}
\caption{$\widehat{G}_2^u$ and $\widehat{G}_2^{t}$}
\end{center}
\end{figure}
\newpage
%


\end{document}